\documentclass[reqno,11pt]{amsart}
\usepackage{geometry}
\usepackage{mathtools,amssymb,amsthm,mathrsfs,color,lineno,paralist,graphicx,float}
\usepackage[colorlinks,
linkcolor=blue,
anchorcolor=green,
citecolor=blue,
]{hyperref}

\usepackage[T1]{fontenc}
\usepackage[utf8]{inputenc}

\usepackage{calc}
\definecolor{bleu1}{RGB}{0,57,128}
\def\bleu1{\color{bleu1}}

\usepackage{etoolbox}
\patchcmd{\section}{\normalfont}{\normalfont \bleu1}{}{}
\patchcmd{\subsection}{\normalfont}{\normalfont \bleu1}{}{}
\patchcmd{\subsubsection}{\normalfont}{\normalfont \bleu1}{}{}
\renewcommand{\proofname}{\it \bleu1 Proof}

\let\newpf\proof \let\proof\relax 
\newenvironment{pf}{\newpf[\proofname]}{\qed\endtrivlist}

\newcommand{\ba}{\overline{A}}

\def\be{\begin{equation}}
\def\ee{\end{equation}}

\def\ba{{\begin{align}}}
\def\ea{{\end{align}}}

\def\bm{\begin{matrix}}
\def\em{\end{matrix}}

\def\SL{{\mathrm{SL}}}

\def\GL{{\mathrm{GL}}}
\def\diag{{\mathrm{diag}}}

\def\0{{\mathbf 0}}

\def\ii{{\textrm i}}

\def\m{{\tilde m}}

\newtheorem{Theorem}{Theorem}[section]
\newtheorem*{Theorem*}{Theorem 1.1}
\newtheorem{Lemma}{Lemma}[section]
\newtheorem{Proposition}{Proposition}[section]

\newtheorem{Remark}{Remark}[section]

\newtheorem{Definition}{Definition}[section]
\newtheorem{Claim}{Claim}

\numberwithin{equation}{section}

\theoremstyle{definition}

\newtheorem{problem}{Problem}[section]
\newtheorem{conjecture}{Conjecture}[section]

\newcommand{\const}{\mathrm{const}}

\renewcommand{\ii}{{\mathrm{i}}}

\newcommand{\A}{{\mathbb A}}
\newcommand{\C}{{\mathbb C}}

\newcommand{\N}{{\mathbb N}}
\newcommand{\Q}{{\mathbb Q}}
\newcommand{\R}{{\mathbb R}}
\newcommand{\T}{{\mathbb T}}

\newcommand{\Z}{{\mathbb Z}}

\newcommand{\la}{\langle}
\newcommand{\ra}{\rangle}

\def\B0{{\bold{0}}}

\catcode`\@=12

\def\Empty{}
\newcommand\oplabel[1]{
  \def\OpArg{#1} \ifx \OpArg\Empty {} \else
    \label{#1}
  \fi}

\newcommand{\comm}[1]{}
\newcommand{\comment}[1]{}

\begin{document}

\title[Absolutely continuity of IDS]{ Absolute continuity of the integrated density of states in the localized regime}

\author{ Jing Wang}

\address{
School of Mathematics and Statistics, Nanjing University of Science and Technology, Nanjing 210094, China
}

\email{jing.wang@njust.edu.cn}

\author{Xu Xu} \address{
Department of Mathematics, Nanjing University, Nanjing 210093, China
}

\email{377056722@qq.com}

\author {Jiangong You}
\address{
Chern Institute of Mathematics and LPMC, Nankai University, Tianjin 300071, China} \email{jyou@nankai.edu.cn}

\author {Qi Zhou}
 \address{
Chern Institute of Mathematics and LPMC, Nankai University, Tianjin 300071, China}

\email{qizhou@nankai.edu.cn}

\setcounter{tocdepth}{1}

\begin{abstract}
We establish the absolute continuity of the integrated density of states (IDS) for quasi-periodic Schr\"odinger operators with a large trigonometric potential and Diophantine frequency. This partially solves Eliasson's open problem in 2002. Furthermore, this result can be extended to a class of quasi-periodic long-range operators on $\ell^2(\Z^d)$. Our proof is based on stratified quantitative almost reducibility results of dual cocycles. Specifically, we prove that a generic analytic one-parameter family of cocycles, sufficiently close to constant coefficients, is reducible except for a zero Hausdorff dimension set of parameters. This result affirms Eliasson's conjecture in 2017.
\end{abstract}

\maketitle

%

\section{Introduction}

\subsection{Quasi-periodic Schr\"odinger operator on $\ell^2(\Z)$}
Consider the one-dimensional quasi-periodic Schr\"odinger operator
\begin{equation}\label{op-11}
(H_{\varepsilon^{-1}V, \alpha, x}u)_n=u_{n+1}+u_{n-1}+\varepsilon^{-1}V(x+n\alpha)u_n,\ \ n\in\Z,
\end{equation}
on $\ell^2(\Z)$, where  the phase $x\in\T$,   the frequency $\alpha \in \R\backslash \Q$ is irrational and the potential $V\in C^\omega(\T,\R)$.
 The Almost Mathieu operator (AMO), given by: \[(H_{2\varepsilon^{-1}\cos,\alpha,x}u)_n=u_{n+1}+u_{n-1}+2\varepsilon^{-1}\cos2\pi(x+n\alpha)u_n,\ \ n\in\Z,\] is the most well-known example of the class of operators described above, where $V(x)=2\cos(2\pi x)$. Peierls \cite{Pe} originally introduced it as a model for an electron on a 2D lattice, subject to a homogeneous magnetic field \cite{Ha,R}. These operators have drawn significant attention, not just because of their importance in physics \cite{Ai,AOS,OA}, but also as fascinating mathematical objects \cite{A2,AJ,AYZ,AYZ2,J,JLiu}.

Our focus is on investigating the regularity of the integrated density of states (IDS). The IDS can be regarded as the average spectral measures of an ergodic family of self-adjoint operators  $\{L_x\}_{x\in X}$ over $x$:
\[\mathcal{N}(E)=\int_{x\in X}\mu_x(-\infty,E]dx,\]
where $\mu_x$  is the associated spectral measure of $L_x$. Understanding the regularity of IDS, including absolute continuity \cite{A1,AD2,GJZ,GS} and H\"older continuity \cite{Amor, aj1, ALSZ,GS0,GS}, is a significant area of research in the spectral theory of quasi-periodic Schr\"odinger operators. It is closely linked to other important topics, such as homogeneous spectrum \cite{DGL1,dgsv,LYZZ}, Parreau–Widom condition \cite{ALSZ}, and Deift's conjecture \cite{BDGL, LYZZ}.

This paper will specifically address the absolute continuity of the IDS. As it is common knowledge, the pure absolutely continuous spectrum  \cite{A2,A1,aj1,DS,Eli92} leads directly to the absolute continuity of the IDS. In the regime where the Lyapunov exponent is zero, Kotani \cite{K} proved a more deeper result: the absolute continuity of the IDS is equivalent to the absolute continuity of the spectral measures for almost every phase. Nevertheless, in the regime where the Lyapunov exponent is positive, the spectral measure is typically singular \cite{B0,BG,Eli97,GS,J}. For instance, Eliasson \cite{Eli97} established that $H_{\varepsilon^{-1}V, \alpha, x}$ has  pure point spectrum for almost all $x$ when $\varepsilon$ is small, and $\alpha$ is Diophantine. Recall that $\alpha\in\T^d$ is Diophantine (denoted as $DC_d:=\cup_{\gamma, \tau}DC_d(\gamma,\tau)$), if there exists $\gamma>0$, $\tau>d-1$, and
\[\inf_{j\in\Z}|\la k,\alpha\ra-j|\geq \frac{\gamma}{|k|^\tau}, \ \ \forall k\in \Z^d\backslash\{0\}.\] 
In such cases, determining whether IDS is absolutely continuous presents a highly challenging obstacle, as Eliasson identified in his famous open problem:

\begin{problem}[\cite{MY}, Problem 4.2.1, Eliasson]\label{problem-Eliasson}
What are the properties of the map $E\mapsto \mathcal{N}(E)$ when $|\varepsilon|$ is small? Is it singular continuous? Absolutely continuous?
\end{problem}

In this paper, we  answer his question for trigonometric polynomial potentials:
\begin{Theorem}\label{main}
Suppose that $\alpha\in DC_1$ and $V(x)$ is a trigonometric polynomial on $\T$. Then there exists $\varepsilon_0=\varepsilon_0(\alpha, V)>0$ such that  $\mathcal{N}(E)$ of (\ref{op-11}) is absolutely continuous if $| \varepsilon|< {\varepsilon_0}$.
\end{Theorem}

Let us take a moment to briefly revisit the previous research concerning the absolute continuity of IDS in the localized regime. There are two categories of results. The first category is for a \textit{fixed} analytic potential. In the positive Lyapunov exponent region, Bourgain and Goldstein \cite{BG} proved that $H_{\varepsilon^{-1}V, \alpha, x}$ has Anderson localization for almost all Diophantine $\alpha$ at any $x$. Suppose that  $V (x)$ is a small perturbation of a trigonometric polynomial; Goldstein and Schlag \cite{GS} further proved that $\mathcal{N}(E)$ of (\ref{op-11}) is absolutely continuous for almost every $\alpha$.

A separate category pertains to a \textit{fixed} $\alpha\in\mathbb{R}\backslash\mathbb{Q}$. For AMO, Jitomirskaya \cite{J} proved that if $|\varepsilon|<1$, then $H_{2\varepsilon^{-1}\cos,\alpha,x}$ displays Anderson localization for almost all $x$ when $\alpha\in\R\backslash\Q$.
Moreover, Avila and Damanik \cite{AD2} demonstrated that $\mathcal{N}(E)$ exhibits absolute continuity iff $\varepsilon\neq \pm 1$. Ge, Jitomirskaya, and Zhao \cite{GJZ} established that when a strong Diophantine $\alpha$ is fixed, the IDS of analytically perturbed non-critical AMO is absolutely continuous provided the perturbation is sufficiently small in a non-perturbative sense.	

In summary, for a fixed Diophantine frequency, previous studies are all restricted to cosine or cosine-like potentials within the region of positive Lyapunov exponents. However, for general analytic potentials, the results are measure-theoretic in $\alpha$, which depend on the potentials  in a very implicit way.   By contrast, our study works for {\it any fixed} Diophantine frequency, and {\it fixed} trigonometric polynomial.

\subsection{Quasiperiodic long-range operator on $\ell^2(\Z^d)$}
It is worth noting that our result is applicable to quasi-periodic operators located on $\ell^2(\Z^d)$, not just one-dimensional Schrödinger operators on $\ell^2(\Z)$. More specifically, we study the following quasi-periodic operator:
  \begin{equation}\label{op-12-1}
 	(L_{ \varepsilon^{-1}V,\alpha,x}^{{ W}} u)_n=\sum_{k\in\Z^d}\hat W_ku_{n-k}+\varepsilon^{-1}V(x+\la n,\alpha\ra)u_n, \ \ n\in\Z^d,
 \end{equation}
where  $V\in C^\omega(\T,\R)$, and $W(\theta)=\sum_{k\in\Z^d}\hat W_k e^{2\pi \ii\la k,\theta\ra}$ is  real analytic   on $\T^d=\R^d/\Z^d$, and  $(1,\alpha)\in\R^{1+d}$ is rationally independent. There has been extensive analysis of the localization problem of \eqref{op-12-1} inspired by the pioneer works of Fröhlich, Spencer, Wittwer \cite{FSW} and Sinai \cite{Sin}. Several studies have focused on $C^2$-$\cos$ type potentials   \cite{CSZ1,ChD}  or just $V(x)=2\cos(2\pi x)$ \cite{BJ, CSZ,GY, GYZo,JK}.  

Though Bourgain \cite{B0} provided evidence, little advancement was observed if $V$ is a general analytic potential. For a {\it fixed} phase $x\in\T$, Bourgain \cite{B0} proved that the operator (\ref{op-12-1}) has Anderson localization for a positive measure set of $\alpha$ if $|\varepsilon|$ is sufficiently small. While the conceptual belief is that the answer is affirmative, it is still unresolved whether the operator (\ref{op-12-1}) undergoes Anderson localization for almost all $x\in\T$, when $\alpha\in DC_d$ is  {\it fixed}.

At this point, it is reasonable to inquire whether $ \mathcal{N}(E)$ exhibits absolute continuity in the localization regime for (\ref{op-12-1}), similar to Problem \ref{problem-Eliasson}. In this paper, we address this inquiry in circumstances where $V$ is a trigonometric polynomial.

 \begin{Theorem}\label{Thm1-1}
 	Fix $\alpha\in DC_{d}$. Suppose that $W(\theta)$ is analytic on $\T^d$, and $V(x)$ is a trigonometric polynomial on $\T^1$. Then there exists $\varepsilon_0=\varepsilon_0(\alpha, d, V, W)>0$ such that   $ \mathcal{N}(E)$ of \eqref{op-12-1}  is  absolutely continuous, provided $| \varepsilon|< {\varepsilon_0}$.
 \end{Theorem}

 \begin{Remark}
 	To the best knowledge of the authors, this gives the first result of the absolute continuity of IDS for quasi-periodic long-range operator on $\ell^2(\Z^d)$ with $d>1$.
 \end{Remark}

\smallskip
We will now briefly outline the main ideas of  the proof. Two methods have previously been established for examining the IDS's absolute continuity in the localized regime. The first one, developed by Goldstein-Schlag \cite{GS}, utilizes the large deviation theorem and avalanche principle for the operator's truncated determinants. However, due to the problem with "double resonances", an unknown zero measure set of frequencies has to be removed \cite{B0,BG,GS}.
Ge-Jitomirskaya-Zhao \cite{GJZ} presented another method, which is supported by the fundamental observation made by Sodin-Yuditskii \cite{sy1, sy2} that the spectral measure is absolutely continuous if the normal boundary's real part of the Borel transform of the measure is integrable, and its topological support is homogeneous. Nevertheless, the aforementioned method is only applicable to Schr\"odinger cocycles with an acceleration of 1. This paper introduces a novel method, which is based on reducibility theory and Aubry duality. As demonstrated, proving the absolutely continuous nature of IDS requires proving the reducibility of its dual quasi-periodic cocycles beyond a set of energies with zero Hausdorff dimension. We also note that all previous research \cite{AD2,GJZ,GS} pertains to one-frequency and is non-perturbative in nature (the smallness of $\varepsilon$ does not depend on $\alpha$), whereas our research includes the multiple-frequency setting. Nonetheless, non-perturbative results should not be anticipated in this case \cite{B2002}.

 \subsection{Stratified quantitative almost reducibility}

Let's explain the concept of reducibility and our precise results. We denote $\mathrm{GL}(m,\C)$ as the set of all $m\times m$ invertible matrices.  For any rationally independent $\alpha\in\R^d$ and $A\in C^\omega(\T^d, \mathrm{GL}(m,\C))$,  the \textit{analytic quasi-periodic $\GL(m,\C)$-cocycle} $(\alpha, A): \T^d\times\C^m\rightarrow \T^d\times \C^m$ is defined as the skew product
  \[(x, u)\mapsto (x+\alpha, A(x)u).\]
 A $\GL(m,\C)$-cocycle $(\alpha, A)$ is said to be \textit{$C^\omega$ reducible}, if there exists $B\in C^\omega(\T^d, \GL(m,\C))$, $\tilde{A}\in \GL(m,\C)$, such that  $$B^{-1}(x+\alpha)A(x)B(x)=\tilde{A}.$$
   From now on, we will always assume $\alpha$ to be Diophantine.
   
The earliest result of local reducibility was due to Dinaburg and Sinai ~\cite{DS}. They demonstrated that if the potential satisfies the assumptions of being analytic and small, then for a positive measure set of energies $E$, the Schr\"odinger cocycle is reducible. Eliasson~\cite{Eli92}, under the same assumptions as Dinaburg and Sinai, showed that the Schr\"odinger cocycle is reducible for a full measure set of energies $E$. The proof is based on what is known "resonance-cancellation'' technique, which originated from Moser and P\"oschel's~\cite{MP} research.  Following this, Eliasson made a well-known conjecture \cite{Eli98,Eli11}: that a generic one-parameter family of analytic cocycles, which are close enough to constant coefficients, is reducible for \textit{a.e.} of the parameters. Krikorian proved this conjecture in linear systems with coefficients in $so(3)$ and Lie algebra of compact semi-simple Lie group in general \cite{K3,K2}, and it was further verified in linear systems, taking values  in $gl(m,\C)$ \cite{HeY}.

Eliasson \cite{Eli17} recently presented a new conjecture during the conference in memory of Jean-Christophe Yoccoz:

 \begin{conjecture}[\cite{Eli17}]\label{conj-eli-3}
The Hausdorff dimension of the complementary set of parameters for reducible cocycles is zero.
\end{conjecture}
Our paper aims to prove this conjecture. We will firstly  provide necessary notations to clarify our approach. Let $\Lambda\subseteq \R$ be a bounded interval and $A: \Lambda\rightarrow \GL(m,\C)$ is analytic in $\lambda\in\Lambda$.
Denote $W_{\delta}(\Lambda):=\{z\in\C:\mathrm{dist}(z,\Lambda)<\delta\}$. For any analytic $*$-valued function $f: \Lambda\rightarrow *$, where $*$ can be $\R, \C, \GL(m,\R), \GL(m,\C)$, we set
\[|f|_\delta:=\sup_{z\in W_\delta(\Lambda)}\|f(z)\|,\]
where $\|\cdot\|$ denotes the absolute value or matrix norm correspondingly,
and we say $f\in C_\delta^\omega(\Lambda, *)$ if $|f|_\delta<\infty$. For an integrable $*$-valued function $f: \T^d\times\Lambda\rightarrow *$, let
\[|f|_{h,\delta}:=\sup_{z\in W_{\delta}(\Lambda)}|f(z)|_h,\]
where $h,\delta>0$, and for any integrable function $g :\T^d\rightarrow *$  we set
\[|g|_h:=\sum_{k\in\Z^d}\|\hat g(k)\|e^{2\pi |k|h},\]
with $\hat g(k)=\int_{\T^d}f(\phi)e^{-2\pi \ii\la k,\phi\ra}d\phi$, and $\|\cdot\|$ denoting the absolute value or matrix norm correspondingly. We say $f\in C_{h,\delta}^\omega(\T^d\times\Lambda, *)$ if $|f|_{h,\delta}<\infty$.

Denote $\Sigma(A(\lambda))=:\Sigma(\lambda)$ as the set of eigenvalues of $A(\lambda)\in \GL(m,\C)$, and for any $u\in\T$, let
\begin{equation}\label{trans-func-2}
g(\lambda,u)=\prod_{\sigma_i,\sigma_j\in \Sigma(\lambda), \atop{i\neq j}}(\sigma_i-e^{2\pi \ii u}\sigma_j).
\end{equation}
\begin{Definition}
We say that $A(\lambda)$ satisfies the $non$-$degeneracy$ condition on an interval $\Lambda$, if there exists $r\in\N^{+}$, $c>0$ such that for $\forall u\in\T$, the following inequality holds for all $\lambda\in\Lambda$,
\begin{equation}\label{def-non-degeneracy}
\max_{0\leq l\leq r}|\frac{\partial^{l}g(\lambda,u)}{\partial\lambda^{l}}|\geq c,
\end{equation}
where $g$ is defined as in (\ref{trans-func-2}).
\end{Definition}


We can now declare our reducibility result.

 \begin{Theorem}\label{Thm2}
Let $h>0$, $\alpha\in\mathrm{DC}_{d}$, and $\Lambda\subset\R$ be an interval. Suppose that $A \in C^{\omega}(\Lambda,\mathrm{GL}(m,\C))$ is non-degenerate on $\Lambda$ in the sense of (\ref{def-non-degeneracy}) with some $r\in\N^+, c>0$. Then there exists $\delta>0$,  $\varepsilon_{0}=\varepsilon_0(\alpha,d,m,\delta,h,r,c,|A|_\delta,|A^{-1}|_\delta)>0$, and $\mathcal S\subseteq \Lambda$ with Hausdorff dimension zero, such that  if $F\in C^{\omega}_{h,\delta}(\T^{d}\times\Lambda,gl(m,\C))$ satisfying $|F|_{h,\delta}\leq\varepsilon_{0}$ and   $\lambda\in\Lambda\backslash\mathcal{S}$, the cocycle $(\alpha,  A(\lambda)+F(\cdot, \lambda))$ is reducible, i.e., there exists $B_{\lambda}\in C_{\frac{h}{4}}^{\omega}(\T^{d},\GL(m,\C))$ such that
$$
B_\lambda^{-1}(\cdot+\alpha)(A(\lambda)+F(\cdot, \lambda))B_\lambda(\cdot)=\tilde{A}(\lambda)\in \mathrm{GL}(m,\C).
$$
In addition, the $\tilde{A}(\lambda)$ has simple eigenvalues.
\end{Theorem}

 \begin{Remark}
 The non-degenerate condition (\ref{def-non-degeneracy}) is generic in $C^\omega(\Lambda, \GL(m,\C))$  \cite{HeY}, and thus Theorem \ref{Thm2} solves Eliasson's conjecture (Conjecture \ref{conj-eli-3}).
 \end{Remark}

We would like to provide some insight on both our result and its proof. Specifically, we in fact prove a bit more,  i.e.,  non-reducible cocycles and reducible cocycles with multiple eigenvalues form a zero Hausdorff dimensional set, and this is  essential for our spectral applications. Additionally, a crucial aspect of our proof involves demonstrating strong almost reducibility in the strip $|\Im z|< \frac{h}{2}$ for any $\lambda\in\Lambda$ to establish the analytic reducibility of the cocycle $(\alpha, A(\lambda)+F(\lambda,\cdot))$  in the strip $|\Im z|< \frac{h}{4}$ for any $\lambda\in\Lambda\backslash\mathcal{S}$.

Recall that
$(\alpha,A)\in C^\omega(\T^d, \GL(m,\C))$ is  almost reducible  if there exist $B_j\in C_{h_j}^\omega(\T^d,$ $\GL(m,\C))$, $A_j\in \GL(m,\C)$ and $F_j\in C_{h_j}^\omega(\T^d,gl(m,\C))$ such that
$$
B_j^{-1}(x+\alpha)A(x)B_j(x)=A_j+F_j(x),
$$
with $|F_j|_{h_j}\rightarrow 0$ and $A_j\rightarrow A_{\infty} \in \GL(m,\C)$. If $h_j \rightarrow 0$, then $(\alpha,A)$ is said to be (weak)-almost reducible or $C^{\infty}$ almost reducible; if $h_j \rightarrow h_*>0$, then $(\alpha,A)$ is said to be (strong)-almost reducible. We in fact have the following:

\begin{Theorem}\label{thm-1}
Suppose all assumptions of Theorem \ref{Thm2}. If  $F\in C^{\omega}_{h,\delta}(\T^{d}\times\Lambda, gl(m,\C))$ satisfies $|F|_{h,\delta}\leq\varepsilon_{0}$, where $\varepsilon_{0}$ depends on $\alpha,d,h,\delta, m, c,r,A$, then  for every $\varsigma>0, \epsilon>0$,  there exists  a partition of $\Lambda$, denoted by $\Pi$, and $\eta=\eta(\epsilon, \varsigma, \alpha, h,\delta, m, c,r,  A)>0$, such that for every $\tilde\Lambda\in\Pi$, there is  
$B\in C^{\omega}_{h/2,\eta}(\T^{d}\times\tilde\Lambda, \GL(m,\C))$ such that
$$
B^{-1} (\cdot+\alpha,\lambda)(A(\lambda)+F(\cdot, \lambda))B(\cdot,\lambda)=\tilde{A}(\lambda)+\tilde{F}(\cdot,\lambda),
$$
where $|\tilde F|_{h/2,\eta}<\epsilon$, $|B|_{h/2,\eta}|\tilde F|_{h/2,\eta}^\varsigma<1$ and $|B^{-1}|_{h/2,\eta}| \tilde F|_{h/2,\eta}^\varsigma<1$.
\end{Theorem}

\begin{Remark}
 We will prove a  \textit{stratified and quantitative} almost reducibility result, and the precise version can be seen in Proposition \ref{pro-iterative}.
\end{Remark}
\begin{Remark}
 With minor modifications,  the conjugation can be defined in any strip with width $h_*<h$, not necessarily $\frac{h}{2}$.

\end{Remark}
We shall provide an overview of previous results about almost reducibility. Eliasson \cite{Eli92} demonstrated that if an $\SL(2,\R)$ cocycle is close to constant, then it is (weak)-almost reducible. Leguil, You, Zhao, and Zhou \cite{LYZZ} later proved the (strong)-almost-reducible version (see also \cite{CCYZ}). While these results are perturbative (dependent on $\alpha$), the non-perturbative version such as the one-frequency case was proved by Avila and Jitomirskaya \cite{aj1}. We must also note that Avila's \textit{Almost Reducibility Conjecture}(ARC) is the global version of almost reducibility. It states that for any \textit{subcritical cocycle}, (strong)-almost reducibility holds. The proof of ARC was announced in \cite{A0} and set to appear in \cite{A1, A3}. To understand its various spectral applications, one may refer to the survey \cite{you}. Furthermore, Eliasson demonstrated (weak)-almost reducibility for quasi-periodic cocycles that take values in higher dimensional groups $\GL(m,\R)$ close to constants. Recently, the (strong)-almost reducibility in this context has also been established \cite{Cha2, GYZ}.

In conclusion, we would like to address the proof. Despite the lack of explicit mention, it has been shown that Conjecture \ref{conj-eli-3} holds for the $\SL(2,\R)$ group case by Avila in \cite{A2}. His proof strongly depends on the existence of fibred rotation number. However, in higher dimensions, the concept of a fibred rotation number is non-existent. While the Maslov index can be defined in the context of symplectic groups \cite{FJN, Pa}, it is not sufficient for reducibility as it
involves all eigenvalues $\{\sigma_i(\lambda)\}_i$ of the constant part.
In higher dimensions, a more detailed comprehension of the non-reducible set, or the exceptional set $\mathcal{S}$, is required. This kind of  \textit{stratified} almost reducibility was initiated by Krikorian \cite{K2}. Different from \cite{K2}, we are able to give quantitative (strong)-almost reducibility for $\GL(m,\C)$ cocyles, while the almost reducibility result in \cite{K2} is the weak version and restricted to the semi-simple case.
As highlighted in \cite{you}, quantitative estimates are crucial for spectral applications. With the scheme we have developed in this paper, it is our expectation that further applications can be explored.

\section{Preliminaries}

\subsection{Linear cocycles.}

We consider the linear cocycle $(\alpha,A):\T^d\times\C^{m}\rightarrow\T^d\times\C^{m}$ defined as
\[(\theta, u)\mapsto (\theta+\alpha, A(\theta)u).\]
The cocycle iterations are given by $(\alpha,A)^{n}=(n\alpha,A(\cdot;n))$, where 
$$
\left\{\begin{array}{lr}A(\theta;n)=A(\theta+(n-1)\alpha)\cdots A(\theta)   &n>0,\\A(\theta;n)=Id   &n=0,\\A(\theta;n)=A^{-1}(\theta+n\alpha)\cdots A^{-1}(\theta-\alpha) &{n<0}.\end{array} \right.
$$
Let $\gamma_{1}(\alpha,A)\geq\cdots\geq \gamma_{m}(\alpha,A)$ be the Lyapunov exponents of $(\alpha,A)$, repeated according to their multiplicity, i.e.,
$$
\gamma_{k}(\alpha,A)=\lim_{n\rightarrow\infty}\frac{1}{n}\int_{\T^d}\ln(\sigma_{k}(A(\theta;n)))d\theta,
$$
where $\sigma_{1}(B)\geq\cdots\geq\sigma_{m}(B)$ are $singular\ values$ of a matrix $B\in \mathrm{GL}(m,\C)$. For a matrix $B\in\mathrm{GL}(m,\mathbb{C})$, $\prod_{j=1}^{k}\sigma_{j}(B)=\|\Lambda^{k}B\|$, where $\Lambda^{k}B$ is the $k$-th exterior product  of $B$. Therefore, we have
$$
\Sigma_{j=1}^{k}\gamma_{j}=\lim_{n\rightarrow\infty}\frac{1}{n}\int_{\T^d}\ln\|\Lambda^{k}A(\theta;n)\|d\theta.
$$


%
%

  \subsection{Aubry duality}

Let $V(x)=\sum_{k=-\ell}^\ell\hat V_k e^{2\pi \ii kx}$. Suppose that the quasi-periodic long-range operator
\begin{equation}\label{op-2}
(L_{ \varepsilon^{-1}V,\alpha,x}^{{ W}} u)_n=\sum_{k\in\Z^d}\hat W_ku_{n-k}+\varepsilon^{-1}V(x+\la n,\alpha\ra)u_n, \ \ n\in\Z^d,
\end{equation}
has a $C^k$ quasi-periodic Bloch wave $u_n=e^{2\pi \ii\la n,\theta\ra}\overline{\psi}(x+\la n,\alpha\ra)$ for some $\overline{\psi}\in C^k(\T,\C)$ and $\theta\in \T^d$. It is easy to see that the Fourier coefficients of $\overline{\psi}$ satisfy the following long-range operator:
\begin{align}\label{op-1}
(L_{\varepsilon W,\alpha,\theta}^{V}u)_{n}=\Sigma_{k=-\ell}^\ell\hat V_{k}u_{n-k}+ \varepsilon W(\theta+n\alpha)u_{n}, \ \ n\in\Z.
\end{align}
We call $L_{\varepsilon W,\alpha,\theta}^{V}$  the dual operator of $L_{ \varepsilon^{-1}V,\alpha,x}^{{ W}}$. Denoting by $\hat{\mathcal{N}}(E)$ the IDS of  $L_{\varepsilon W,\alpha,\theta}^{V}$, it is  well-known  \cite{Pu}:
\begin{equation}\label{Thm1-2-1}
 \mathcal{N}(E)=\hat{\mathcal{N}}(E).
\end{equation}

The eigenvalue equations of (\ref{op-1}) are
\[
\Sigma_{k=-\ell}^{\ell}\hat V_{k}u_{n-k}+ \varepsilon W(\theta+n\alpha)u_{n}=E u_{n},\quad n\in\mathbb{Z}.
\]
Without loss of generality, we assume $\hat V_\ell\neq 0$. Then the eigenvalue equations can be viewed as a  $\mathrm{GL}(2\ell,\C)$-cocycle
$(\alpha, A_\varepsilon): \T^d\times \C^{2\ell}\rightarrow \T^d\times \C^{2\ell}$, where 
%
%
\[
A_{\varepsilon}(E,\theta)=\begin{pmatrix} -\frac{\hat V_{\ell-1}}{\hat V_{\ell}} & \cdots & -\frac{\hat V_{1}}{\hat V_{\ell}} & \frac{E-\hat V_0- \varepsilon W(\theta)}{\hat V_{\ell}} & -\frac{\hat V_{-1}}{V_{\ell}} & \cdots & -\frac{\hat V_{\ell+1}}{\hat V_{\ell}} & -\frac{\hat V_{-\ell}}{\hat V_{\ell}} \\ 1 &&&&&&&
\\& \ddots &&&&&&
\\&& 1 &&&&&
\\&&& 1&&&&
\\&&&& 1&&&
\\&&&&& \ddots&&
\\&&&&&& 1& 0
\end{pmatrix}.
\]
Let $\hat \gamma$ be the fibred entropy (i.e. sum of the positive Lyapunov exponents) of the corresponding cocycle $(\alpha, A_\varepsilon)$. Then it relates the IDS by 
Thouless formula \cite{HP}:
\begin{equation}\label{thou}
\hat \gamma (E)=\int_{\R}\ln|E-E'|d\hat {\mathcal N}(E')-\ln|\hat V_{\ell}|.
\end{equation}

\subsection{Hausdorff measure and Hausdorff dimension}
Let $(\mathcal X, \rho)$ be a metric space.  We denote by diam$(\Omega)$ the diameter of $\Omega$ for  any subset $\Omega\subseteq \mathcal X$.
\begin{Definition}\label{def-hausdorff}
For any $\Omega\subseteq \mathcal X$, any $\eta\in (0, \infty]$ and any $\varrho\in [0,\infty)$, let
\[
H_\eta^\varrho(\Omega):=\inf\{\sum_{i=1}^\infty (\mathrm{diam} \Omega_i)^\varrho \ : \ \Omega\subseteq \cup_i \Omega_i,\ \mathrm{and}\ \mathrm{diam}(\Omega_i)<\eta\}.
\]
Then $H^\varrho(\Omega):=\lim_{\eta\downarrow 0} H_\eta^\varrho(\Omega)$ is called the \textit{Hausdorff $\varrho$-dimensional measure} of $\Omega$, and $\mathrm{dim}_H(\Omega):=\inf\{\varrho: \ H^\varrho(\Omega)=0\}$ is called the Hausdorff dimension.
\end{Definition}

\section{Preparation lemmas}

In this section, we first introduce some useful concepts and preparation lemmas:

\subsection{Roots of an algebraic equation.}

In the following we consider polynomials in $C_{N}[X]$ of the form
$$
\chi(X)=a_0X^{n}+a_{1}X^{n-1}+\cdots+a_{n}, \quad n\leq N.
$$
Let $|\chi|:=\sup_i{|a_{i}|}$. We say $\chi(\lambda)( X)\in C_{\delta}^{\omega}(\Lambda,C_{N}[X])$ if its coefficients $a_{i}\in C_{\delta}^{\omega}(\Lambda)$. Given $\chi(\lambda)(X) \in C_{\delta}^{\omega}(\Lambda,C_{N}[X])$, we set
$$
|\chi|_{\delta}=\sup_i(|a_{i}|_{\delta}).
$$
Denote by $\Sigma_{\chi}$ be the multiset of zeroes of $\chi$ (counting the multiplicity). Conversely, for a given finite multiset $\Sigma$ whose elements $z\in\C$ may repeat, we call $\chi_{\Sigma}(X)=\prod_{z\in\Sigma}(X-z)$ the \textit{characteristic polynomial} of $\Sigma$. Supposing $\chi_{1},\chi_{2}\in C_N[X]$ with degrees $m_{1}, m_{2}\leq N$, we define their \textit{resultant} as
$$
\mathrm{Res}(\chi_{1},\chi_{2}):=\mathrm{Res}(\Sigma_{\chi_{1}},\Sigma_{\chi_{2}})
=\prod_{\sigma_{i}\in\Sigma_{\chi_{1}},\tau_{j}\in\Sigma_{\chi_{2}}}(\sigma_{i}-\tau_{j}).
$$
More generally, we denote
$$
\mathrm{Res}(\chi_{1},\chi_{2};u):=\mathrm{Res}(\Sigma_{\chi_{1}},e^{2\pi \ii u}\Sigma_{\chi_{2}})=\mathrm{Res}(\chi_{1},e^{2\pi \ii m_{2}u}\chi_{2}(e^{-2\pi \ii u}\cdot)).
$$
Suppose that there is a decomposition of $\Sigma$, i.e. $\Sigma=\Sigma_{1}\cup\cdots\cup\Sigma_{l}$.
Let $\chi:=\chi_{\Sigma }$ and $\chi_{i}:=\chi_{\Sigma_{i}}$.
For $u\in\T$, we have
\[
\mathrm{Res}(\chi,\chi;u)=\prod_{1\leq i,j\leq l}\mathrm{Res}(\chi_{i},\chi_{j};u).
\]

By direct computations, we can obtain the following estimates.
\begin{Lemma}[\cite{K2}]\label{resl}
Let $\chi,\chi_{1},\chi_{2}\in C_N[X]$ of degrees $n,m_{1},m_{2}\leq N$. Their corresponding multisets of zeroes  are $\Sigma,\Sigma_{1},\Sigma_{2}$. For $u\in\T$, we have the following:
\begin{eqnarray*}
&&|\chi_{e^{2\pi \ii u}\Sigma}|=|\chi|, \\
&&|\mathrm{Res}(\chi_{1},\chi_{2})|\leq (m_{1}+m_{2})!(1+|\chi_{1}|)^{m_{2}}(1+|\chi_{2}|)^{m_{1}}, \\
&&|\mathrm{Res}(\chi_{1},\chi_{2};u)|\leq (m_{1}+m_{2})!(1+|\chi_{1}|)^{m_{2}}(1+|\chi_{2}|)^{m_{1}}.
\end{eqnarray*}
If the degrees of $\chi_1',\chi_2'$ are also $m_1, m_2$, then  we have:
\begin{align*}
&|\mathrm{Res}(\chi_{1}',\chi_{2}';u)-\mathrm{Res}(\chi_{1},\chi_{2};u)|\\
\leq &(m_{1}+m_{2}+1)!(1+|\chi_{1}|)^{m_{2}}(1+|\chi_{2}|)^{m_{1}}\max\{|\eta_1|,|\eta_2|\},
\end{align*}
where $\eta_{1}=\chi_{1}'-\chi_{1}$, $\eta_{2}=\chi_{2}'-\chi_{2}$, provided $|\eta_1|, |\eta_1|\leq 1$.
\end{Lemma}

In this paper we always consider the case $A\in C^\omega_\delta(\Lambda, \GL(m,\C))$. Denote its spectrum $\Sigma(\lambda):=\Sigma(A(\lambda))$ which surely does not contain zero.  Thus there is $R>0$, such that $\Sigma(\lambda)\subseteq D(R)$ for $\lambda\in W_\delta(\Lambda)$, where 
 $$D(R):=\{z\in\C:\frac{1}{R}\leq|z|\leq R\}.$$ 
 Then we have the following basic estimates:
 
 \begin{Lemma}\label{lem-difference-character-poly}
 Suppose that $\Sigma(\lambda)\subseteq D(R)$ for $A\in C^\omega_\delta(\Lambda, \GL(m,\C))$, where $\lambda\in W_\delta(\Lambda)$ and   $\Sigma(\lambda):=\Sigma(A(\lambda))$. Then  
 \begin{equation}\label{lem-difference-1}
 |\chi_{\Sigma}|_\delta\leq m! R^m,\qquad 
 \left|\frac{\partial \chi_\Sigma}{\partial X}\right|_\delta\leq m|\chi_{\Sigma}|_\delta.
 \end{equation}
 Furthermore, we have 
 \begin{equation}\label{lem-difference-2}
 |\chi_{\Sigma}-\chi_{\Sigma'}|_\delta\leq m! M^{m-1}|A-A'|_\delta,
  \end{equation}
 where $\Sigma'(\lambda):=\Sigma(A'(\lambda))$ and $M=\max\{1, |A|_\delta, |A'|_\delta\}$.
 \end{Lemma}
 \begin{pf}
 Suppose the characteristic polynomials of $A(\lambda)$ and $A'(\lambda)$ are
 \begin{align*}
 \chi_{\Sigma(\lambda)}(X)=X^{m}+a_{1}(\lambda)X^{m-1}+\cdots+a_{m}(\lambda), \\
 \chi_{\Sigma'(\lambda)}(X)=X^{m}+b_{1}(\lambda)X^{m-1}+\cdots+b_{m}(\lambda).  \end{align*}
 Then we have 
 \[
 \frac{\partial\chi_{\Sigma(\lambda)}}{\partial X}=m X^{m-1} +(m-1)a_1(\lambda) X^{m-2}+\cdots+ a_{m-1}(\lambda).
 \]
 Since $\Sigma(\lambda)\subseteq D(R)$ for $\lambda\in W_\delta(\Lambda)$, by Vieta's Formula, we can obtain that for any $1\leq j\leq m$,
 \[
 |a_j(\lambda)|\leq C_m^jR^j\leq m!R^j,
 \]
then \eqref{lem-difference-1} follows. 
To prove \eqref{lem-difference-2}, we recall the following result:

 \begin{Lemma}[\cite{Bh}, Proposition 20.3]\label{lem-bh}
Denote by $a_k$ and $b_k$  $(k=1,\cdots, m)$ the $k$-th coefficients of their characteristic polynomials of  $A$ and $B$ in $\GL(m,\C)$  respectively. Then we have 
 \[
 |a_k-b_k|\leq k C_m^k M^{k-1} \|B-A\|,
 \]
 where $M=\max\{\|A\|, \|B\|\}$.
 \end{Lemma} 
 As a consequence, for any $\lambda\in W_\delta(\Lambda)$ and $1\leq j\leq m$,
 \begin{eqnarray*}
 |a_j(\lambda)-b_j(\lambda)|\leq j C_m^j M^{j-1} \|A(\lambda)-A'(\lambda)\|\leq m! M^{j-1} \|A(\lambda)-A'(\lambda)\|.
  \end{eqnarray*}
 Then \eqref{lem-difference-2} follows directly.
  \end{pf}

Now we analyze the function $g(\lambda,u)$ defined in (\ref{trans-func-2}). In the case $u\in\Z$, we have
\begin{eqnarray}\label{equ-g-analytic-1}
g(\lambda, u)&=& \prod_{\sigma_i,\sigma_j\in \Sigma(\lambda), \atop{i\neq j}}(\sigma_i- \sigma_j) = \prod_{\sigma_i \in \Sigma(\lambda) } (\prod_{ \sigma_j\in \Sigma(\lambda), \atop{i\neq j}}(\sigma_i- \sigma_j) )\nonumber\\
&=&  \prod_{\sigma_i \in \Sigma(\lambda) } \frac{\partial\chi_{\Sigma(\lambda)}}{\partial X}(\sigma_i)=\textrm{Res}(\chi_{\Sigma(\lambda)}, \frac{\partial\chi_{\Sigma(\lambda)}}{\partial X}),\label{equ-g-0}
\end{eqnarray}
where the coefficients of $\chi_{\Sigma(\lambda)},\frac{\partial\chi_{\Sigma(\lambda)}}{\partial X}$ are all analytic in $\lambda$, and hence $g(\lambda,u)$.
 
 In the case $u\in\T \backslash  \Z $, it is obvious that the function $g(\lambda,u)$ can be expressed as
\begin{equation}\label{trans-func-1}
g(\lambda,u)=\frac{\mathrm{Res}(\chi_{\Sigma(\lambda)},\chi_{\Sigma(\lambda)};u)}{(\prod_{z\in\Sigma(\lambda)}z)(1-e^{2\pi \ii u})^{m}}=\frac{\mathrm{Res}(\chi_{\Sigma(\lambda)},\chi_{\Sigma(\lambda)};u)}{\det A(\lambda)(1-e^{2\pi \ii u})^{m}},
\end{equation}
where $\#\Sigma=m$. Moreover,  by the analyticity of $\mathrm{Res}(\chi_\Sigma(\lambda),\chi_\Sigma(\lambda);u)$ and $\det A(\lambda)$ in $\lambda$, together with $\det A(\lambda)\neq 0$, we have $g(\lambda,u)$ is also analytic in $\lambda$. 
In this case, we have the following observation:
\begin{Lemma}\label{lem-g-neq-0}
Suppose $\chi(X)=X^m+a_1(\lambda)X^{m-1}+\cdots+a_m(\lambda)$. For $u\in\T\backslash\Z$, we have
\[
\mathrm{Res }(\chi,\chi;u)= (e^{2\pi \ii   u}-1)^m \mathrm{Res}(\chi,\tilde \chi_u) .
\]
where \[\tilde \chi_u(X)= a_1(\lambda)X^{m-1}+  \frac{e^{4\pi \ii u}-1}{e^{2\pi \ii u}-1} a_2(\lambda)X^{m-2}+\cdots+ \frac{e^{2m\pi \ii u}-1}{e^{2\pi \ii u}-1} a_m(\lambda).\]
\end{Lemma}
\begin{pf}
Since 
\begin{align*}
e^{2\pi \ii mu}\chi (e^{-2\pi \ii u}X)&=X^m+a_1(\lambda)e^{2\pi \ii u}X^{m-1}+\cdots +a_{m-1}(\lambda)e^{2\pi \ii (m-1)u}X+e^{2\pi \ii mu}a_m(\lambda)\\
&=:\acute\chi_u(X),
\end{align*}
then by the definition of $\mathrm{Res}(\chi,\chi;u)$, we have 
\[
\mathrm{Res}(\chi,\chi;u)=\mathrm{Res}(\chi,\acute \chi_u).
\]
Now by the fact that for any two polynomials $f(X)=a_0X^n+a_1X^{n-1}+\cdots +a_n$, 
$g(X)=b_0X^m+b_1X^{m-1}+\cdots+b_m$, their resultant 
\begin{eqnarray*}
 \mathrm{Res}(f,g)  =  \left|\begin{array}{ccccccccc}a_0 & a_1 & a_2 & \cdots   & 0 & 0 &0\\
0 & a_0 & a_1 &    \cdots &0&0&0 \\
\vdots & \vdots & \vdots & \vdots&     &\vdots&\vdots\\
0&0&0&\cdots& a_{n-1}& a_n&0 \\
0&0&0&\cdots& a_{n-2}& a_{n-1}&a_n \\
b_0 &b_1 &  b_2 & \cdots&   0 &  0&0 \\0 &b_0 & b_1 &    \cdots &0&0&0 \\
\vdots & \vdots & \vdots & \vdots&     &\vdots&\vdots\\
0&0&0&\cdots&    b_{m-1}&   b_m&0 \\
0&0&0&\cdots& b_{m-2}& b_{m-1}&  b_m\end{array}\right|,
\end{eqnarray*}
we can obtain that 
\begin{align*}
 &\mathrm{Res}(\chi,\acute\chi_u) \\
 &=   (e^{2\pi \ii u}-1)^m \left|\begin{array}{ccccccccc}1 & a_1 & a_2 & \cdots   & 0 & 0 &0\\
0 & 1 & a_1 &    \cdots &0&0&0 \\
\vdots & \vdots & \vdots & \vdots&     &\vdots&\vdots\\
0&0&0&\cdots& a_{m-1}& a_m&0 \\
0&0&0&\cdots& a_{m-2}& a_{m-1}&a_m \\
0 & a_1 & \frac{e^{4\pi \ii u}-1}{e^{2\pi \ii u}-1}a_2 & \cdots&   0 &  0&0 \\0 & 0 & a_1 &    \cdots &0&0&0 \\
\vdots & \vdots & \vdots & \vdots&     &\vdots&\vdots\\
0&0&0&\cdots&   \frac{e^{2(m-1)\pi \ii u}-1}{e^{2\pi \ii u}-1}a_{m-1}&   \frac{e^{2m\pi \ii u}-1}{e^{2\pi \ii u}-1}a_m&0 \\
0&0&0&\cdots&   \frac{e^{2(m-2)\pi \ii u}-1}{e^{2\pi \ii u}-1}a_{m-2}&  \frac{e^{2(m-1)\pi \ii u}-1}{e^{2\pi \ii u}-1} a_{m-1}&  \frac{e^{2m\pi \ii u}-1}{e^{2\pi \ii u}-1}a_m\end{array}\right|\\
& =  (e^{2\pi \ii u}-1)^m\mathrm{Res}(\chi,\tilde\chi_u),
\end{align*}
where $\tilde \chi_u(X)= a_1(\lambda)X^{m-1}+  \frac{e^{4\pi \ii u}-1}{e^{2\pi \ii u}-1} a_2(\lambda)X^{m-2}+\cdots+ \frac{e^{2m\pi \ii u}-1}{e^{2\pi \ii u}-1} a_m(\lambda) $.
\end{pf}

\subsection{Transversality}
 
 As we will see in the proof, the transversality  of the function $g(\lambda,u)$ is crucial for us. Here we  first  generalize   the notion of transversality   introduced in \cite{Eli97,Eli02, K3,K2}  following Pyartii \cite{Py}.

\begin{Definition}
A function $f:(a,b)\rightarrow\C$ is  said to be $(C,c,r)$-Pyartli, if $f\in C^{r+1}$ and for any $ x\in(a,b)$, we have
\begin{align*}
&\sup_{0\leq j\leq r+1}|\partial^{j}f(x)|\leq C, \\
&\sup_{0\leq j\leq r}|\partial^{j}f(x)|\geq c>0.
\end{align*}
\end{Definition}

Pyartli function will imply some good estimates  of preimages:
\begin{Lemma}\label{Lem1}
Let $f$ be $(C,c,r)$-Pyartli on $(a,b)$. Then for $0<\varsigma\leq \frac{c}{2}$, there is a disjoint union of intervals $\cup_{i\in J}I_{i}$ such that
\begin{align*}
\# J&\leq 2^{r}(\frac{2C(b-a)}{c}+1), \\
\max_{i\in J}|I_{i}|&\leq 2( \frac{2\varsigma}{c})^{\frac{1}{r}}, \\
|f(x)|&\geq \varsigma, \quad \forall x\in(a,b)\backslash\cup_{i\in J} I_{i}.
\end{align*}
\end{Lemma}

\begin{pf}
Similar estimates appeared in \cite{Eli02,K2}, we include the proof just for completeness.
If $r=0$, then the result is obvious. For $r\geq 1$, we assume that there exists $x_0\in (a, b)$ with $|x_0-a|\leq \frac{c}{2C}$, and $1\leq r_0\leq r$ such that 
\[
|\partial^{r_0}f(x_0)|\geq c.
\]
Otherwise, for all $x\in (a, b)\cap (a, a+\frac{c}{2C}]$, we have $|f(x)|\geq c$, and we only need to consider the interval $(a+\frac{c}{2C}, b)$ instead of $(a,b)$. Here, without loss of generality, we assume that $b-a>\frac{c}{2C}$. 

Now,   let 
$$
\tilde{f}(x):=\Re (e^{-\ii\mathrm{arg}\partial^{r_0} f(x_0)}f(x)).
$$
Then we have $\partial^{r_0}\tilde{f}(x)\geq \frac{c}{2}$ for  $x\in I=(a, x_0+\frac{c}{2C}]\cap (a,b)$. Now we consider $\partial^{r_0-1}\tilde{f} $ in $I$. Due to the fact that $\partial^{r_0-1}\tilde{f} $ is monotonic in $I$, there is at most one interval $I_{1,1}\subseteq I$ with minimal length, such that
\[
|\partial^{r_0-1}\tilde{f}(x)| \geq \frac{c}{2}(\frac{2\varsigma}{c})^{\frac{1}{r_0}}, \quad \forall x\in I\backslash I_{1,1}.
\]
Then by the fact that $\partial^{r_0}\tilde{f}(x)\geq\frac{c}{2}$ for $x\in I$,
we can get that 
\[
|I_{1,1}| \leq 2(\frac{2\varsigma}{c})^{\frac{1}{r_0}}.
\]
We continue to consider $\partial^{r_0-2}\tilde{f}$. It is monotonic on each component of $I\backslash I_{1,1}$. So there are at most two intervals $I_{2,1},I_{2,2}\subseteq I\backslash I_{1,1}$ such that
\begin{align*}
|\partial^{r_0-2}\tilde{f}(x)|&\geq \frac{c}{2}(\frac{2\varsigma}{c})^{\frac{2}{r_0}}, \quad \forall x\in I\backslash I_{1,1}\cup I_{2,1}\cup I_{2,2},\\
|I_{2,1}|,|I_{2,1}|&\leq 2(\frac{2\varsigma}{c})^{\frac{1}{r_0}}.
\end{align*}
Iterating the process for $r_0$ times, we get at most $2^{r_0}-1$ intervals $I_{i,j}$ such that
\begin{align*}
|\tilde{f}(x)|&\geq \varsigma, \quad \forall x\in I\backslash \cup_{1\leq i\leq r_0,1\leq j\leq 2^{i-1}} I_{i,j},\\
|I_{i,j}|&\leq 2(\frac{2\varsigma}{c})^{\frac{1}{r_0}}.
\end{align*}
It is obvious that
$$
|f(x)|= |e^{-\ii\mathrm{arg}\partial^{r_0}f(x_0)}f(x)|\geq |\Re (e^{-\ii\mathrm{arg}\partial^{r_0}f(x_0)}f(x))|=|\tilde{f}(x)|.
$$
Now we consider $f(x)$ on the interval $(a+\frac{c}{2C}, b)$, and repeat the above process as for $I$. In the end, the length of the considered interval is no more than $\frac{c}{2C}$ and the we do the above operation the last time.  Then the result follows. 
\end{pf}

Following \cite{K3,K2}, we also introduce a slightly stronger definition of transversality for a  better control of  the derivatives of functions (also its products).

\begin{Definition}
Let $f:(a,b)\rightarrow\C$ be an analytic function. We say $f$ is $(M,\delta,c,r)$-transverse if:
\begin{enumerate}
\item $f$ is $(M, \delta)$-bounded, i.e., $f\in C_{\delta}^{\omega}((a,b),\C)$ and $|f|_{\delta}\leq M$;
\item for any $x\in(a-\frac{\delta}{2},b+\frac{\delta}{2})$, $\sup_{0\leq j\leq r}|\partial^{j}f(x)|\geq c>0$.
\end{enumerate}
\end{Definition}

\begin{Remark}\label{tran-pya}
By Cauchy's estimate, if $f$ is $(M,\delta,c,r)$-transverse, then automatically  $f$ is $(\frac{(r+1)!M}{\min\{1,\delta^{r+1}\}},c,r)$-Pyartli.
\end{Remark}

The following  lemmas are  quite important for us: Lemma \ref{Lem2}  says that a product of transverse functions is also  transverse; Lemma \ref{Lem3} shows that if a product of "not too large" functions is transverse, then so is each of them.

\begin{Lemma}[{\cite{K3,K2}}]\label{Lem2}
Suppose $f_{1},\cdots,f_{l}$ are $(M_{i},\delta_{i},c_{i},r_{i})$-transverse. Let $M=\sup M_{i}$, $\delta=\inf\delta_{i}$, $c=\inf c_{i}$ and $r=r_{1}+\cdots+r_{l}$. Then $f=f_{1}\cdots f_{l}$ is $(M',\delta,c',r)$-transverse where
$$
M'=M^{l}, \quad\quad c'=((\frac{\delta}{4r^{2}M})^{rl}c)^{l^{r+1}}.
$$
\end{Lemma}
\begin{Lemma}[{\cite{K3,K2}}]\label{Lem3}
Suppose $f_{1},\cdots,f_{l}$ are functions belong to $C_{\delta}^{\omega}(\Lambda)$ with $|f_{i}|_{\delta}\leq M$. If $f=f_{1}\cdots f_{l}$ is $(M,\delta,c,r)$-transverse, then each $f_{i}$ is $(M,\delta,c',r)$-transverse with $c'=(\frac{2lM}{\delta})^{-rl}c$.
\end{Lemma}

If we restrict the function to $g(\lambda,u),$ then we extend the notation of transversality to the multiset $\Sigma(\lambda)$.

\begin{Definition}
We say a finite multiset $\Sigma(\lambda)$  is $(M,\delta,c,r)$-transverse  on $\Lambda$ if for any $u\in\T$, $\lambda\mapsto g(\lambda,u)$ is $(M,\delta,c,r)$-transverse on $\Lambda$.
\end{Definition}

Suppose that  $\Sigma(\lambda)$ can be decomposed into  $\Sigma(\lambda)=\Sigma_{1}(\lambda)\cup\cdots\cup\Sigma_{l}(\lambda)$.  
 Let $g_i(\lambda,u) \ (i=1,\cdots, l)$ be functions defined as in (\ref{trans-func-2}) with respect to $\Sigma_i(\lambda)$. Then one has 
\begin{equation}\label{equ-relation-g-g-i}
g(\lambda,u)=\left(\prod_{i=1}^l g_{i}(\lambda,u)\right)\left(\prod_{i, j=1,i\neq j}^l\mathrm{Res}( \chi_{\Sigma_{i}(\lambda)},\chi_{\Sigma_{j}(\lambda)};u)\right).
\end{equation}
This motivated the following definition:

\begin{Definition}
We say a decomposition $\Sigma(\lambda)=\Sigma_{1}(\lambda)\cup\cdots\cup\Sigma_{l}(\lambda)$ on $\Lambda$ is $(M,\delta,c,r)$-transverse if
\begin{enumerate}[(1)]
\item For all $1\leq i\leq l$, the multiset $\Sigma_{i}(\lambda)$ is $(M,\delta,c,r)$-transverse on $\Lambda$;
\item For all $i\neq j$ and all $u\in\T$, $\lambda\mapsto\mathrm{Res}(\chi_{\Sigma_{i}(\lambda)},\chi_{\Sigma_{j}(\lambda)};u)$ is $(M,\delta,c,r)$-transverse on $\Lambda$.
\end{enumerate}
\end{Definition}

Now, thanks to Lemma \ref{Lem2} and Lemma \ref{Lem3}, 
transversality of decomposition  will imply transversality of the multiset   $\Sigma(\lambda)$, and it is also true conversely. 

\begin{Lemma}\label{hebin}
Suppose $A\in C^\omega_\delta(\Lambda, \GL(m,\C))$ with $A=\mathrm{diag}\{A_{11}, \cdots, A_{ll}\}$ and $\Sigma(A(\lambda))=:\Sigma(\lambda)\subseteq D(R)$ for $\lambda\in W_\delta(\Lambda)$, where $A_{ii} (1\leq i\leq l)$ are block matrices, and $\Sigma_i(\lambda):=\Sigma(A_{ii}(\lambda))$ for  $1\leq i\leq l$. \begin{enumerate}
\item If the decomposition $\Sigma=\Sigma_1\cup\cdots\cup\Sigma_l$ is $(M,\delta, c,r)$-transverse on $\Lambda$, then $\Sigma$ is $((2R)^{m^2}, \delta, c', l^2r)$-transverse on $\Lambda$ with $c'= ( \ (\frac{4l^4r^2(2R)^{m^2}}{\delta} )^{-rl^4}c\  )^{l^{2l^2r+2}}$.
\item If the finite multiset $\Sigma$ is $(M, \delta, c, r)$-transverse on $\Lambda$, then the decomposition $\Sigma=\Sigma_1\cup\cdots\cup\Sigma_l$ is $((2R)^{m^2}, \delta, c'', r)$-transverse on $\Lambda$ with $c''= (\frac{2l^2(2R)^{m^2}}{\delta} )^{-rl^2}c$.
\end{enumerate}

\end{Lemma}
\begin{pf}
Notice that $\Sigma(\lambda)\subseteq D(R)$ implies 
 \begin{equation}\label{equ-esti-g}
|g(\cdot, u)|_\delta, |g_i(\cdot, u)|_\delta, |\textrm{Res}(\chi_{\Sigma},\chi_{\Sigma};u)|_{\delta}, |\textrm{Res}(\chi_{\Sigma_i},\chi_{\Sigma_j};u)|_{\delta}\leq (2R)^{m^2}.
\end{equation}
If   $\Sigma=\Sigma_1\cup\cdots\cup\Sigma_l$ is $(M,\delta, c,r)$-transverse on $\Lambda$, by (\ref{equ-esti-g}), it is also $((2R)^{m^2}, \delta, c,r)$-transverse on $\Lambda$. 
 Then by  (\ref{equ-relation-g-g-i}) and  Lemma \ref{Lem2}, we obtain that $g(\lambda, u)$ is $((2R)^{m^2l^2}, \delta, c', l^2r)$-transverse on $\Lambda$, which is also $((2R)^{m^2}, \delta, c', l^2r)$-transverse by (\ref{equ-esti-g}). By similar reasoning, (2) follows from Lemma \ref{Lem3}.

%
%
%
%
\end{pf}

\subsection{Separation of spectrum and normal form}

\begin{Definition}
Given $\nu>0$, we say a partition $\Sigma=\Sigma_{1}\cup\cdots\cup\Sigma_{l}$ is $\nu$-separated if   $|\sigma-\tau|>\nu$ for $\forall(\sigma,\tau)\in\Sigma_{i}\times\Sigma_{j}$ with $i\neq j$. In addition, if there exists $\zeta>0$ such that $\mathrm{diam}(\Sigma_{i})\leq\zeta$ for $1\leq i\leq l, $\footnote{Recall that for a set $\Sigma\subseteq \C$, we denote $\textrm{diam}(\Sigma)$ as the diameter of $\Sigma$.} then we say the decomposition of $\Sigma$ is $(\nu,\zeta)$-separated.
\end{Definition}
For any multiset $\Sigma$ with $\#\Sigma\leq m$ and $\mu>0$, there is always a $(\mu,\zeta)$-separated decomposition with $\zeta\leq m\mu$. We call it the \textit{maximal} $(\mu,\zeta)$-separated decomposition.

It is well-known that if the eigenvalues of $A(\lambda)$ is separated, then it  can be block diagonalized:

\begin{Lemma}[\cite{Eli02,HeY}]\label{lem-similar-transform}
Suppose $A\in C_\delta^\omega(\Lambda, \GL(m,\C))$ and the decomposition $\Sigma(A(\lambda))=:\Sigma(\lambda)=\Sigma_1(\lambda)\cup\cdots\cup\Sigma_l(\lambda)$, is $\nu$-separated for $\lambda\in W_\delta(\Lambda)$. If $\nu<\mathrm{const}\cdot |A|_\delta$, where $\mathrm{const}$ is a constant only depending on $m$, then there exists $S\in C_\delta^\omega(\Lambda, \GL(m,\C))$, such that for any $\lambda\in W_\delta(\Lambda)$,
\[
S^{-1}(\lambda)A(\lambda)S(\lambda)=\mathrm{diag}\{A_{11}(\lambda),\cdots,A_{ll}(\lambda)\},
\]
with   $\Sigma_{i}(\lambda)=\Sigma(A_{ii}(\lambda))$. Moreover, we have   estimates
\[
|A_{ii}|_\delta,\ |S|_{\delta}, \ |S^{-1}|_{\delta}\leq   \acute b(\frac{|A|_\delta}{\nu})^{m^2(m+2)},
\]
where the constant $\acute b> (120m)^{m^2+4m}$ only depends on $m$.
\end{Lemma} 

In the sequel of the article we will fix $\acute b> (120m)^{m^2+4m}$ to be a constant depending only on $m$. 
As a consequence, we have the following:

\begin{Proposition}\label{Nth}
Let $A\in C_{\delta}^{\omega}(\Lambda,\mathrm{GL}(m,\C))$ with $|A|_{\delta} \leq \tilde M$, where $\tilde M\geq 1$. Suppose that $\Sigma(\lambda):=\Sigma(A(\lambda))$ is $(M,\delta,c,r)$-transverse on $\Lambda$ and $\Sigma(\lambda)\subseteq D(R)$ for $\lambda\in W_\delta(\Lambda)$. For any $0<\nu'\leq 1$, we let
\[
\begin{matrix}
\begin{array}{l}
M'=(2R)^{m^2}, \\ c'=(\acute b R\nu'^{-1}\delta^{-1}\tilde M)^{-rm^3(m+6)}c,\\ \zeta'=10m\nu',
\end{array}
&
\begin{array}{l}
\delta'=\acute b^{-1}((R^2\tilde M)^{-1}\nu')^{3m}\delta, \\ r'=r,\\ \tilde{M}'=\acute b (\nu'^{-1}\tilde{M})^{m^2(m+2)},
\end{array}
\end{matrix}
\]
where $\acute b> (120m)^{m^2+4m}$ is the constant in Lemma \ref{lem-similar-transform}  depending only on $m$.
If $\nu'<\const\cdot \tilde M$, where $\const$ is a constant only depending on $m$, then there exists a partition $\Pi=\{\tilde\Lambda\}$ of $\Lambda$ with 
\[
\#\Pi\leq   \frac{|\Lambda|}{\delta'}+1,
\]
such that  for any $\lambda\in W_{\delta'}(\tilde{\Lambda})$,  $\Sigma(\lambda)$ has a decomposition 
$$\Sigma(\lambda)=\Sigma_{1}(\lambda)\cup\cdots\cup\Sigma_{l}(\lambda),$$
where  $l$ only depends on  $\tilde\Lambda$, satisfying the following properties:
\begin{enumerate}[(1)]
\item The decomposition $\Sigma(\lambda)=\Sigma_{1}(\lambda)\cup\cdots\cup\Sigma_{l}(\lambda)$ is $(\nu',\zeta')$-separated for all $\lambda \in W_{\delta'}(\tilde{\Lambda})$ and $(M',\delta',c',r')$-transverse in $\tilde \Lambda$.
\item There exists $\lambda_{0}\in\tilde{\Lambda}$ such that the decomposition $\Sigma(\lambda_0)=\Sigma_{1}(\lambda_0)\cup\cdots\cup\Sigma_{l}(\lambda_0)$ is  $(8\nu',\zeta')$-separated.   
Moreover, for any $\lambda\in W_{\delta'}(\tilde{\Lambda})$,  $\#\Sigma_{i}(\lambda)=\#\Sigma_{i}(\lambda_{0})$ and 
\[d_H(\Sigma_{i}(\lambda),\Sigma_{i}(\lambda_{0}))<\nu'. \footnote{We denote $d_H(\Omega_1,\Omega_2)$ as the Hausdorff distance between $\Omega_1$ and $\Omega_2$. }\]
 \item There exists $S\in C_{\delta'}^{\omega}(\tilde{\Lambda},\mathrm{GL}(m,\C))$ such that
\[
    S^{-1}(\lambda)A(\lambda)S(\lambda)=\mathrm{diag}\{A_{11}(\lambda),\cdots,A_{ll}(\lambda)\},
\]
with  $\Sigma_{i}(\lambda)=\Sigma(A_{ii}(\lambda))$ and  $|A_{ii}|_{\delta'},  |S|_{\delta'},|S^{-1}|_{\delta'}\leq \tilde{M}'$.
\end{enumerate}
\end{Proposition}

\begin{pf}
Due to $\Sigma(\lambda)\subseteq D(R)$, for any $u\in\T$, we have 
\[
|g(\lambda,u)|_\delta\leq(2R)^{m^{2}},
\]  
and by Lemma \ref{lem-difference-character-poly}, we also have $|\chi_{\Sigma(\lambda)}|_{\delta}\leq m!R^{m}$.
We divide the interval $\Lambda$ into $1+[\frac{|\Lambda|}{\delta'}]$ intervals with equal length no more than $\delta'$. Then  we have
\[\#\Pi\leq 1+\frac{|\Lambda|}{\delta'}=1+\acute b (R^2\tilde M\nu'^{-1})^{3m}\frac{|\Lambda|}{\delta}.\]

On each subinterval $\tilde{\Lambda}$, we now give the decomposition of the multiset $\Sigma(\lambda)$ on $W_{\delta'}(\tilde\Lambda)$.
Fix some $\lambda_{0}\in\tilde\Lambda$, take the   maximal $(8\nu',8m\nu')$-separated decomposition $\Sigma(\lambda_{0})=\Sigma_{1}(\lambda_{0})\cup\cdots\cup\Sigma_{l}(\lambda_{0})$, and define
\begin{align}\label{Nsd}
\Sigma_{i}(\lambda):=W_{\nu'}(\Sigma_{i}(\lambda_{0}))\cap\Sigma(\lambda)\subseteq \Sigma(\lambda).
\end{align}
Now for any $ \lambda\in W_{\delta'}(\tilde{\Lambda})$ and $\zeta\in\Sigma(\lambda)$, by Cauchy's estimate and (\ref{lem-difference-1}), we have
\begin{eqnarray}\label{equ-est-different-parameter}
\prod_{\sigma\in\Sigma(\lambda_0)}|\sigma-\zeta| &=& |\chi_{\Sigma(\lambda_0)}(\zeta)-\chi_{\Sigma( \lambda)}(\zeta)|\\
\nonumber&\leq& m \max\{1, |\zeta| ^{m-1}\} \cdot |\chi_{\Sigma(\lambda)}-\chi_{\Sigma(\lambda_{0})}|_{\delta'} \\
\nonumber&\leq& m  R^{m-1}|\frac{d}{d\lambda}\chi_{\Sigma(\lambda)}|_{\delta'}|\lambda-\lambda_{0}| \\
\nonumber &\leq& \frac{2m R^{m-1}}{\delta }|\chi_{\Sigma(\lambda)}|_{\delta}|\lambda-\lambda_{0}|\\
\nonumber&\leq& \frac{4m\cdot m!R^{2m-1}\delta'}{\delta}<\nu'^{m}.
\end{eqnarray}
This implies that for any $\zeta\in\Sigma(\lambda)$, there exists $\sigma\in\Sigma(\lambda_0)=\Sigma_{1}(\lambda_{0})\cup\cdots\cup\Sigma_{l}(\lambda_{0})$ such that $|\zeta-\sigma|<\nu'$,  which means that $\zeta\in W_{\nu'}(\Sigma_i(\lambda_0))$ for some $i=1,\cdots, l$, that is, $\Sigma(\lambda)\subseteq  \Sigma_{1}(\lambda)\cup\cdots\Sigma_l(\lambda)$. On the other hand, 
 $\Sigma(\lambda_{0})=\Sigma_{1}(\lambda_{0})\cup\cdots\cup\Sigma_{l}(\lambda_{0})$ is $8\nu'$-separated, which implies that $\Sigma_i(\lambda)\cap\Sigma_j(\lambda)=\emptyset$ for $i\neq j$. Consequently,  we get that 
\[\Sigma(\lambda)=\Sigma_{1}(\lambda)\cup\cdots\cup\Sigma_{l}(\lambda),\]
is a  decomposition for any $ \lambda\in W_{\delta'}(\tilde{\Lambda})$.


By similar arguments, for any $1\leq i\leq l$, $\tilde\sigma\in \Sigma_i(\lambda_0)\subseteq \Sigma(\lambda_0)$,  there exists at least one $\sigma\in\Sigma(\lambda)$ such that $|\sigma-\tilde\sigma|<\nu'$.
This just means  $\sigma\in \Sigma_i(\lambda)$ by definition, and thus
$\Sigma_i(\lambda_0)\subseteq W_{\nu'}(\Sigma_i(\lambda))$. 
Together with the fact that $\Sigma_i(\lambda)\subseteq W_{\nu'}(\Sigma_i(\lambda_0))$, we have 
\[d_H(\Sigma_i(\lambda), \Sigma_i(\lambda_0))<\nu'.\]
%

Next we show  this decomposition is indeed $(6\nu', (8m+2)\nu')$-separated.  For any $  \lambda_1,\lambda_2\in W_{\delta'}(\tilde{\Lambda})$, $\zeta\in\Sigma_i(\lambda_1)$ and $\zeta'\in\Sigma_j(\lambda_2)$, there exist
$\sigma \in\Sigma_i(\lambda_0)$ and $\sigma'\in\Sigma_j(\lambda_0)$, such that $|\zeta-\sigma|<\nu'$, $|\zeta'-\sigma'(\lambda_0)|<\nu'$. Then for $i\neq j$, we have
\begin{equation}\label{inequ-thm3-1}
|\zeta-\zeta'|
 \geq |\sigma-\sigma'|-|\zeta-\sigma|-|\sigma'-\zeta'|>6\nu',
\end{equation}
and for $i=j$, we have
\begin{equation*}
 |\zeta-\zeta'|
 \leq |\sigma-\sigma'|+|\zeta-\sigma|+|\sigma'-\zeta'|<(8m+2)\nu',
\end{equation*}
and they also hold if $\lambda_1=\lambda_2$. Furthermore, by $(\ref{inequ-thm3-1})$ and the continuity of $\Sigma(\lambda)$, we have
 $\#\Sigma_i(\lambda)=\#\Sigma_i(\lambda_{0})$ for any $\lambda\in W_{\delta'}(\tilde\Lambda)$ and $i=1,\cdots,l$. 

Once we have these, (3) directly follows from Lemma \ref{lem-similar-transform}, and  the transversality of decomposition follows from 
 Lemma \ref{hebin}. We finish the whole proof.
\end{pf}

In the following, we will prove that separability and transversality of the decomposition is stable with good estimates:

\begin{Lemma}\label{ndh}
Suppose that   
\[
A=\diag\{A_{11}, \cdots, A_{ll}\} \in C_\delta^\omega(\Lambda, \GL(m,\C)) ,
\]
 \[
 A'=\mathrm{diag}\{A_{11}', \cdots, A'_{ll}\}\in C_\delta^\omega(\Lambda, \GL(m,\C)),
 \]
  where $A_{ii}, A_{ii}'$ are block matrices with the same order for $1\leq i\leq l$. 
    Let $R, \tilde M>1$.
  If  $|A|_\delta\leq \tilde M$, $|A_{ii}-A_{ii}'|_\delta\leq\epsilon <1$, $\Sigma(A(\lambda))=:\Sigma(\lambda)\subseteq D(R)$ for $\lambda\in W_\delta(\Lambda)$, and
 \begin{equation}\label{equ-lem-8-1}
64m^2\tilde M^2R\epsilon^{\frac{1}{m}}<1,
\end{equation}
then the following hold:
\begin{enumerate}
\item For any $\lambda\in W_\delta(\Lambda)$, 
\begin{equation}\label{equ-distance-hausdorff}
d_H(\Sigma_i(\lambda), \Sigma_i'(\lambda))\leq 4m^2\tilde M^2\epsilon^{\frac{1}{m}}, \ (\forall 1\leq i\leq l),
\end{equation}
where $\Sigma_i(\lambda):=\Sigma(A_{ii}(\lambda))$, and $\Sigma_{i}'(\lambda):=\Sigma( {A}_{ii}'(\lambda))$. Consequently, 
\begin{equation}\label{diff}  \Sigma'(\lambda):=\Sigma(A'(\lambda))  \subseteq D(R')
\end{equation}
  with $R'=R+8m^2 \tilde M^2R^2\epsilon^{\frac{1}{m}}$.
\item If   $\Sigma=\Sigma_{1}\cup\cdots\cup\Sigma_{l}$ is $(\nu,\zeta)$-separated in $W_\delta(\Lambda)$, then   $\Sigma'=\Sigma_{1}'\cup\cdots\cup\Sigma_{l}'$ is $(\nu',\zeta')$-separated in $W_{\delta}(\Lambda)$, with 
\[
\nu'=\nu-8m^2  \tilde{M}^2\epsilon^{\frac{1}{m}}, \ \  \zeta'=\zeta+8m^2   \tilde{M}^2\epsilon^{\frac{1}{m}}.
\]
\item If the decomposition $\Sigma=\Sigma_{1}\cup\cdots\cup\Sigma_{l}$ is   $(M,\delta,c,r)$-transverse on $\Lambda$, then $\Sigma'=\Sigma_{1}'\cup\cdots\cup\Sigma_{l}'$ is $(M',\delta,c',r)$-transverse on $\Lambda$, with 
\[M'=(2R')^{m^{2}}, \ \ c'=c-2^{8m}m^{10m^2}R^{3m^2}\tilde M^m(\frac{2r}{\delta})^{r}\epsilon.\]
\item For any $1\leq i\leq l$, $u\in\T$, we have 
\[
|g_i(\cdot,u)-g_i'(\cdot, u)|_\delta\leq 2^{8m}m^{10m^2}R^{3m^2}\tilde M^m\epsilon,
\]
where we denote 
\[ g_i(\lambda, u)=\prod_{\sigma_{\ell_1},\sigma_{\ell_2}\in \Sigma_i(\lambda) \atop{\ell_1\neq \ell_2}}(\sigma_{\ell_1}-e^{2\pi\ii u}\sigma_{\ell_2}),\ \  g_i'(\lambda, u)=\prod_{\sigma_{\ell_1},\sigma_{\ell_2}\in \Sigma_i'(\lambda) \atop{\ell_1\neq \ell_2}}(\sigma_{\ell_1}-e^{2\pi\ii u}\sigma_{\ell_2}).
\]
\end{enumerate}
\end{Lemma}
\begin{pf}
  Due to $|A_{ii}- {A}_{ii}'|_{\delta}\leq\epsilon$ and $|A_{ii}|_\delta\leq \tilde M$, by (\ref{lem-difference-2}), we have
\begin{equation}\label{equ-esti-diffrence-sigma}
|\chi_{\Sigma_i'}-\chi_{\Sigma_i}|_{\delta}\leq m!(\tilde{M}+\epsilon)^{m-1}\epsilon.
\end{equation}
To prove \eqref{equ-distance-hausdorff}, by symmetry, we only need to prove for any $\sigma'\in\Sigma_{i}'(\lambda)$, 
 there exists  $\sigma\in\Sigma_{i}(\lambda)$ such that 
\begin{equation}\label{es1}|\sigma-\sigma'|\leq 4m^2\tilde M^2\epsilon^{\frac{1}{m}}.\end{equation} 
Indeed, this follows directly from the following estimate 
\begin{eqnarray*}
\nonumber|\prod_{\sigma\in\Sigma_{i}(\lambda)}(\sigma'-\sigma)|&=&|\chi_{\Sigma_{i}(\lambda)}(\sigma')-\chi_{\Sigma_{i}'(\lambda)}(\sigma')|\\
\nonumber&\leq&
m\cdot m!(\tilde{M}+\epsilon)^{m-1}\epsilon\max\{1,|\sigma'|^{m-1}\}\\
\nonumber&\leq& (4m^2\tilde M^2\epsilon^{\frac{1}{m}})^m,
\end{eqnarray*}
  where the first inequality holds by (\ref{equ-esti-diffrence-sigma}), and the second inequality holds because $|\sigma'|\leq |A'|_\delta\leq \tilde M+\epsilon$. As a consequence, by \eqref{equ-lem-8-1} and \eqref{es1}, 
 for any $\sigma'\in\Sigma'(\lambda)$, we have
\begin{eqnarray*}
|\sigma'| &\leq& R+4m^2 \tilde M^2\epsilon^{\frac{1}{m}}< R'\\
|\sigma'| &\geq& \frac{1}{R}-4m^2 \tilde M^2\epsilon^{\frac{1}{m}}>\frac{1}{R(1+8m^2\tilde M^2R\epsilon^{\frac{1}{m}})}=\frac{1}{R'},
\end{eqnarray*}
which just means $ \Sigma'(\lambda) \subseteq D(R')$. These finish the proof of (1), and (2) follows directly from (\ref{equ-distance-hausdorff}).

Now, we will prove (3) and (4).   For convenience we denote $\chi_{i}=\chi_{\Sigma_{i}}$ and $\chi_{i}'=\chi_{\Sigma_{i}'}$.  Since for $\lambda\in W_\delta(\Lambda)$ and  any $1\leq i\leq l$, $\Sigma_i(\lambda)\subseteq D(R)$,  by Lemma \ref{lem-difference-character-poly}, we have 
\begin{equation}\label{equ-chi-1}
|\chi_i|_\delta\leq m! R^m.
\end{equation}
For any $i, j\in\{1,2,\ldots,l\}$ and for any $u\in\T$, by \eqref{diff} we have  
 \begin{equation}\label{es2}
|\textrm{Res}(\chi_i', \chi_j';u)|_\delta
 =\sup_{\lambda\in W_\delta(\Lambda)}\prod_{\substack{\sigma_{k_1}\in\Sigma_i'(\lambda), \sigma_{k_2}\in \Sigma_j'(\lambda)}}|\sigma_{k_1}-e^{2\pi \ii u}\sigma_{k_2}|\leq (2R')^{m^2}.
\end{equation}
 Moreover, by Lemma \ref{resl} and (\ref{equ-esti-diffrence-sigma}), (\ref{equ-chi-1}), we have
\begin{eqnarray}\label{joint}
\nonumber\lefteqn{|\mathrm{Res}(\chi_{i}',\chi_{j}';u)-\mathrm{Res}(\chi_{i},\chi_{j};u)|_{\delta}}\\
\nonumber&\leq& (2m+1)!(1+|\chi_i|_\delta)^m (1+|\chi_j|_\delta)^m\max\{|\chi_i-\chi_i'|_\delta, |\chi_j-\chi_j'|_\delta\}\\
\nonumber&\leq & (2m+1)!(1+m!R^m)^{2m}m!(\tilde M+\epsilon)^{m-1}\epsilon\\
&\leq& 2^{8m}m^{6m^2}R^{2m^2}\tilde M^m\epsilon.
\end{eqnarray}
Then for any $\lambda\in W_{\frac{\delta}{2}}(\Lambda)$ and any $\ell\in\{0,1,\ldots,r\}$,  this implies that
\begin{eqnarray*}
|\partial_\lambda^\ell\textrm{Res}(\chi_i', \chi_j';u)|
&\geq& |\partial_\lambda^\ell\textrm{Res}(\chi_i, \chi_j;u)|-|\partial_\lambda^\ell(\textrm{Res}(\chi_i', \chi_j';u)-\textrm{Res}(\chi_i, \chi_j;u))|\\
&\geq& |\partial_\lambda^\ell\textrm{Res}(\chi_i, \chi_j;u)|-r!(\frac{2}{\delta})^r\cdot 2^{8m}m^{6m^2}R^{2m^2}\tilde M^m\epsilon.
\end{eqnarray*}
By the assumption, for $i\neq j$,  $\mathrm{Res}(\chi_i, \chi_j;u)$ is $(M,\delta,c,r)$-transverse on $\Lambda$, and thus for $\lambda\in W_{\frac{\delta}{2}}(\Lambda)\cap\R$,
\[\sup_{0\leq\ell\leq r }|\partial_\lambda^\ell\textrm{Res}(\chi_i', \chi_j';u)|\geq c-r!(\frac{2}{\delta})^r\cdot 2^{8m}m^{6m^2}R^{2m^2}\tilde M^m\epsilon\geq c'.\]
Combining with \eqref{es2},  we obtain $\mathrm{Res}(\chi_i', \chi_j';u)$ is $(M',\delta,c',r)$-transverse on $\Lambda$ for $i\neq j$.

We are left to prove the transversality of  $\Sigma_i(\lambda)$ and (4).
First, by \eqref{diff}, we have
$|g_i'(\lambda,u)|_{\delta}\leq (2R')^{m^2}.$ Now we denote $m_{i}:=\#\Sigma_{i}=\#\Sigma'_i$.

\textbf{Case a: ${u\in\Z}$.}  By (\ref{equ-g-0}), Lemma \ref{resl} and Lemma \ref{lem-difference-character-poly}, we have for any $\lambda\in W_\delta(\Lambda)$,
\begin{eqnarray*}
 |g_i(\lambda,u)-g_i'(\lambda,u)|&=&|\mathrm{Res}(\chi_i, \frac{\partial \chi_i}{\partial X})-\mathrm{Res}(\chi_i', \frac{\partial \chi_i'}{\partial X})|\\
&\leq & (2m_i)!(1+|\chi_i|_\delta)^{m_i-1}(1+|\frac{\partial \chi_i}{\partial X}|_\delta)^{m_i} m_i|\chi_i-\chi_i'|_\delta \\
&\leq &(2m)!(1+m\cdot m! R^m)^{2m-1}m\cdot m!(\tilde M+\epsilon)^{m-1}\epsilon\\
&<& 2^{8m}m^{10m^2}R^{2m^2}\tilde M^m\epsilon.
\end{eqnarray*}

\textbf{Case b: $u\in\T\backslash\Z$.} We denote 
\[
f_i(\lambda,u):=\frac{\mathrm{Res}(\chi_i(\lambda),\chi_i(\lambda);u)}{ (1-e^{2\pi \ii u})^{m_i}},\ \ \ f_i'(\lambda,u):=\frac{\mathrm{Res}(\chi_i'(\lambda),\chi_i'(\lambda);u)}{ (1-e^{2\pi \ii u})^{m_i}}. \] 
Then by Lemma \ref{lem-g-neq-0}, Lemma \ref{resl}, and Lemma \ref{lem-difference-character-poly}, we have for any $\lambda\in W_\delta(\Lambda)$,
\begin{eqnarray*}
|f_i(\lambda,u)-f_i'(\lambda,u)|&=&|\mathrm{Res}(\chi_i(\lambda),\tilde\chi_{i,u}(\lambda))-\mathrm{Res}(\chi_i'(\lambda),\tilde\chi_{i,u}'(\lambda))|\\
&\leq & (2m )!(1+|\chi_i|_\delta)^m(1+|\tilde\chi_{i,u}|_\delta)^m\max\{|\chi_i-\chi_i'|_\delta, |\tilde\chi_i-\tilde\chi_{i,u}'|_\delta\}\\
&\leq& (2m )!(1+(m+1)!R^m)^{2m}(m+1)!(\tilde M+\epsilon)^{m-1}\epsilon\\
&<& 2^{7m}m^{2m^2+6m}R^{2m^2}\tilde M^{m-1}\epsilon.
\end{eqnarray*}
By  \eqref{trans-func-1}, we can obtain that
 for any $\lambda\in W_\delta(\Lambda)$
\begin{eqnarray*}
|g_i(\lambda,u)-g_i'(\lambda,u)| &\leq& \frac{1}{|\det A_{ii}|}|f_i(\lambda,z)-f_i'(\lambda,z)|+\left|\frac{\det A_{ii}'-\det A_{ii}}{\det A_{ii}}\right||g_i'(\lambda,u)|\\
&\leq& R^m|f_i(\lambda,z)-f'_i(\lambda,z)|+R^m|\chi_i-\chi_i'|_\delta(2R')^{m^2}\\
&\leq& 2^{8m}m^{8m^2}R^{3m^2}\tilde M^{m}\epsilon.
\end{eqnarray*}
Then (4) follows in both cases. Furthermore, in both case, for any $\lambda\in W_{\frac{\delta}{2}}(\Lambda)\cap \R$ and any $\ell\in\{0,1,\ldots,r\}$, we can get
\begin{eqnarray*}
|\partial_\lambda^\ell g_i'(\lambda,u)|
&\geq& |\partial_\lambda^\ell g_i(\lambda,u)|-|\partial_\lambda^\ell(g_i'(\lambda,u)-g_i(\lambda,u))|\\
&\geq& |\partial_\lambda^\ell g_i(\lambda, u)|-r!(\frac{2}{\delta})^r\cdot 2^{8m}m^{10m^2}R^{3m^2}\tilde M^{m}\epsilon,
\end{eqnarray*}
implying $
\sup_{0\leq \ell\leq r}|\partial_\lambda^\ell g_i'(\lambda,u)|\geq c'.$
Therefore, $\Sigma'=\Sigma_{1}'\cup\cdots\cup\Sigma_{l}'$ is $(M',\delta,c',r)$-transverse on $\Lambda$, which finishes the proof.

\end{pf}

\section{Key iteration steps}\label{seckam}

In this section, we prove the key KAM iteration steps for stratified quantitative almost reducibility.

\subsection{Structure of resonance}\label{sub-sec-4-1}

In this subsection, we always assume that $l\leq m$, where $l$ is the number of the multisets $\{\Sigma_i\}$.
We now give the definitions of resonance and non-resonance between two sets.

\begin{Definition}
Fix $N,\varsigma>0$. We say two different multisets $\Sigma_{1}$ and $\Sigma_{2}$ are $(N,\varsigma)$-resonant if there are $\sigma_{1}\in\Sigma_{1}$ and $\sigma_{2}\in\Sigma_{2}$ such that for some $|k|\leq N$, 
\begin{equation}\label{inequ-def-1}
|\sigma_{1}-e^{2\pi \ii\langle k,\alpha\rangle}\sigma_{2}|<\varsigma.
\end{equation}
If such $k$ is unique, we say $\Sigma_{1}$ and $\Sigma_{2}$ are $(N, \varsigma)$-resonant with unique $k$, and we call $k$ as the resonant integer between $\Sigma_1$ and $\Sigma_2$. If it does not exist  such $k$, we say $\Sigma_{1}$ and $\Sigma_{2}$ are $(N,\varsigma)$-nonresonant. In addition, we say a multiset $\Sigma$ is $(N, \varsigma)$-nonresonant, if for any $\sigma_1, \sigma_2\in\Sigma$ and   $0<|k|\leq N$, (\ref{inequ-def-1}) does not hold.
\end{Definition}

It is obvious that $-k$ is the resonant integer between $\Sigma_2$ and $\Sigma_1$, if $k$ is a resonant integer between $\Sigma_1$ and $\Sigma_2$.\\

Before giving the structure of resonance, we introduce the following notation  for convenience.

\begin{Definition}\label{def-resonant-structrue}
Let $0<N\leq N'$, $K>0$. We say the decomposition
$$\Sigma(\lambda)=\Sigma_{1}(\lambda)\cup\cdots\cup\Sigma_{l}(\lambda)$$
 satisfies $\mathcal{H}(N,N',K)$ on $W_{\delta}(\Lambda)$, if  $\{1,\cdots,l\}$ can be divided into disjoint union:
$$ \{1,\cdots,l\} =  S_{1}\cup  S_{2} \cup \cdots \cup  S_{\tilde l},$$
and the following hold:
\begin{enumerate}[(a)]
\item $\Sigma_{i}(\lambda)$ is $(N',(2K)^{-1})$-nonresonant  for any $i\in\{1,2,\cdots, l\}$.
\item $\Sigma_{i}(\lambda)$, $\Sigma_{j}(\lambda)$ are $(mN',2mK^{-1})$-resonant with the same unique $k_{ij}$ for all $\lambda$, where $|k_{ij}|\leq mN$, if $i, j \in S_{l_1}$ with $i\neq j$.
\item $k_{ij}+k_{js}=k_{is}$, if  $i, j, s \in S_{l_1}$ are different.
\item $\Sigma_{i}(\lambda)$, $\Sigma_{j}(\lambda)$ are $(N',(2K)^{-1})$-nonresonant, if $i \in S_{l_1}$, $j\in S_{l_2}$ with $l_1\neq l_2$.
\end{enumerate}
\end{Definition}

For $N\in\N$, $j\in\N_0$ and $a>1$, let $N_{j}=a^{j}N$. The basic observation is that if the  decomposition is separated, then the decomposition will essentially satisfy some $\mathcal{H}(N,N',K)$ condition. The precise result is the following:

\begin{Lemma}\label{deco}
Let $\alpha\in\mathrm{DC}_{d}(\gamma,\tau)$. Suppose that the decomposition $\Sigma(\lambda)=\Sigma_{1}(\lambda)\cup\cdots\cup\Sigma_{l}(\lambda )$ is $(\nu,\zeta)$-separated and $\Sigma(\lambda)\subseteq D(R)$ for  $\lambda\in W_{\delta}( {\Lambda})$, and one can find $\lambda_{0}\in {\Lambda}$ such that 
\begin{equation}\label{equ-lem-condition-hausdorff-distance}
d_H(\Sigma_{i}(\lambda), \Sigma_{i}(\lambda_{0}))<\nu',
\end{equation}
 $\#\Sigma_{i}(\lambda)=\#\Sigma_{i}(\lambda_{0})$ for $\lambda\in W_{\delta}({\Lambda})$.  Then the decomposition satisfies $\mathcal H(N_p, N_{p+1}, K)$ in $W_\delta(\Lambda)$ for some integer $  p\in [0, m]$, if
\begin{align}\label{group}
8m\nu'<8m\zeta<mK^{-1}<\frac{\gamma}{10R(3mN_{m+1})^{\tau}}.
\end{align}
 
\end{Lemma}

%
\begin{pf}
We first check ($a$). For any $\lambda\in W_\delta(\Lambda)$, $i\in \{1,\cdots, l\}$, $\sigma_{1},\sigma_{2}\in\Sigma_{i}(\lambda)$, since $\alpha\in DC_d(\gamma,\tau)$, $\Sigma(\lambda)\subseteq D(R)$ and the decomposition is $(\nu,\zeta)$-separated, then for $0<|k|\leq mN_{m+1}$,   we have
\begin{eqnarray*}
|\sigma_{1}-e^{2\pi \ii\langle k,\alpha\rangle}\sigma_{2}|
&\geq& |\sigma_2||1-e^{2\pi \ii\la k,\alpha\ra}|-|\sigma_1-\sigma_2|\\
&\geq&\frac{\gamma}{2R(mN_{m+1})^{\tau}}-\zeta>K^{-1}.
\end{eqnarray*}
This just means $\Sigma_{i}(\lambda)$ is $(mN_{m+1},K^{-1})$-nonresonant with itself for any $\lambda\in W_\delta(\Lambda)$.

To verify ($b$)-($d$), we fist need the following observation:

\begin{Claim} \label{claim} For any $\lambda\in W_{\delta}( \Lambda)$,   if  $i\neq j$ and $\Sigma_{i}(\lambda)$, $\Sigma_{j}(\lambda)$ are $(mN_{m+1},2mK^{-1})$-resonant, then they are resonant with the unique $|k|\leq mN_{m+1}$. 
\end{Claim}
\begin{pf}
 If   there exist $\sigma_{1},\sigma_{1}'\in\Sigma_{i}(\lambda)$, $\sigma_{2},\sigma_{2}'\in\Sigma_{j}(\lambda)$ and distinct $|k_1|, |k_2|\leq mN_{m+1}$, such that
\[|\sigma_1-e^{2\pi \ii \la k_1,\alpha\ra}\sigma_2|, |\sigma_1'-e^{2\pi \ii\la k_2,\alpha\ra}\sigma_2'|<2mK^{-1},\]
then we obtain that 
\begin{align*}
\lefteqn{\frac{\gamma}{2R|k_{1}-k_{2}|^{\tau}}<|\sigma_{2}e^{2\pi \ii\langle k_{1}-k_{2},\alpha\rangle}-\sigma_2|
\leq |\sigma_{1}-e^{2\pi \ii\langle k_1,\alpha\rangle}\sigma_{2}|}\\
&&+|\sigma_{1}-\sigma_{1}'|+|\sigma_{1}'-e^{2\pi \ii\langle k_2,\alpha\rangle}\sigma_{2}'|+|\sigma_{2}'-\sigma_{2}|< 4mK^{-1}+2\zeta.
\end{align*}
This contradicts to (\ref{group}).  \end{pf}

Now we  denote $\Sigma_i(\lambda_0)$ by $\Sigma_i$ for convenience, and give the division of the set $\{1,\cdots,l\}$ according to $\lambda_0$. 
First of all, we need to introduce two concepts: $(L, \eta)$-connected and 
maximal $(L,\eta)$-connected component. 
For $L,\eta>0$, we say that $\Sigma_{i}$, $\Sigma_{j}$ are $(L,\eta)$-{\it connected} if there exists a $(L,\eta)$-resonant path of length $t$:
$$
\Sigma_{i_{0}},\Sigma_{i_{1}},\cdots,\Sigma_{i_{t}},\quad \textrm{with}\ \ i_{0}=i,i_{t}=j,
$$
such that $(\Sigma_{i_{0}},\Sigma_{i_{1}})$, $(\Sigma_{i_{1}},\Sigma_{i_{2}})$, $\cdots$, $(\Sigma_{i_{t-1}},\Sigma_{i_{t}})$ are all $(L,\eta)$-resonant. Notice that if $\Sigma_{i}$, $\Sigma_{j}$ are $(L,\eta)$-connected of length $t$ then they are $(tL,t\eta+(t-1)\zeta)$-resonant. Moreover, if $\Sigma_i, \Sigma_j$ and $\Sigma_j, \Sigma_s$ are both $(L,\eta)$-connected, then $\Sigma_i, \Sigma_s$ are also $(L,\eta)$-connected. Let $S:=\{1,\cdots,l\}$. We say a set $S'\subset S$ is a \textit{maximal} $(L,\eta)$-{\it connected component} of $S$ if the following hold:
\begin{enumerate}[1)]
\item For all $i,j\in S'$ with $i\neq j$, we have $\Sigma_{i},\Sigma_{j}$ are $(L,\eta)$-connected.
\item For all $i\in S'$ and $j\notin S'$, $\Sigma_{i},\Sigma_{j}$ are not $(L,\eta)$-connected.
\end{enumerate}

\smallskip

Let's finish our proof. We divide $S$ into $\ell_{i}$ maximal $(N_{i},K^{-1})$-connected components where $0\leq i\leq l$. Notice that if $S'$ is a maximal $(N_{i},K^{-1})$-connected component, then there exists some maximal $(N_{i+1},K^{-1})$-connected component $\tilde{S}$ such that $S'\subset \tilde{S}$. So we have $\ell_{i+1}\leq \ell_{i}$. Combining $1\leq \ell_{i}\leq l$, there exists $p$ with $0\leq p\leq l$ such that $\ell_{p}=\ell_{p+1}$. We let $\{S_{i}\}_{i=1}^{\tilde l}$ be the maximal $(N_{p},K^{-1})$-connected components. Then we have the following: 

\medskip

(1) If $i,j$ with $i\neq j$ belong to the same subset $S_{l_1}$, then $\Sigma_i, \Sigma_j$ are $(N_p,K^{-1})$-connected, and hence $(mN_p, mK^{-1}+(m-1)\zeta)$-resonant with a unique  $k_{ij}$, where the uniqueness follows from Claim \ref{claim}.  Then there exist $\sigma_1\in\Sigma_i, \sigma_2\in\Sigma_j$ such that 
\[
|\sigma_1-e^{2\pi \ii\la k_{ij},\alpha\ra}\sigma_2|<mK^{-1}+(m-1)\zeta.
\]
It is obvious that $k_{ji}=-k_{ij}$. For any $\lambda\in W_\delta(\Lambda)$,  by (\ref{equ-lem-condition-hausdorff-distance}), there exist $\tilde\sigma_1\in\Sigma_i(\lambda), \tilde\sigma_2 \in \Sigma_j(\lambda)$, so that 
$|\sigma_1-\tilde\sigma_1|, |\sigma_2-\tilde\sigma_2|<\nu'$, and thus
\[
|\tilde\sigma_1-e^{2\pi \ii\la k_{ij},\alpha\ra}\tilde\sigma_2|<mK^{-1}+(m+1)\zeta<2mK^{-1}.
\]
Therefore, again by Claim \ref{claim}, we have $\Sigma_i(\lambda)$ and $\Sigma_j(\lambda)$ are $(mN_{m+1}, 2mK^{-1})$-resonant with the same unique $k_{ij}$ for all  $\lambda\in W_\delta(\Lambda)$, where $|k_{ij}|\leq mN_p$.  We have verified ($b$) in $\mathcal H(N_p, N_{p+1}, K)$.

\medskip


\medskip

(2) If distinct $i,j,s$ belong to the same subset $S_{l_1}$, then
$\Sigma_{i}$, $\Sigma_{j}$ and $\Sigma_{s}$ are $(mN_p, mK^{-1}+(m-1)\zeta)$-resonant with each other,   with the resonant terms $k_{ij}, k_{js}, k_{si}$ respectively. This means that there exist $\sigma_{i},\sigma_{i}'\in\Sigma_{i}$, $\sigma_{j},\sigma_{j}'\in\Sigma_{j}$ and $\sigma_{s},\sigma_{s}'\in\Sigma_{s}$ such that
\begin{align*}
|\sigma_{i}-e^{2\pi \ii\langle k_{ij},\alpha\rangle}\sigma_{j}'|&< mK^{-1}+(m-1)\zeta, \\
|\sigma_{j}-e^{2\pi \ii\langle k_{js},\alpha\rangle}\sigma_{s}'|&<mK^{-1}+(m-1)\zeta, \\
|\sigma_{s}-e^{2\pi \ii\langle k_{si},\alpha\rangle}\sigma_{i}'|&< mK^{-1}+(m-1)\zeta,
\end{align*}
with $|k_{ij}|,|k_{js}|,|k_{si}|\leq mN_{p}$. So if $k_{ij}+k_{js}+k_{si}\neq 0$, we have
\begin{align*}
\frac{\gamma}{2R|k_{ij}+k_{js}+k_{si}|^{\tau}}\leq |\sigma_{i}(1-e^{2\pi \ii\langle k_{ij}+k_{js}+k_{si},\alpha\rangle})|
\leq 3m(K^{-1}+\zeta).
\end{align*}
This contradicts our assumption   \eqref{group}, and we verify ($c$) in $\mathcal H(N_p, N_{p+1}, K)$.


\medskip

(3) If  $i \in S_{l_1}$, $j\in S_{l_2}$, i.e. $i,j$ belong to different subsets, then $\Sigma_i, \Sigma_j$ are not $(N_{p+1},K^{-1})$-connected. In particular, they are $(N_{p+1},K^{-1})$-nonresonant. For any $\lambda\in W_\delta(\Lambda)$, $\tilde\sigma_i\in\Sigma_i(\lambda)$, $\tilde\sigma_j\in\Sigma_j(\lambda)$, by (\ref{equ-lem-condition-hausdorff-distance}), there exist $\sigma_i\in\Sigma_i$, $\sigma_j\in\Sigma_j$, such that $|\sigma_i-\tilde\sigma_i|, |\sigma_j-\tilde\sigma_j|<\nu'$. Then for  $|k|\leq N_{p+1}$, we have 
\[|\tilde\sigma_i-e^{2\pi \ii\la k,\alpha\ra}\tilde\sigma_j|\geq K^{-1}-2\nu'>(2K)^{-1}.\]
This means $\Sigma_i(\lambda), \Sigma_j(\lambda)$ are $(N_{p+1}, (2K)^{-1})$-nonresonant for any $\lambda\in W_\delta(\Lambda)$, and verifies ($d$) in $\mathcal H(N_p, N_{p+1}, K)$.
 We finish the proof.
\end{pf}

\subsection{Eliminate the non-resonant terms}\label{sub-sec-4-2}
For given $h,\delta,\eta>0$, $\alpha\in\R^d$, and $A\in C_\delta^\omega(\Lambda, \GL(m,\C))$, we decompose $C_{h,\delta}^\omega(\T^d\times\Lambda, gl(m,\C))=:\mathcal B_{h,\delta}=\mathcal B_{h,\delta}^{(nre)}(\eta)\oplus \mathcal B_{h,\delta}^{(re)}(\eta)$ (depending on $A, \alpha, \eta$) in such a way that for any $Y\in \mathcal B_{h,\delta}^{(nre)}(\eta)$, we have
\begin{equation}
AY, Y^+A\in \mathcal B_{h,\delta}^{(nre)}, \ \ |AY-Y^+A|_{h,\delta}\geq \eta |Y|_{h,\delta},
\end{equation}
where we use $Y^+$ to represent the function $Y(\cdot+\alpha)$ here and in the sequel.
Moreover, we denote $\mathbb P_{nre}$ and $\mathbb P_{re}$ as the standard projective operators from $\mathcal B_{h,\delta}$ onto $\mathcal B_{h,\delta}^{(nre)}(\eta)$ and $\mathcal B_{h,\delta}^{(re)}(\eta)$ respectively.
For any $N>0$, we define the truncating operators $\mathcal T_N$ on $C_{h,\delta}^\omega(\T^d\times\Lambda, gl(m,\C))$ as
\[(\mathcal T_Nf)(\theta,\lambda)=\sum_{k\in\Z^d, |k|\leq N}\hat F(k,\lambda)e^{2\pi \ii \la k,\theta\ra},\]
and $\mathcal R_N$ on $C_{h,\delta}^\omega(\T^d\times\Lambda, gl(m,\C))$ as
\[(\mathcal R_Nf)(\theta,\lambda)=\sum_{k\in\Z^d, |k|>N}\hat F(k,\lambda)e^{2\pi \ii \la k,\theta\ra}.\]
Furthermore, we assume $a\geq m$ in the sequel.

Using the above notations, we can apply the quantitative implicit function theorem in \cite{CCYZ} to get the following result:
\begin{Lemma} [\cite{CCYZ}] \label{lem-implicit function}
Assume that $\epsilon<(4\max\{1, |A|_\delta\})^{-4}$ and $\eta\geq 48\max\{1, |A|_\delta\}\epsilon^{1/2}$. Then for any $F\in \mathcal B_{h,\delta}$ with $|F|_{h,\delta}\leq \epsilon$, there exist $Y\in \mathcal B_{h,\delta}$ and $F^{(re)}\in \mathcal B_{h,\delta}^{(re)}(\eta)$ such that
\[e^{-Y^+}(A+F)e^Y=A+F^{(re)},\]
where $|Y|_{h,\delta}\leq \epsilon^{1/2}$ and $|F^{(re)}|_{h,\delta}\leq 3\epsilon\max\{1, |A|_\delta\}$.
\end{Lemma}


Then one can prove the following:

\begin{Proposition}\label{ehe}
Given $N, K, \tilde M>1, a\geq 1, h>h'>0$, assume that $\alpha\in\mathrm{DC}_{d}(\gamma,\tau)$, $A \in C_{\delta}^{\omega}(\Lambda,\mathrm{GL}(m,\C))$, $A(\lambda)=\mathrm{diag}\{A_{11}(\lambda),\cdots,A_{ll}(\lambda)\}$ and $|A|_{\delta}<\tilde{M}$. Let $\Sigma(\lambda):=\Sigma(A(\lambda))$, $\Sigma_{i}(\lambda):=\Sigma(A_{ii}(\lambda))$. Suppose that the decomposition $$\Sigma(\lambda)=\Sigma_{1}(\lambda)\cup\cdots\cup\Sigma_{l}(\lambda)$$
satisfies $\mathcal{H}(N, aN, K)$ on $W_{\delta}(\Lambda)$.
Then for any $F\in C_{h,\delta}^{\omega}(\T^{d}\times\Lambda, gl(m,\C))$  with
\[|F|_{h,\delta}\leq\epsilon<(12m^5\tilde M^2 K)^{-2m^2},\]
there exist $Y, f^{(re)},  F'\in C_{h,\delta}^{\omega}(\T^{d}\times\Lambda, gl(m,\C))$, such that
\begin{eqnarray}\label{nhe}
e^{-Y(\cdot+\alpha)} (A+ {F(\cdot)}) e^{Y(\cdot)} &= &A'+f^{(re)}(\cdot)+F'(\cdot)\\
 \nonumber&=& \diag \{ A_{11}', \cdots,   A_{ll}'\} +f^{(re)}(\cdot)+F'(\cdot)
\end{eqnarray}
 with estimates: 
 \[
\begin{matrix}
\begin{array}{l}
|A_{ii}-  A_{ii}'|_\delta\leq 3\tilde M\epsilon, \\ |Y|_{h,\delta}\leq \epsilon^{1/2},
\end{array}
&
\begin{array}{l}
| F'|_{h',\delta}\leq 3\tilde Me^{-2\pi(h-h')aN}\epsilon, \\ |f^{(re)}|_{h,\delta}\leq 3\tilde M\epsilon.
\end{array}
\end{matrix}
\]
 Moreover,
viewing $f^{(re)}$ as a block matrix with the same partition as $A$, we have  
$f_{ij}^{(re)}(\theta)=\hat{{f}}^{(re)}_{ij} (k_{ij})e^{2\pi \ii\la k_{ij},\theta\ra}$, where $k_{ij}$ is the resonant integer between $\Sigma_i(\lambda)$ and $\Sigma_j(\lambda)$ defined in $\mathcal H(N,aN,K)$ with $|k_{ij}|\leq mN$, and $f_{ij}^{(re)}(\theta)\equiv 0$, if $\Sigma_i$ and $\Sigma_j$ are non-resonant or $i=j$. 

\end{Proposition}
\begin{pf}
In the sequel, we  take the matrices as   block matrices with the same partition as $A$.
Since $\Sigma(\lambda)=\Sigma_1(\lambda)\cup\cdots\cup\Sigma_l(\lambda)$ satisfies $\mathcal H(N, aN, K)$ for $\lambda\in W_\delta(\Lambda)$, then the set
$\{1,\cdots,l\}$ can be divided into disjoint union:
$$ \{1,\cdots,l\} =  S_{1}\cup  S_{2} \cup \cdots \cup  S_{\tilde l},$$
where properties $(a)-(d)$ hold. Let
\begin{eqnarray*}
\lefteqn{\mathcal B_{h,\delta}^1=}\\
&\left\{ G\in C_{h,\delta}^\omega(\T^d\times \Lambda, gl(m, \C)) : \begin{array}{ll} G=(G_{ij})_{1\leq i, j\leq l}, \ \mathcal R_{aN}G=0; \hat G_{ii}(0)=0;  \hat G_{ij}(k_{ij})=0\\ \textrm{if} \ i\neq j \ \textrm{belong to the same}\ S_{l_1} \textrm{for some $1\leq l_1\leq \tilde l$} \end{array}\right\}.
\end{eqnarray*}
\begin{Claim}\label{claim-non-resonant}
We have $\mathcal B_{h,\delta}^1\subseteq \mathcal B_{h,\delta}^{(nre)}(\eta)$, with $\eta=48\tilde M\epsilon^{1/2}$.
\end{Claim}
\begin{pf}
For any $G\in \mathcal B_{h,\delta}^1$, one can check that $AG=(A_{ii}G_{ij})_{1\leq i,j\leq l}\in \mathcal B_{h,\delta}^1$, $G^+A=(G^+_{ij}A_{jj})_{1\leq i,j\leq l}\in \mathcal B_{h,\delta}^1$.
Suppose 
\[
AG-G^+A=H\in\mathcal B_{h,\delta}^1.
\]
 With the Fourier expansion of the $(i,j)$-th matrix block, we have for any $\lambda\in W_\delta(\Lambda)$, $1\leq i, j\leq l$, $|k|\leq aN$,
\[
A_{ii}(\lambda)\hat{G}_{ij}(k,\lambda)-\hat{G}_{ij}(k,\lambda)\cdot e^{2\pi \ii\langle k,\alpha\rangle}A_{jj}(\lambda)=\hat{H}_{ij}(k,\lambda).
\]
Let
\[
L_{ij}(k;\lambda)=\mathrm{Id}\otimes A_{ii}(\lambda)-e^{2\pi \ii\langle k,\alpha\rangle}(A_{jj}(\lambda))^{T}\otimes\mathrm{Id},
\]
where $\otimes$ represents Kronecker product between matrices. If we view $\hat{G}_{ij}(k,\lambda)$ and $\hat{H}_{ij}(k,\lambda)$ as vectors, then
\[
\hat{G}_{ij}(k,\lambda)=L_{ij}^{-1}(k;\lambda)\hat H_{ij}(k,\lambda).
\]

If $i\neq j$ belong to the same subset $S_{l_1}$ for some $1\leq l_1\leq \tilde l$, then for $k\neq k_{ij}$, $|k|\leq aN$, we have
\begin{eqnarray}\label{inequ-pro-2-1}
\|L_{ij}^{-1}(k;\lambda)\|
&\leq& \frac{(m^{2}\|L_{ij}(k;\lambda)\|)^{m^{2}}}{\inf_{\lambda\in W_{\delta}(\Lambda)}|\det L_{ij}(k;\lambda)|}\\
\nonumber&=&\frac{(m^{2}\|L_{ij}(k;\lambda)\|)^{m^{2}}}{\inf_{\lambda\in W_{\delta}(\Lambda)}|\textrm{Res}(\chi_{\Sigma_i}(\lambda), \chi_{\Sigma_j}(\lambda); \la k,\alpha\ra)|}\\
\nonumber&\leq& \frac{(2m^4\tilde M)^{m^2}}{  K^{-m^2}}=(2m^4\tilde M  K)^{m^2}.
\end{eqnarray}
where the inequality follows from property $(b)$ of  $\mathcal H(N, aN, K)$, that is $\Sigma_i(\lambda)$ and $\Sigma_j(\lambda)$ are $(maN,2mK^{-1} )$-resonant with the same resonant term $k_{ij}$. 

If $i\neq j$ belong to different subsets $S_{l_1}, S_{l_2}$, similarly, by property $(d)$ of $\mathcal H(N, aN, K)$, for any $|k|\leq aN$, we obtain that
\[
\|L_{ij}^{-1}(k;\lambda)\|\leq (4m^4 \tilde M K )^{m^2}.
\]

If $i=j$, by property $(a)$ of $\mathcal H(N, aN, K)$, then for $0<|k|\leq aN$, we can also get
\[
\|L_{ij}^{-1}(k;\lambda)\|\leq (4m^4 \tilde M K)^{m^2}.
\]

Since $G, H\in\mathcal B_{h,\delta}^1$, then for any $\lambda\in W_\delta(\Lambda)$ we have
\[|G_{ij}(\lambda)|_h=\sum_{|k|\leq aN}\|\hat G_{ij}(k,\lambda)\|e^{2\pi |k|h}\leq  (4m^4\tilde M K)^{m^2}| H_{ij}(\lambda)|_h, \]
and thus
\[|G|_{h,\delta}\leq m(4m^4\tilde MK)^{m^2}|H|_{h,\delta}<\frac{1}{48\tilde M\epsilon^{1/2}}|H|_{h,\delta},\]
which implies that
\[|AG-G^+A|_{h,\delta}\geq 48\tilde M\epsilon^{1/2}|G|_{h,\delta}.\]
Therefore, we have $\mathcal B_{h,\delta}^1\subseteq \mathcal B_{h,\delta}^{(nre)}(\eta)$ with $\eta=48\tilde M\epsilon^{1/2}$. 
\end{pf}

\medskip

Applying Lemma \ref{lem-implicit function}, there exist $Y\in C_{h,\delta}^\omega(\T^d\times\Lambda,gl(m,\C))$ and $F^{(re)}\in \mathcal B_{h,\delta}^{(re)}(\eta)$ such that
\[e^{-Y^+}(A+F)e^Y=A+F^{(re)},\]
with $|Y|_{h,\delta}\leq \epsilon^{1/2}$, $|F^{(re)}|_{h,\delta}\leq 3\tilde M\epsilon$. Then by Claim \ref{claim-non-resonant},  for $1\leq i\leq l$, 
 \[
 \mathcal T_{aN}F_{ii}^{(re)}(\theta)=\hat{F}_{ii}^{(re)}(0),
 \]
  and if $i\neq j$ belong to the same $S_{l_1}$, then
\[\mathcal T_{aN}F_{ij}^{(re)}(\theta)=\hat{F}_{ij}^{(re)}(k_{ij})e^{2\pi \ii\la k_{ij},\theta\ra};\]
otherwise, $\mathcal T_{aN}F_{ij}^{(re)}(\theta)\equiv 0$.

Finally we let $ F':=\mathcal R_{aN} F^{(re)}$, and
\begin{eqnarray*}A'&:=&A+\textrm{diag}\{\hat{ F}_{11}^{(re)}(0), \cdots, \hat{ F}_{ll}^{(re)}(0)\},\\
 f^{(re)}&:=&\mathcal T_{aN}F^{(re)}-\textrm{diag}\{\hat{ F}_{11}^{(re)}(0), \cdots, \hat{ F}_{ll}^{(re)}(0)\}.\end{eqnarray*}
 Then one can compute
\begin{eqnarray*}
 | F'|_{h',\delta}&=&\sup_{\lambda\in W_\delta(\Lambda)}\sum_{|k|>aN}\|\hat{F}^{(re)}(k,\lambda)\|e^{2\pi |k|h'} \\
&\leq& \sup_{\lambda\in W_\delta(\Lambda)}e^{-2\pi a N(h-h')}\sum_{|k|>aN}\|\hat{F}^{(re)}(k,\lambda)\|e^{2\pi |k|h}\\
&\leq & 3\tilde Me^{-2\pi a N(h-h')}\epsilon.
\end{eqnarray*}
We finish the whole proof.
\end{pf}
%

\subsection{Remove the resonances}\label{sub-sec-4-3}
\begin{Lemma}\label{rmr}
Let  $A'\in C_{\delta}^{\omega}(\Lambda,\mathrm{GL}(m,\C))$, $f^{(re)} \in C_{h,\delta}^{\omega}(\T^{d}\times\Lambda,gl(m,\C))$, where 
$$A'(\lambda)=\mathrm{diag}\{A_{11}'(\lambda),\cdots,A_{ll}'(\lambda)\}$$ for $\lambda\in W_{\delta}(\Lambda),$ and view $f^{(re)}$ as block matrices with the same block partition as $A'$. Suppose that there exists a disjoint partition of $\{1,\cdots,l\}$ 
 $$ \{1,\cdots,l\} =  S_{1}\cup  S_{2} \cup \cdots \cup  S_{\tilde l}, $$
 such that the following hold:
\begin{enumerate}[a)]
\item If $i, j \in S_{l_1}$ with $i\neq j$, then there exists a unique $k_{ij}$ such that $\hat {f}_{ij}^{(re)}(k)=0$ for $k\neq k_{ij}$, and
we have $k_{ij}+k_{ji}=0$.
\item If  $i,j$ and $t$ are different and belong to   the same subset $S_{l_1}$, we have $k_{ij}+k_{jt}=k_{it}$.
\item If $i=j$, then $\hat {f}_{ij}^{(re)}(k)=0$ for any $k$, and we denote $k_{ii}=0$.
\item If $i \in S_{l_1}$, $j\in S_{l_2}$ with $l_1\neq l_2$, then  $\hat {f}_{ij}^{(re)}(k)=0$ for any $k$, and we denote $k_{ij}=0$.
\end{enumerate}
Then   there exists $H\in C_{h}^{\omega}(\T^{d},\GL(m,\C))$, independent of $\lambda$, such that
\[
H^{-1}(\cdot+\alpha)(A'+f^{(re)} )H(\cdot)=A''.
\]
where $|H|_{h}, |H^{-1}|_h\leq \max_{1\leq i,j\leq l}\{e^{2\pi |k_{ij}|h}\}$, $A''\in C_{\delta}^{\omega}(\Lambda,\mathrm{GL}(m,\C))$,  $|A_{ij}''|_{\delta}=|\hat {f}_{ij}^{(re)}(k_{ij})|_{\delta}$ for  $i\neq j$, and for $1\leq i\leq l$,   $A_{ii}''=e^{-2\pi \ii\langle k_{ij},\alpha\rangle}A_{ii}'$ for some $j$ that belongs to the same subset $S_{l_1}$ as $i$.
\end{Lemma}
\begin{pf}
For each $1\leq i\leq l$ we can find $S_{m(i)}$ such that $i\in S_{m(i)}$. Let 
\[
n(i)=\min_{j\in S_{m(i)}}\{j\}.\]
 We view $H$ as a block matrix with the same block partition as $A'$, and let $H_{ii}=e^{2\pi \ii\la k_{in(i)},\theta\ra}\textrm{Id}$, $H_{ij}:=0$ for $i\neq j$. Then $H\in C^\omega_h(\T^d, \GL(m,\C))$ does not depend on $\lambda$.  Moreover, we have
$$H^{-1}(\theta+\alpha)A'H(\theta)=\textrm{diag}\{e^{-2\pi \ii\la k_{1n(1)},\alpha\ra}A_{11}',\cdots, e^{-2\pi \ii\la k_{ln(l)},\alpha\ra}A_{ll}'\}.$$ Furthermore, the $(i,j)$-th block of $H^{-1}(\theta+\alpha)f^{(re)}(\theta)H(\theta)$ is 0 if $i=j$ or $i,j$ belong to different subsets $S_{l_1}, S_{l_2}$; otherwise, we have $i\neq j$ belong to the same subset $S_{l_1}$ for some $1\leq l_1\leq \tilde l$, implying $n(i)=n(j)$, and then the $(i,j)$-th block is
$$e^{-2\pi \ii\la k_{in(i)},\alpha\ra}\hat {f}_{ij}^{(re)}(k_{ij})e^{2\pi \ii\la k_{ij}+k_{jn(j)}-k_{in(i)},\theta\ra},$$ where
\[
k_{ij}+k_{jn(j)}-k_{in(i)}=
\left\{\begin{array}{lll} k_{ij}+k_{ji}=0,& \textrm{if}\ n(j)=n(i)=i,\\
k_{ij}-k_{ij}=0, &\textrm{if}\ n(j)=n(i)=j,\\
k_{ij}+k_{jn(j)}-k_{in(j)}=0, & \textrm{otherwise},
\end{array}\right.
\]
and thus equals $e^{-2\pi \ii\la k_{in(i)},\alpha\ra}\hat {f}_{ij}^{(re)}(k_{ij})$.
Let 
\[
A'':=H^{-1}(\theta+\alpha)(A'+f^{(re)}(\theta))H(\theta).
\]
 and we finish the proof.

%
%
%
\end{pf}

\subsection{Block diagonalization}\label{sub-sec-4-4}

The following lemma assures that if the spectrum of the diagonal block elements of the new constant matrix is separated again, then the new constant matrix can be diagonalized by a matrix close to identity.  
\begin{Lemma}\label{norm2}
Let $R, \tilde M'>1$. Given $A, G$ where $A\in C_{\delta}^{\omega}(\Lambda,\mathrm{GL}(m,\C)), G\in C_\delta^\omega(\Lambda, {gl}(m,\C))$ with $A(\lambda)=\mathrm{diag}\{A_{11}(\lambda),\cdots,A_{ll}(\lambda)\}$ for $\lambda\in W_{\delta}(\Lambda)$, and $A_{ii} (1\leq i\leq l)$ being block matrices, suppose $|A|_{\delta}\leq\tilde{M}'$, $|G|_{\delta}\leq\epsilon$, the decomposition $\Sigma(\lambda)=\Sigma_{1}(\lambda)\cup\cdots\cup\Sigma_{l}(\lambda)$ is $(\nu,\zeta)$-separated in $W_{\delta}(\Lambda)$ with $0<\nu\leq 1$,  $(M,\delta,c,r)$-transverse on $\Lambda$, and $\Sigma(\lambda)\subseteq D(R)$ for $\lambda\in W_{\delta}(\Lambda)$, where $\Sigma(\lambda):=\Sigma(A(\lambda))$ and $\Sigma_i(\lambda):=\Sigma(A_{ii}(\lambda))$ for $1\leq i\leq l$. If
\begin{equation}\label{equ-lem-10-1}
(6m^{4}\tilde{M}'\nu^{-1})^{m^{2}+1}\epsilon^{\frac{1}{2}}\leq 1,
\end{equation}
and $2^7m^2\tilde M'^3R\epsilon^{1/m}<1$, then there exists $S\in C_{\delta}^{\omega}(\Lambda,\mathrm{GL}(m,\C))$ such that
\begin{align}\label{n2eq}
S^{-1}(A+G)S=\tilde{A},
\end{align}
where $\tilde{A}\in C_{\delta}^{\omega}(\Lambda,\mathrm{GL}(m,\C))$, $\tilde{A}=\mathrm{diag}\{\tilde{A}_{11},\cdots,\tilde{A}_{ll}\}$ in $W_{\delta}(\Lambda)$,
\begin{align*}
|A_{ii}-\tilde{A}_{ii}|_{\delta}&\leq 3\tilde M'\epsilon,\\
|S-\mathrm{Id}|_{\delta}&\leq 2\epsilon^{\frac{1}{2}},
\end{align*}
and the decomposition $\tilde\Sigma(\lambda)=\tilde\Sigma_1(\lambda)\cup\cdots\cup\tilde\Sigma_l(\lambda)$ is $(\nu',\zeta')$-separated in $W_\delta(\Lambda)$, with $\tilde\Sigma(\lambda):=\Sigma(\tilde A(\lambda))$ and $\tilde\Sigma_i(\lambda):=\Sigma(\tilde A_{ii}(\lambda))$ for $1\leq i\leq l$.
In addition, we have $\tilde\Sigma(\lambda)=\tilde\Sigma_{1}(\lambda)\cup\cdots\cup\tilde\Sigma_{l}(\lambda)$ is $(M',\delta,c',r)$-transverse on $\Lambda$, $\tilde\Sigma(\lambda)\subseteq D(R')$ for $\lambda\in W_{\delta}(\Lambda)$, and for any $1\leq i\leq l$, $\lambda\in W_\delta(\Lambda)$,
\[
d_H(\Sigma_i(\lambda), \tilde\Sigma_i(\lambda))<8m^2\tilde M'^3\epsilon^{\frac{1}{m}},
\]
where
$$
\begin{matrix}
\begin{array}{l}
R'=R+16m^2 R^{2}\tilde{M}'^3\epsilon^{\frac{1}{m}},\\ M'=(2R')^{m^{2}}, \\ c'=c-2^{8m+2}m^{10m^2}R^{3m^2} \tilde M'^{m+1} (\frac{2r}{\delta})^{r}\epsilon.
\end{array}
&
\begin{array}{l}
\nu'=\nu-16m^2  \tilde{M}'^3\epsilon^{\frac{1}{m}} \\ \zeta'=\zeta+16m^2  \tilde{M}'^3\epsilon^{\frac{1}{m}},\\
{}
\end{array}
\end{matrix}
$$
\end{Lemma}
\begin{pf}
Let $\mathcal B_\delta:=C^\omega_\delta(\Lambda, gl(m,\C))$, and it can be viewed as a special case of $\mathcal B_{h,\delta}$ in Section \ref{sub-sec-4-2}, where the functions are constant for the variable in $\T^d$ and we omit $h$ here. Let 
\[\mathcal B_{\delta}^1:=\{Y\in C_\delta^\omega(\Lambda, gl(m,\C))\ :\ Y_{ii}=0,\ \forall 1\leq i\leq l\},\]
\[\mathcal B_{\delta}^2:=\{Y\in C_\delta^\omega(\Lambda, gl(m,\C))\ :\ Y_{ij}=0,  \ \forall 1\leq i\neq j \leq l\}.\]
Then $B_\delta=\mathcal B_\delta^1\oplus\mathcal B_\delta^2$, and for any $Y\in\mathcal B_\delta^1$, we have the $(i,j)$-th block matrices of $AY$ and $YA$ are $A_{ii}Y_{ij}$ and $Y_{ij}A_{jj}$ respectively. It is obvious that $AY, YA\in\mathcal B_\delta^1$. Moreover, if we let 
\[
AY-YA=H\in \mathcal B_\delta^1,
\]
 then for $i\neq j$, the $(i,j)$-th block matrix equation is 
\[A_{ii}Y_{ij}-Y_{ij}A_{jj}=H_{ij}.\]
Let 
\[L_{ij}(\lambda)=\textrm{Id}\otimes A_{ii}(\lambda)-(A_{jj}(\lambda))^T\otimes\textrm{Id}.\]
Then for any $\lambda\in W_\delta(\Lambda)$ and $i\neq j$, if we view $Y_{ij}$ and $H_{ij}$ as vectors, we have 
\[Y_{ij}(\lambda)=L_{ij}^{-1}(\lambda)H_{ij}(\lambda).\]
By the fact that $\mathrm{dist}(\Sigma_{i}(\lambda),\Sigma_{j}(\lambda))>\nu$ for $i\neq j$ and $\lambda\in W_{\delta}(\Lambda)$, we have
 \[|\det L_{ij}(\lambda)|=|\textrm{Res}(\chi_{\Sigma_i(\lambda)},\chi_{\Sigma_j(\lambda)})|\geq \nu^{m^2},\]
 which implies that
 \[|L^{-1}_{ij}|_\delta\leq (2m^4\nu^{-1}\tilde M')^{m^2}.\]
Therefore, we have
\[
|Y|_{\delta}\leq m (2m^4\nu^{-1}\tilde M')^{m^2}|H|_{\delta}\leq \frac{|H|_{\delta}}{48\tilde M' \epsilon^{1/2}}.
\]
which means $|AY-YA|_\delta\geq 48\tilde M'\epsilon^{1/2}|Y|_\delta$, and thus $\mathcal B_\delta^1\subseteq \mathcal B_\delta^{(nre)}(\eta)$ with $\eta=48\tilde M'\epsilon^{1/2}$. Then by Lemma \ref{lem-implicit function}, there exist $Y\in C_\delta^\omega(\Lambda, gl(m,\C))$ and $G^{(re)}\in\mathcal B_\delta^{(re)}(\eta)\subseteq \mathcal B_\delta^2$ with $|Y|_\delta\leq \epsilon^{1/2}$ and $|G^{(re)}|_\delta\leq 3\tilde M'\epsilon$, such that 
\[e^{-Y}(A+G)e^Y=A+G^{(re)}.\]
Denote $\tilde A=A+G^{(re)}$ and $S=e^Y$. Then $\tilde A=\textrm{diag}\{\tilde A_{11},\cdots,\tilde A_{ll}\}$, and 
\[|A_{ii}-\tilde A_{ii}|_\delta=|G_{ii}^{(re)}|_\delta\leq 3\tilde M'\epsilon,\]
\[|S-\textrm{Id}|_\delta=|e^Y-\textrm{Id}|_\delta<2\epsilon^{1/2}.\]
Since $|A_{ii}-\tilde A_{ii}|_\delta\leq 3\tilde M'\epsilon$, then by Lemma \ref{ndh}, 
 the rest of the lemma holds, and we finish the proof.
\end{pf}

\section{Stratified quantitative almost reducibility and reducibility }\label{sec-4}
We will use the KAM steps developed in Section \ref{seckam} to prove the stratified quantitative almost reducibility for all parameters, as well as the quantitative reducibility for full measure set of parameters. 

Let $\Lambda\subseteq \R$ be an interval, $\alpha\in\mathrm{DC}_{d}(\gamma,\tau)$ with $\gamma>0, \tau>d$,  $\delta_1, c_1>0$, $r_1\in \N^+$, $M_1, R_1>1, R=2R_1$, and $\Pi_1$ be a partition of $\Lambda$. For each $\tilde\Lambda \in \Pi_1$ and $A_1\in C^\omega_{\delta_1}(\tilde\Lambda, \GL(m,\C))$, 
we
suppose $\Sigma(A_1(\lambda))\subseteq D(R_1)$ for $\lambda\in\tilde\Lambda$, and $\Sigma(A_{1}(\lambda))$ is $(M_{1},\delta_{1},c_{1}, r_{1})$-transverse on $\tilde\Lambda$.

\subsection{Stratified quantitative almost reducibility} In this subsection, we state the iterative lemma (Proposition \ref{pro-iterative}) for the stratified quantitative almost reducibility. We first give the iteration parameters.
 Let $0<\varepsilon<1$, $\epsilon_1=\varepsilon^{1/2}$, and for $n\geq 1$,
\begin{align}
\epsilon_{n+1}=\epsilon_{n}^{4^{n}}=\epsilon_{1}^{2^{n^{2}+n}}, \quad h_{n}=(\frac{1}{2}+\frac{1}{2^{n}})h_{1}.
\end{align}
We now define
\begin{align*}
N_{n}=[2^{n+1}\frac{|\log\epsilon_{n}|}{2\pi h_{1}}]+1, \ \ a_{n}=8^{n+1}m,
\end{align*}
\begin{align*}
N_{n,p}=a_{n}^{p}N_{n},\ \ K_{n}=144mR\gamma^{-1}(3mN_{n,m+1})^{\tau}.
\end{align*}
Let $b_{1}(m,R,\tau,\gamma)=160m^{(m+2)(\tau+1)}12^{\tau+2}R\gamma^{-1}$, $b_{2}(m,\tau)=(6m+7)\tau$ and 
\begin{equation}\label{equ-def-u-n}
u_{n}=b_{1}^{-1}(\frac{4\pi h_{1}}{|\log\varepsilon|})^{\tau}e^{-b_{2}e^{4\sqrt{n}}}.
\end{equation}
Set $s_{1}=1$, and we define $s_{n}$ iteratively:
\begin{align}\label{pset}
s_{i+1}=\min\{n: K_{n}^{-1}<d_1u_{s_{i}}\},
\end{align}
where $d_1=160m^{m+1}$.
We denote 
\[
\xi_{n}=\#\{i: s_{i}\leq n\}.
\]
 By \eqref{pset} we have
$$
u_{s_{i}}^{-1}<d_1K_{s_{i+1}}\leq b_{1}(\frac{|\log\epsilon_{1}|}{2\pi h_{1}})^{\tau}2^{b_{2}s_{i+1}^{2}}=b_{1}(\frac{|\log\varepsilon |}{4\pi h_{1}})^{\tau}2^{b_{2}s_{i+1}^{2}}.
$$
Then we can get
\[
\sqrt{s_{i+1}}>e^{\sqrt{s_{i}}}.
\]
Similarly, by the fact $d_1K_{s_{i+1}-1}\leq u_{s_i}^{-1}$, we obtain that
\[
\sqrt{s_{i+1}}\leq d_2e^{\sqrt{s_{i}}},
\]
where $d_2=d_2(m)$ is a constant.
Thus for any $k\in\N$ we have
\[\xi_n<2\ln^{(k)}(n)+k+1,\]
where $\ln^{(k)}(n)$ means $\underbrace{\ln\cdots\ln}_{k}(n)$, and if $\xi_{n}\leq k$ we denote $\ln^{(k)}(n)=1$.

Suppose $\epsilon_1$ is small enough such that 
\begin{equation}\label{equ-epsilon-add-1}
(\frac{1}{2\pi h_1}\ln \frac{1}{\epsilon_1})^{6m(m+2)\tau}<\frac{1}{\epsilon_1},\ \ 36(m+7)^3\tau<\ln \frac{1}{\epsilon_1}.
\end{equation}
Moreover, it is obvious that for $\epsilon_1$ sufficiently small, we have for all $k\in\N$
\begin{equation}\label{equ-epsilon-add-2}
6m(m+2)b_2e^{4\sqrt s_k}<2^{s_k^2-3s_k+2}|\ln \epsilon_1|,  
\end{equation} 
In the sequel, we assume $\epsilon_1$ is small enough such that (\ref{equ-epsilon-add-1}) and (\ref{equ-epsilon-add-2}) hold.


\medskip


Suppose $\Pi_{n} (n\geq 1)$ is a partition of $\Lambda$. We denote by $b,\kappa$ the constants only depending on $m$ with $b\geq \max\{(120m)^{8m^3}, \acute b\}$ and $\kappa\geq m^{2m^2+ 10}$, where $\acute b$ is the constant only depending on $m$ in Lemma \ref{lem-similar-transform}. For all $n\geq 2$, we introduce the following notation for convenience.

\medskip

\noindent\textbf{Property $\mathcal P(n,\Pi_n,\Pi_{n-1})$}: For matrices-valued functions $A_n, A_{n-1}$ on $\Lambda$, on each $\tilde \Lambda\in\Pi_n$, where $\Pi_{n}$ is a sub-partition of $\Pi_{n-1}$, for sequences of numbers $R_n, r_n, \nu_n, \zeta_n, \delta_n, \tilde M_n, c_n$ defined according to $\tilde\Lambda$,  the following hold:
\begin{enumerate}[$(\mathcal P 1)_n$]
\item We have $A_n \in C_{\delta_{n}}^{\omega}(\tilde{\Lambda},\mathrm{GL}(m,\C))$,  $|A_n|_{\delta_{n}}<\tilde{M}_{n}$, and for any $\lambda\in W_{\delta_{n}}(\tilde\Lambda)$,  
\[
    A_n (\lambda)=\mathrm{diag}\{A_{n,11}(\lambda),\cdots,A_{n, l_{n}l_{n}}(\lambda)\},
\]
    where $\Sigma^{(n)}(\lambda):=\Sigma(A_n(\lambda))\subseteq D(R_n)$, and $\Sigma_{i}^{(n)}(\lambda):=\Sigma(A_{n,ii}(\lambda))$.
\item
The decomposition $
\Sigma^{(n)}(\lambda)=\Sigma_{1}^{(n)}(\lambda)\cup\cdots\cup\Sigma_{l_{n}}^{(n)} (\lambda)$ is $(\nu_{n},\zeta_{n})$-separated in $W_{\delta_{n}}(\tilde{\Lambda})$ and $((2R_{n})^{m^{2}},\delta_{n},c_{n},r_{n})$-transverse on $\tilde{\Lambda}$.

\item For $n=s_{k}$  there exists $\lambda_{0}\in\tilde\Lambda$ such that  $\Sigma_{1}^{(n)}(\lambda_{0})\cup\cdots\cup\Sigma_{l_{n}}^{(n)}(\lambda_{0})$ is $(8\nu_{n}, \zeta_{n})$-separated and 
\begin{eqnarray*}
&d_H(\Sigma_{i}^{(n)}(\lambda),\Sigma_{i}^{(n)}(\lambda_{0}))<\nu_n,\ \forall\lambda\in W_{\delta_n}(\tilde \Lambda),\\
&\#\Pi_n \leq \#\Pi_{n-1}+\frac{bR_{n-1}^{6m}\tilde M_{n-1}^{3m}| \Lambda|}{\nu_n^{3m}\delta_{n-1}}\max_{\tilde\Lambda^-\in\Pi_{n-1}}\frac{1}{\delta_{n-1}(\tilde\Lambda^-)},
\end{eqnarray*}
 where $\delta_{n-1}$ may be different for different $\tilde\Lambda^-\in\Pi_{n-1}$, and  
\begin{eqnarray*}
&& R_n =R_{n-1}+  b R_{n-1}^2\tilde M_{n-1}^4\epsilon_{n-1}^{\frac{1}{m}},\\  
&& \nu_{n}  =u_{s_{k}},\ \   \zeta_{n} =10mu_{s_{k}},\ \   r_{n} =m^{2}r_{n-1}, \\
&&   \delta_{n} = {b}^{-1}R_{n-1}^{-6m}(\nu_{n}\tilde M_{n-1}^{-1})^{3m}\delta_{n-1}, \ \
  \tilde{M}_{n} =b(\nu_{n}^{-1}\tilde{M}_{n-1})^{m^2(m+2)},\\
&& c_{n} =( {b}^{-1}R_{n-1}^{-1}r_{n-1}^{-1}\delta_{n-1}\nu_{n}\tilde M_{n-1}^{-1}c_{n-1})^{\kappa^{r_{n-1}}}- r_{n-1}^{\kappa r_{n-1}}\epsilon_{n-1}.
\end{eqnarray*}

\item For $s_{k}<n<s_{k+1}$, we have for any $\tilde\Lambda^-\in\Pi_{n-1}$, $\lambda\in W_{\delta_{n-1}}(\tilde\Lambda^-)$
\[A_{n-1}(\lambda)=\textrm{diag}\{A_{n-1,11}(\lambda),\cdots,A_{n-1,l_{n-1}l_{n-1}}(\lambda) \},\]
and the following:
\begin{enumerate}[  \hspace{-3em} $(\mathcal P 4.1)_n$]
\item
If $l_{n}<l_{n-1}$, then there exists a partition $\Pi_n(\tilde\Lambda^-)$ of 
$\tilde\Lambda^-$ with 
\[
 \#\Pi_n(\tilde\Lambda^-)\leq 1+\frac{bR_{n-1}^{2m-1}|\tilde\Lambda^-|}{\nu_n^m\delta_{n-1}(\tilde\Lambda^-)},
\]
where 
$\delta_{n-1}$ may change with $\tilde\Lambda^-\in\Pi_{n-1}$, and for any $\tilde\Lambda\in \Pi_n(\tilde\Lambda^-)$ there exists $\lambda_{0}\in\tilde{\Lambda}$ such that  for any $\lambda\in W_{\delta_n}(\tilde\Lambda)$
\[
d_H(\Sigma_{i}^{(n)}(\lambda),\Sigma_{i}^{(n)}(\lambda_{0}))<\nu_n,
\]
with
\begin{eqnarray*}
 &&R_n=R_{n-1}+  b R_{n-1}^2\tilde M_{n-1}^4\epsilon_{n-1}^{\frac{1}{m}},\\
&&  \nu_{n} =\nu_{n-1}-b\tilde{M}_{n-1}^4\epsilon_{n-1}^{\frac{1}{m}},\ \
\zeta_{n} =m(\zeta_{n-1}+7\nu_{n-1}),\\
&& r_{n} =m^{2}r_{n-1},\ \ \delta_{n} = {b}^{-1}R_{n-1}^{-(2m-1)} \nu_{n}^{m}\delta_{n-1},\ \
\tilde{M}_{n} = \tilde{M}_{n-1}+20m\tilde M_{n-1}^2\epsilon_{n-1},\\
&& c_{n} =( {b}^{-1}R_{n-1}^{-1}r_{n-1}^{-1}\delta_{n-1} c_{n-1})^{\kappa^{r_{n-1}}}- ( bR_{n-1}  \tilde{M}_{n-1}\delta_{n-1}^{-1}r_{n-1})^{\kappa r_{n-1}}\epsilon_{n-1}.
\end{eqnarray*}

\item If $l_{n}=l_{n-1}$, then $\#\Pi_n(\tilde\Lambda^-)=1$ and there exists  $k_i^{(n)}\in \Z^d (i=1,2,\cdots, l_n)$ such that for any $\lambda\in W_{\delta_n}(\tilde\Lambda^-)$ we have
\[
d_H(\Sigma_i^{(n)}(\lambda), e^{2\pi \ii\la k_i^{(n)},\alpha\ra}\Sigma_i^{(n-1)}(\lambda))<\frac{\nu_{n-1}-\nu_n}{2},
\]
with
\begin{align*}
&R_n=R_{n-1}+  b R_{n-1}^2\tilde M_{n-1}^4\epsilon_{n-1}^{\frac{1}{m}},\\
&\nu_{n} =\nu_{n-1}-b \tilde{M}_{n-1}^4\epsilon_{n-1}^{\frac{1}{m}},\ \
\zeta_{n} =\zeta_{n-1}+b \tilde{M}_{n-1}^4\epsilon_{n-1}^{\frac{1}{m}},\\
&r_{n}=r_{n-1},\ \ \delta_{n} =\delta_{n-1},\ \
\tilde{M}_{n} =\tilde{M}_{n-1}+20m\tilde M_{n-1}^2\epsilon_{n-1},\\
& c_{n} =c_{n-1}-b  R_{n-1}^{3m^2} \tilde{M}_{n-1}^{m+2} (2 r_{n-1} \delta_{n-1}^{-1})^{ r_{n-1}}\epsilon_{n-1}.
\end{align*}
\item We have 
\[
\#\Pi_n\leq \#\Pi_{n-1}+\frac{bR_{n-1}^{2m-1}|\Lambda|}{\nu_n^m}\max_{\tilde\Lambda\in\Pi_{n-1}}\frac{1}{\delta_{n-1}(\tilde\Lambda)},
\]
 where $\Pi_n=\cup_{\tilde\Lambda^-\in\Pi_{n-1}}\Pi_n(\tilde\Lambda^-)$. 
\end{enumerate}
\end{enumerate}

\medskip

We say $\mathfrak P(1)$ holds if there exists some  partition $\Pi_1$ of $\Lambda$, that the properties $(\mathcal P 1)_1$-$(\mathcal P 3)_1$  hold for $A_1$. For $\tilde m\geq 2$, we say $\mathfrak P(\tilde m)$ holds if $\mathfrak P(1)$ holds and for any $2\leq n\leq \tilde m$, there exists a partition $\Pi_n$ of $\Lambda$ such that $\mathcal P(n,\Pi_n,\Pi_{n-1})$ holds for $A_n, A_{n-1}$.
  
 \medskip
  
Now we are in the position to state the iterative lemma for the stratified quantitative almost reducibility.

\begin{Proposition}[Iterative Lemma] \label{pro-iterative}
Assuming $\mathfrak P(\tilde m)$ holds for any $\tilde m\geq 1$,  $|F_{\tilde m}|_{h_{\tilde m},\delta_{\tilde m}}<\epsilon_{\tilde m}$, 
 \begin{eqnarray}
\label{equ-pro-3-con-2}& K_{\tilde m}^{-1}<\frac{\gamma}{144mR_{\tilde m} (3mN_{\tilde m,m+1})^\tau},\\
\label{equ-pro-3-con-3} & \epsilon_{\tilde m}^{\frac{1}{2m}}<\frac{\gamma}{48R_{\tilde m}(2N_{\tilde m,m+1})^\tau},
\end{eqnarray}
and 
\begin{align}\label{con1}
C(d)(  b R_n\tilde{M}_{n}K_{n}\nu_{n}^{-1})^{ m^{2} +1}(\frac{ 3\cdot 2^nN_{n,m+2}  }{h_1 })^{ d}(\frac{2m^2r_n}{\delta_n})^{m^2r_n/3} c_n^{-1/3}\leq \epsilon_{n}^{-\frac{1}{6}},
\end{align}
  for all $n\leq \tilde m$,
then there exists a sub-partition $\Pi_{\tilde m+1}$  of $\Pi_{\tilde m }$,  such that on each $\tilde\Lambda\in\Pi_{\tilde m+1}$, there exists
 $Z_{\tilde m}\in C_{h_{\tilde m+1},\delta_{\tilde m+1}}^{\omega}(\T^{d}\times\tilde\Lambda,\mathrm{GL}(m,\C))$ such that
 \[
    Z_{\tilde m}^{-1}(\cdot+\alpha)(A_{\tilde m}+F_{\tilde m}(\cdot))Z_{\tilde m}(\cdot)=A_{\tilde m+1}+F_{\tilde m+1}(\cdot),
\]
  satisfying $\mathfrak P(\tilde m+1)$. Moreover, we have estimates:
    \begin{itemize}
    \item  $  |F_{\tilde m+1}|_{h_{\tilde m+1},\delta_{\tilde m+1}}<\epsilon_{\tilde m+1}  $, $|F_{\tilde m+1}|_{h_{\tilde m+1},\delta_{\tilde m+1}}\leq |F_{\tilde m }|_{h_{\tilde m },\delta_{\tilde m }}$,
    \item 
  $|Z_{\tilde m}|_{h_{\tilde m+1},\delta_{\tilde m+1}},\ \ |Z_{\tilde m}^{-1}|_{h_{\tilde m+1},\delta_{\tilde m+1}}<\epsilon_{\tilde m+2m+2}^{-1},$
  \item  
  $ |Z_{\tilde m}|^{ \tilde m^2/2}_{h_{\tilde m+1},\delta_{\tilde m+1}} |F_{\tilde m+1}|_{h_{\tilde m+1},\delta_{\tilde m+1}}, \ \  |Z_{\tilde m}^{-1}|^{ \tilde m^2/2}_{h_{\tilde m+1},\delta_{\tilde m+1}}|F_{\tilde m+1}|_{h_{\tilde m+1},\delta_{\tilde m+1}} <1.$ 
     \end{itemize}
\end{Proposition}

We leave the proof to Section \ref{sec-5-3}.

\subsection{Stratified quantitative reducibility}
Now we state the key proposition (Proposition \ref{pro-iterative-lemma-full}) to prove stratified full measure reducibility (Theorem \ref{Thm2}), indeed reducibility out of a zero Hausdorff dimensional set. We first give the iteration parameters.

Let $0<\tilde \epsilon_1, \tilde h_1, \tilde \delta_1, \tilde c_1'<1$, $\tilde M_1', \tilde R_1, \tilde C_1>1$, $\tilde r_1\in\N^+$, $\tilde r_1'=m^2\tilde r_1$, and  
define  sequences inductively for $j\geq 1$:
\begin{eqnarray*}
&&\tilde\epsilon_{j+1}=\tilde\epsilon_1^{(\frac{3}{2})^{j}}, \  \ \tilde h_{j+1}=\tilde h_j-\frac{\tilde h_1}{2^{j+1}}, \ \ \tilde N_j=2^{j }\frac{|\ln \tilde\epsilon_j|}{\pi\tilde h_1},\ \  K_j=\tilde\epsilon_j^{-\frac{1}{10\tilde r_1'}},\\
&& \tilde M_{j+1}'=\tilde M_j'+3\tilde M_j'\tilde\epsilon_j,\ \ \tilde R_{j+1}=\tilde R_j  +16m^2\tilde M_j'^3\tilde R_j^2\tilde\epsilon_j^{\frac{1}{m}},  \ \ \tilde\delta_{j+1}=\frac{\tilde\delta_j}{4(2\tilde R_j)^{m^2}K_j}, \\
&& \tilde C_{j+1}=\tilde C_j+\frac{(\tilde r_1'+1)!}{\tilde\delta_{j+1}^{\tilde r_1'+1}}(8m^{4}\tilde R_j \tilde M_j')^{ 3m^2 }\tilde\epsilon_j, \ \ \tilde c_{j+1}'=\tilde c_j'-\frac{ \tilde r_1'  !}{(\tilde\delta_{j+1}/2)^{\tilde r_1' }}(8m^{4}\tilde R_j \tilde M_j')^{ 3m^2 }\tilde\epsilon_j.
 \end{eqnarray*}

 \begin{Proposition}\label{pro-iterative-lemma-full}
Let $0<\varrho, \gamma, \tilde \epsilon_1, \tilde h_1, \tilde \delta_1, \tilde c_1'<1$, $\tau>d$, $\tilde M_1', \tilde R_1, \tilde C_1>1$, $\tilde r_1\in\N^+$, $\tilde r_1'=m^2\tilde r_1$, $\tilde \Lambda\subseteq \R$ is an interval, and $\tilde\epsilon_j, \tilde h_j, \tilde N_j, K_j, \tilde M_j', \tilde R_j, \tilde\delta_j, \tilde C_j, \tilde c_j'$ are defined as above. If $\tilde\epsilon_1, \tilde M_1', \tilde R_1, \tilde\delta_1, \tilde r_1, \tilde h_1, \tilde\delta_1$ satisfy
\begin{eqnarray}
 \label{equ-est-epsilon-1}  \tilde \epsilon_1^{\frac{\varrho}{100 \tilde r_1'^2}}
  &<& \min\{  \ (9\tilde R_1)^{-(\tau+d)},\  \tilde c_1',\  (2|\ln\tilde \epsilon_1|)^{-(\tau+d)}\},\\
\label{equ-est-epsilon-2}  \tilde \epsilon_1 &<&  (\gamma\tilde h_1^\tau)^{4m^2} (16m^2\tilde R_1\tilde M_1'\tilde r_1'\tilde\delta_1^{-1})^{-2^{8}m^2\tilde r_1'},
 \end{eqnarray}
 then the following holds for $j\geq 1$:  
 
 Suppose there exists a disjoint union of intervals $\Lambda^{(j)}\subseteq \tilde\Lambda$, such that for any $\bar\Lambda^{(j)}\in\mathcal C(\Lambda^{(j)})$, where  $\mathcal C(\Lambda^{(j)})$ denotes the set of all the connected components of $\Lambda^{(j)}$,  we have  gotten 
$\tilde A_j\in C^\omega_{\tilde\delta_j}(\bar\Lambda^{(j)},\mathrm{GL}(m,\C))$, $\tilde F_j\in C^\omega_{\tilde h_j, \tilde\delta_j}(\T^d\times \bar\Lambda^{(j)}, gl(m,\C))$, with $|\tilde A_j|_{\tilde \delta_j}\leq\tilde M_j'$, $|\tilde F_j|_{\tilde h_j, \tilde\delta_j}\leq \tilde\epsilon_j$,  $\tilde\Sigma_j(\lambda):=\Sigma(\tilde A_j(\lambda))$ being $((2\tilde R_j)^{m^2}, \tilde\delta_j, \tilde c_j', \tilde r_1')$-transverse on $\bar \Lambda^{(j)}$,  and $\tilde\Sigma_j(\lambda)\subseteq D(\tilde R_j)$ for $\lambda\in W_{\tilde\delta_j}(\bar\Lambda^{(j)})$.

Furthermore, we assume that for any $\lambda\in\bar\Lambda^{(j)}$, $u\in\T$, we have
\[\sup_{0\leq l\leq \tilde r_1'+1}|\frac{\partial^l g_j(\lambda,u)}{\partial\lambda^l}|\leq \tilde C_j,\]
where $g_j(\lambda,u)=\prod_{\sigma_{\ell_1},\sigma_{\ell_2}\in\tilde\Sigma_j(\lambda) \atop{\ell_1\neq \ell_2}}(\sigma_{\ell_1}-e^{2\pi\ii u}\sigma_{\ell_2})$. Then there exists 
$\Lambda^{(j+1)}\subseteq \Lambda^{(j)}$, such that 
for any $\bar\Lambda^{(j+1)}\in\mathcal C(\Lambda^{(j+1)})$, 
there are $\tilde A_{j+1}\in C^\omega_{\tilde\delta_{j+1}}(\bar\Lambda^{(j+1)}, \GL(m,\C))$, $Y_j, \tilde F_{j+1}\in C_{\tilde h_j,\tilde \delta_{j+1}}^\omega(\T^d\times\bar\Lambda^{(j+1)}, gl(m,\C))$, with $|Y_j|_{\tilde h_j,\tilde\delta_{j+1}}\leq \tilde\epsilon_j^{\frac{1}{2}}$, such that 
\[e^{-Y_j(\cdot+\alpha)}(\tilde A_j+\tilde F_j(\cdot))e^{Y_j(\cdot)}=\tilde A_{j+1}+\tilde F_{j+1}(\cdot),\]
$|A_{j+1}-A_j|_{\tilde \delta_{j+1}}\leq 3\tilde M_j'\tilde\epsilon_j$,
and $\tilde A_{j+1}, \tilde F_{j+1}$ satisfy the same assumptions of $\tilde A_j, \tilde F_j$ with $j$ being replaced by $j+1$,
where 
\[
\Lambda^{(j+1)}=\Lambda^{(j)}\backslash\mathcal R_j(\tilde\Lambda),\ \ 
\mathcal R_j(\tilde\Lambda)=\cup_{i\in J_j(\tilde\Lambda)} I^{(j)}_i,
\]
with $\{I_{i}^{(j)}\}_i\subseteq \Lambda^{(j)}$ being disjoint intervals, $|I_{i}^{(j)}|<\tilde\epsilon_j^{\frac{2}{25\tilde r_1'^2}}$, and   
\begin{eqnarray*}
&&\#J_j(\tilde\Lambda)\leq 2^{\tilde r_1' +d}\tilde N_j^d(\frac{8\tilde C_1}{\tilde c_1'}| \Lambda^{(j)}|+\#\mathcal C (\Lambda^{(j)})),\\
 && \#\mathcal C(\Lambda^{(j+1)})\leq 2\#\mathcal C(\mathcal R^{(j)}(\tilde\Lambda)).
\end{eqnarray*}
  \end{Proposition}

\subsection{Proof of Proposition \ref{pro-iterative}}\label{sec-5-3}

By the assumption, we have obtained $\Pi_{\tilde m}$ and
$A_{\tilde m}, F_{\tilde m}$,
where for any $\tilde \Lambda\in\Pi_{\tilde m}$, we have $A_{\tilde m}=\mathrm{diag}\{A_{\tilde m,11},\cdots, A_{\tilde m,l_{\tilde m}l_{\tilde m}}\}$ in $W_{\delta_{\tilde m}}(\tilde{\Lambda})$. In the sequel, we fix $\tilde\Lambda\in\Pi_{\m}$. Without loss of generality, we assume $|F_{\tilde m}|_{h_{\tilde m},\delta_{\tilde m}}\geq \epsilon_{\tilde m+1}$; otherwise, we let $A_{\tilde m+1}=A_{\tilde m}$, $F_{\tilde m+1}=F_{\tilde m}$, and if necessary, we divide $\tilde\Lambda$ into $1+\left[\frac{\tilde\Lambda}{\delta_{\m+1}}\right]$ intervals with equal length less than $\delta_{\m+1}$ as the partition $\Pi_{\m+1}(\tilde\Lambda)$ of $\tilde\Lambda$. Then the result follows.

Without loss of generality, we assume that $s_k\leq \tilde m<s_{k+1}$.
\subsubsection{Step 1: Eliminate the non-resonant terms}
\begin{Claim}
There exists $p\in\{0,1,\cdots,m\}$ such that the decomposition $\Sigma_{1}^{(\tilde m)}\cup\cdots\cup\Sigma_{l_{\tilde m}}^{(\tilde m)}$ satisfies $\mathcal{H}( N_{\tilde m,p},N_{\tilde m,p+1}, K_{\tilde m} )$ in $W_{\delta_{\tilde m}}(\tilde \Lambda)$.  
\end{Claim}
\begin{pf}
Let  
\[m_{1}=\max\{n\in \N\cap [s_{k}, \tilde m) \ :\ l_{n}>l_{\tilde m}\},\]
and if no such $n$ exists, let $m_{1}=s_{k}$. Then by properties $ (\mathcal P3)_{s_k}$ and $ (\mathcal P4)_n$ with $m_1+1\leq n\leq \tilde m$, we have
\begin{equation}\label{equ-pro-3-5}
d_H(\Sigma_{i}^{(\tilde m)}(\lambda),  \Sigma_{i}^{(\tilde m)}(\lambda_{0}))<2\nu_{m_{1}}-\nu_\m,
\end{equation}
 for some $\lambda_{0}\in\tilde{\Lambda}$.
  By the definition of $\nu_n$ and $\zeta_n$ for $s_k\leq n\leq \tilde m$ and (\ref{con1}), we obtain that
\begin{equation}\label{equ-est-mu-zeta}
\frac{u_{s_k}}{2}<\nu_n\leq u_{s_k}, \ \ 10mu_{s_k}\leq \zeta_n<20 m^{m+1}u_{s_k},
\end{equation}
where $u_{s_k}$ is defined as in (\ref{equ-def-u-n}). Then by (\ref{pset}) and (\ref{equ-pro-3-con-2}),
 we have
\[
16m\nu_{m_1}<8m\zeta_{\tilde m}< mK_{\tilde m}^{-1}<\frac{\gamma}{10R_{\tilde m}(3mN_{\tilde m,m+1})^{\tau}}.
\]
Together with $(\mathcal P1)_{\tilde m}$-$( \mathcal P2)_{\m}$,  we can apply Lemma \ref{deco}, and then the result follows. 
\end{pf}

Then by (\ref{con1}), we can apply Proposition  \ref{ehe} to $A_\m+F_\m(\cdot)$, and get  $Y_{\tilde m}, f_{\tilde m}^{(re)}, F_{\tilde m,1}\in C^\omega_{h_{\tilde m}, \delta_{\tilde m}}(\T^d\times\tilde\Lambda, gl(m,\C))$ such that 
\[e^{-Y_{\tilde m}(\cdot+\alpha)}(A_{\tilde m}+F_{\tilde m}(\cdot))e^{Y_{\tilde m}(\cdot)}=A_{\tilde m,1}+f_{\tilde m}^{(re)}(\cdot)+F_{\tilde m,1}(\cdot),\]
where $A_{\tilde m,1}=\textrm{diag}\{A_{\tilde m,1,11}, \cdots, A_{\tilde m,1, l_{\tilde m}l_{\tilde m}}\}$, $|Y_{\tilde m}|_{h_{\tilde m}, \delta_{\tilde m}}\leq \epsilon_{\tilde m}^{1/2}$, $|f_{\tilde m}^{(re)}|_{h_{\tilde m},\delta_{\tilde m}}\leq 3\tilde M_{\tilde m}\epsilon_{\tilde m}$, 
\begin{eqnarray*}
& |A_{\tilde m,1,ii}-A_{\tilde m,ii}|_{\delta_{\tilde m}}\leq 3\tilde M_{\tilde m}\epsilon_{\tilde m}, \\
& |F_{\tilde m,1}|_{h_{\tilde m+1},\delta_{\tilde m}}\leq 3\tilde M_{\tilde m}\epsilon_{\tilde m}e^{-2\pi h_1a_{\tilde m}N_{\tilde m,p}/2^{\tilde m+1}},
\end{eqnarray*}
 and $f_{\tilde m}^{(re)}$ satisfies the properties of $f^{(re)}$ in Proposition \ref{ehe} with $N=N_{\m,p}, a=a_{\m}, K=K_\m$.
Moreover, if we denote $\tilde\Sigma(\lambda):=\Sigma(A_{\tilde m,1}(\lambda))$,  and $\tilde \Sigma_i(\lambda):=\Sigma(A_{\tilde m,1,ii}(\lambda))$  for $1\leq i\leq l_{\tilde m}$,
then by (\ref{con1}), we can apply Lemma \ref{ndh}, and obtain that 
\begin{enumerate}[1)]
\item $|A_{\tilde m,1}|_{\delta_{\tilde m}}\leq \tilde M'$,
\item $d_H(\tilde\Sigma_i(\lambda),  \Sigma_i^{(\m)}(\lambda))<8m^2\tilde M_{\m}^3\epsilon_\m^{1/m}$, $\tilde\Sigma (\lambda)\subseteq D(R')$ for $\lambda\in W_{\delta_{\tilde m}}(\tilde\Lambda)$,
\item $\tilde\Sigma=\tilde\Sigma_{1} \cup\cdots\cup\tilde\Sigma_{l_{\tilde m}} $ is $(\nu',\zeta')$-separated in $W_{\delta_{\tilde m}}(\tilde{\Lambda})$,
\item $\tilde\Sigma =\tilde\Sigma_{1} \cup\cdots\cup\tilde\Sigma_{l_{\tilde m}} $ is $(M',\delta_{\tilde m},c',r_{\tilde m})$-transverse on $\tilde{\Lambda}$,
\end{enumerate}
where
$$
\begin{matrix}
\begin{array}{l}
R'=R_{\tilde m}+16m^{2}R_{\tilde m}^{2}\tilde{M}_{\tilde m}^3\epsilon_{\tilde m}^{\frac{1}{m}},\\  \nu'=\nu_{\tilde m}-16m^{2} \tilde{M}_{\tilde m}^3\epsilon_{\tilde m}^{\frac{1}{m}},\\
c'=c_{\tilde m}-2^{8m+2}m^{10m^2}R_{\tilde m}^{3m^2}\tilde{M}_{\tilde m}^{m+1} (\frac{2r_{\tilde m}}{\delta_{\tilde m}})^{ r_{\tilde m}}\epsilon_{\tilde m},
\end{array}
&
\begin{array}{l}
\tilde{M}'=\tilde{M}_{\tilde m}+3\tilde M_{\tilde m}\epsilon_{\tilde m}, \\ \zeta'=\zeta_{\tilde m}+16m^{2} \tilde{M}_{\tilde m}^3\epsilon_{\tilde m}^{\frac{1}{m}}, \\ M'=(2R')^{m^{2}}.
\end{array}
\end{matrix}
$$

\subsubsection{Step 2: Remove the resonances}

Now  applying Lemma \ref{rmr} to $A_{\tilde m,1}$ and $f_{\tilde m}^{(re)}$, there exists $H_{\tilde m}\in C^\omega_{h_{\tilde m}}(\T^d, \GL(m,\C))$, independent of $\lambda$, such that 
\[H_{\tilde m}^{-1}(\cdot+\alpha)(A_{\tilde m,1}+f_{\tilde m}^{(re)}(\cdot))H_{\tilde m}(\cdot)=A_{\tilde m,2},\]
with $A_{\tilde m,2,ii}=e^{2\pi \ii\la k_i',\alpha\ra}A_{\tilde m,1,ii}$ for some $k_i'\in\Z^d$, where $k_i'$ is independent of $\lambda$ for all $i=1,\cdots,l_{\tilde m}$, and
\begin{eqnarray*}
&&|H_{\tilde m}|_{h_{\tilde m+1}}\leq e^{2\pi mN_{\tilde m,p}h_{\tilde m+1}},\\
&&|A_{\tilde m,2}-\textrm{diag}\{A_{\tilde m,2,11},\cdots, A_{\tilde m,2,l_{\tilde m}l_{\tilde m}}\}|_{\delta_{\tilde m}}\leq 3(m-1)\tilde M_{\tilde m}\epsilon_{\tilde m}.
\end{eqnarray*}
Let $F_{\tilde m,2}(\cdot)=H^{-1}_{\tilde m}(\cdot+\alpha)F_{\tilde m,1}(\cdot)H_{\tilde m}(\cdot)$. Then we have 
\begin{eqnarray*}
 |F_{\tilde m,2}|_{h_{\tilde m+1},\delta_{\tilde m}} \leq  |H_{\tilde m}|_{h_{\tilde m+1}}^{2}|F_{\tilde m,1}|_{h_{\tilde m+1},\delta_{\tilde m}} 
 \leq  3\tilde M_{\tilde m}\epsilon_{\tilde m}^{1+( \frac{a_{\tilde m}}{2^{\tilde m+1}}-2m) 2^{\tilde m+1} a_{\tilde m}^p}\leq 3\tilde M_{\tilde m}\epsilon_{\tilde m}^{1+4ma_{\tilde m}^p\cdot 8^{\tilde m}}.
\end{eqnarray*}
Furthermore, if we denote $ A_{\tilde m,2}':=\textrm{diag}\{A_{\tilde m,2,11},\cdots, A_{\tilde m,2,l_{\tilde m}l_{\tilde m}}\}$,  $\hat{\Sigma}(\lambda):=\Sigma(A_{\tilde m,2}'(\lambda))$,   and $\hat{\Sigma}_{i}(\lambda):=\Sigma(A_{\tilde m,2,ii}(\lambda))$ for $1\leq i\leq l_{\tilde m}$, then the following hold: 
\begin{Lemma}\label{claim-est-A-m-2}
The decomposition $\hat \Sigma(\lambda)=\hat{\Sigma}_{1}(\lambda)\cup\cdots\cup\hat{\Sigma}_{l_{\tilde m}}(\lambda)$ is   $(M',\delta_{\tilde m},c',r_{\tilde m})$-transverse on $\tilde{\Lambda}$, for $\lambda\in W_{\delta_{\tilde m}}(\tilde\Lambda)$, $\hat \Sigma(\lambda)\subseteq D(R')$, $\mathrm{diam}(\hat\Sigma_i(\lambda))\leq \zeta'$, and 
\begin{equation}\label{est-A-m-2}
|A_{\tilde m,2}'|_{\delta_{\tilde m}}\leq \tilde M',\ \ \hat\Sigma_i(\lambda)=e^{2\pi \ii \la  k_i',\alpha\ra}\tilde\Sigma_i(\lambda).\end{equation}
\end{Lemma}
\begin{pf}
Since $A_{\tilde m,2,ii}=e^{2\pi \ii\la k_i',\alpha\ra}A_{\tilde m,1,ii}$, we can get (\ref{est-A-m-2}) and then by the properties of $\tilde\Sigma_i(\lambda)$, we have
for $\lambda\in W_{\delta_{\tilde m}}(\tilde\Lambda)$,  $\hat \Sigma(\lambda)\subseteq D(R')$, $\mathrm{diam}(\hat\Sigma_i(\lambda))\leq \zeta'$. Moreover, one can check that for any $u\in\T$ and $\lambda\in W_{\delta_{\tilde m}}(\tilde\Lambda)$,
\begin{eqnarray*}
&\hat g_i(\lambda,u)=e^{2\pi \ii m_i(m_i-1)\la k_i', \alpha\ra} \tilde g_i(\lambda,u),\\
&\mathrm{Res}(\chi_{\hat{\Sigma}_{i}(\lambda)},\chi_{\hat{\Sigma}_{j}(\lambda)};u)=e^{2\pi \ii m_im_j\la  k_i',\alpha\ra}\mathrm{Res}(\chi_{\tilde\Sigma_{i}(\lambda)},\chi_{\tilde\Sigma_{j}(\lambda)}; u+\la   k_j'-  k_i',\alpha\ra),
\end{eqnarray*}
where $m_i=\# \hat\Sigma_i$ is independent of $\lambda$ and $\theta$,  
$
\tilde g_i(\lambda,u)=\prod_{\sigma_{\ell_1},\sigma_{\ell_2}\in\tilde\Sigma_i(\lambda) \atop{\ell_1\neq \ell_2}}(\sigma_{\ell_1}-e^{2\pi\ii u}\sigma_{\ell_2})$, and $\hat g_i(\lambda,u)=\prod_{\sigma_{\ell_1},\sigma_{\ell_2}\in\hat\Sigma_i(\lambda) \atop{\ell_1\neq \ell_2}}(\sigma_{\ell_1}-e^{2\pi\ii u}\sigma_{\ell_2})$.
Then the decomposition $\hat \Sigma(\lambda)=\hat{\Sigma}_{1}(\lambda)\cup\cdots\cup\hat{\Sigma}_{l_{\tilde m}}(\lambda)$ is   $(M',\delta_{\tilde m},c',r_{\tilde m})$-transverse on $\tilde{\Lambda}$.
\end{pf}

\subsubsection{Step 3: Block diagonalization} We will divide this step into three cases according to the spectrum of $A_{\tilde m,2}'$, and  one of the following situations will happen:
 
 \smallskip
\textbf{Case $\mathbf{a: s_{k}\leq \tilde m<s_{k+1}-1}$ and $\mathbf{\hat{\Sigma}_{1}\cup\cdots\cup\hat{\Sigma}_{l_{\tilde m}}}$ is $\mathbf{(\nu',\zeta')}$-separated in $\mathbf{W_{\delta_{\tilde m}}(\tilde\Lambda)}$.}
In this case, since by (\ref{con1})
\[(6m^4\tilde M'\nu'^{-1})^{m^2+1}(3(m-1)\tilde M_{\tilde m}\epsilon_{\tilde m})^\frac{1}{2}<1,\ \
2^7m^2 \tilde  M'^3 R'(3(m-1)\tilde M_{\tilde m}\epsilon_{\tilde m})^{\frac{1}{m}}<1,\]
then Lemma \ref{norm2} applies and there exists $Z_{\tilde m,1}\in C_{\delta_{\tilde m}}^{\omega}(\tilde{\Lambda},\mathrm{GL}(m,\C))$ such that
\[
    Z_{\tilde m,1}^{-1}A_{\tilde m,2}Z_{\tilde m,1}=\mathrm{diag}\{ A_{\tilde m,3, 11},\cdots, A_{\tilde m,3, l_{\tilde m }l_{\tilde m} }\},
\]
 with  $|Z_{\tilde m,1}-\mathrm{Id}|_{\delta_{\tilde m}}\leq 4m\tilde M_{\tilde m}^{\frac{1}{2}}\epsilon_{\tilde m}^{\frac{1}{2}}$ and
 $|  A_{\tilde m,3,ii}-A_{\tilde m, 2,ii}|_{\delta_{\tilde m}}< 18m\tilde M_{\tilde m}^2\epsilon_{\tilde m}$. Let 
 $l_{\tilde m+1}=l_{\tilde m}$, 
 \begin{eqnarray*}
 &&\#\Pi_{\tilde m+1}(\tilde\Lambda)=1,\  \ A_{\tilde m+1}=\mathrm{diag}\{A_{\tilde m,3, 11},\cdots, A_{\tilde m,3, l_{\tilde m }l_{\tilde m} }\}, \\
&&F_{\tilde m+1}=Z_{\tilde m,1}^{-1}F_{\tilde m,2}Z_{\tilde m,1},\ \ Z_{\tilde m}=e^{Y_{\tilde m}}H_{\tilde m}Z_{\tilde m,1}.
\end{eqnarray*}
 Then
\begin{eqnarray*}
&&|A_{\tilde m+1}|_{\delta_{\tilde m}}\leq\tilde M_{\tilde m}+20m\tilde M_{\tilde m}^2\epsilon_{\tilde m}=\tilde M_{\tilde m+1},\\
&&|F_{\tilde m+1}|_{h_{\tilde m+1},\delta_{\tilde m}}\leq 6\tilde M_{\tilde m}\epsilon_{\tilde m}^{1+4m8^{\tilde m}a_{\tilde m}^p}<\epsilon_{\tilde m}^{4^{\tilde m}}=\epsilon_{\tilde m+1},\\
&&|Z_{\tilde m}|_{h_{\tilde m+1},\delta_{\tilde m}}\leq 2|H_{\tilde m}|_{h_{\tilde m+1}}\leq 2\epsilon_{\tilde m}^{-m2^{\tilde m+2}a_{\tilde m}^p}<\epsilon_{\tilde m+2m+2}^{-1},\\
&&|Z_{\tilde m} |^{\tilde m^2/2}_{h_{\tilde m+1},\delta_{\tilde m}}|F_{\tilde m+1}|_{h_{\tilde m+1},\delta_{\tilde m}}\leq 3\tilde M_{\tilde m}  \epsilon_\m^{1+4m a_{\m}^p(8^\m- \m^22^\m)}<\epsilon_\m^{4ma_\m^p2^\m} <1,
\end{eqnarray*}
and the same estimates hold for $Z_\m^{-1}$. 

Moreover, $\Sigma^{(\m+1)}(\lambda)\subseteq D(R_{\m+1})$ for $\lambda\in W_{\delta_{\m+1 }}(\tilde{\Lambda})$, $\Sigma^{(\m+1)}=\Sigma_1^{(\m+1)}\cup\cdots\cup\Sigma_{l_{\m+1}}^{(\m+1)}$
 is $((2R_{\m+1})^{m^{2}},\delta_{\m+1},c_{\m+1},r_{\m+1})$-transverse on $\tilde{\Lambda}$ and $(\nu_{\m+1},\zeta_{\m+1})$-separated in $W_{\delta_{\m+1 }}(\tilde{\Lambda})$, where $R_{\m+1}, r_{\m+1}$, $\nu_{\m+1}$, $\zeta_{\m+1}$, $\delta_{\m+1}$ and $c_{\m+1}$ are defined as in $(\mathcal P4.2)_{\tilde m+1}$. 
 
 Furthermore, 
 for $\lambda\in W_{\delta_{\m+1}}(\tilde\Lambda)$, $1\leq i\leq l_{\m+1}$,
 \begin{equation*}
d_H( \hat \Sigma_i(\lambda), \Sigma_i^{(\m+1)}(\lambda))<40m^2\tilde M_{\m}^4\epsilon_{\m}^{\frac{1}{m}}<\frac{\nu_\m-\nu_{\m+1}}{4}.
   \end{equation*}
 By the fact that $\hat\Sigma_i=e^{2\pi \ii \la  k_i',\alpha\ra}\tilde\Sigma_i$ and $d_H(\tilde\Sigma_i(\lambda),  \Sigma_i^{(\m)}(\lambda))<8m^2\tilde M_{\m}^3\epsilon_\m^{1/m}$ in $W_{\tilde\delta_\m}(\tilde\Lambda)$, we get that for $\lambda\in W_{\delta_{\m+1}}(\tilde\Lambda)$
  \begin{eqnarray*}
  d_H(\Sigma_i^{(\m+1)}(\lambda), e^{2\pi \ii\la  k_i',\alpha\ra}\Sigma_i^{(\m)}(\lambda))<\frac{\nu_\m-\nu_{\m+1}}{2}.
   \end{eqnarray*}

 \smallskip
 
\textbf{Case $\mathbf{b: s_{k}\leq \m<s_{k+1}-1}$ and $\mathbf{\hat{\Sigma}_{1}\cup\cdots\cup\hat{\Sigma}_{l_{\m}}}$ is not $\mathbf{(\nu',\zeta')}$-separated in $\mathbf{W_{\delta_{\tilde m}}(\tilde\Lambda)}$}. 
By (\ref{con1}), (\ref{equ-pro-3-5}), (\ref{equ-est-mu-zeta}), and the fact that $d_H(\tilde\Sigma_i(\lambda),  \Sigma_i^{(\m)}(\lambda))<8m^2\tilde M_{\m}^3\epsilon_\m^{1/m}$ and $\hat\Sigma_i=e^{2\pi i \la  k_i',\alpha\ra}\tilde\Sigma_i$ with $i=1,\cdots,l_\m$, we have for any $\lambda\in W_{\delta_\m}(\tilde\Lambda)$,
\begin{equation}\label{equ-pro-3-6}
d_H(\hat\Sigma_i(\lambda),\hat\Sigma_i(\lambda_0))=d_H(\tilde\Sigma_i(\lambda),\tilde\Sigma_i(\lambda_0))<2\nu_{m_1},
\end{equation}
 for some $\lambda_{0}\in\tilde{\Lambda}$. Then we can regathering the subsets $\hat \Sigma_1(\lambda), \cdots, \hat \Sigma_{l_\m}(\lambda)$ so that the following holds:
 \begin{Lemma}
 The set $\{1, 2, \cdots, l_\m\}$ can be divided into disjoint union
 \[
  \{1, 2, \cdots, l_\m\}=\tilde S_{1}\cup\cdots\cup\tilde S_{l_{\m+1}}, \]
  such that the decomposition $\hat\Sigma(\lambda)=\acute\Sigma_1(\lambda)\cup\cdots\cup\acute \Sigma_{l_{\m+1}}(\lambda)$ is $(\nu', \bar\zeta)$-separated in $W_{\delta_\m}(\tilde\Lambda)$,  and $(M', \delta_\m, c'', m^2r_\m)$-transverse on $\tilde\Lambda$, where $l_{\m+1}<l_{\m}$,  $\   \bar\zeta=m(\zeta'+\nu'+4\nu_{m_1})$, and
  \[
  \acute{\Sigma}_{j}=\bigcup_{i\in\tilde S_{j}}\hat{\Sigma}_{i},\ \ c''=((\frac{\delta_\m}{4m^4r_\m^2(2R')^{m^2}})^{2m^4r_\m}c')^{m^{2m^2r_\m+2}}
  \]
 \end{Lemma}
 \begin{pf}

 We divide the set $\{1,\cdots,l_{\m}\}$ into subsets $\tilde S_{1},\cdots,\tilde S_{k}$ such that
    \begin{enumerate}[(1)]
    \item If $\ell\neq \ell'$ belong to the same $\tilde S_{i}$, then there exist $\ell_{1}, \ell_2,\cdots,\ell_{t}$ with $\ell_1=\ell, \ell_t=\ell'$, such that $\mathrm{dist}(\hat{\Sigma}_{\ell_{j}}(\lambda_{j}),\hat{\Sigma}_{\ell_{j+1}}(\lambda_{j}))\leq\nu'$ for some $\lambda_{j}\in W_{\delta_{\m}}(\tilde{\Lambda})$.\footnote{Here, for any  $\Omega_1,\Omega_2\subseteq \C$, we denote $\mathrm{dist}(\Omega_1,\Omega_2)=\textrm{inf}_{x\in\Omega_1, y\in\Omega_2}|x-y|$.}
    \item If $\ell$ and $\ell'$ belong to different subsets $\tilde S_{i_1}, \tilde S_{i_2}$, then there exists no such path, implying that  $\mathrm{dist}(\hat{\Sigma}_{\ell}(\lambda),\hat{\Sigma}_{\ell'}(\lambda))>\nu'$ for all $\lambda\in W_{\delta_{\m}}(\tilde{\Lambda})$.
    \end{enumerate}
Since $\hat{\Sigma}_{1}\cup\cdots\cup\hat{\Sigma}_{l_{\m}}$ is not $(\nu',\zeta')$-separated in $W_{\delta_{\tilde m}}(\tilde\Lambda)$ and   $\forall \lambda\in W_{\delta_\m}(\tilde\Lambda)$,  $\textrm{diam}(\hat\Sigma_i(\lambda))\leq \zeta'$,  then $k<l_\m$, and we denote $l_{\m+1}=k$. 

 Let $\acute{\Sigma}_{j}=\bigcup_{i\in\tilde S_{j}}\hat{\Sigma}_{i}$.
Then  $\hat\Sigma=\acute\Sigma_1\cup\cdots\cup\acute \Sigma_{l_{\m+1}}$. By (\ref{equ-pro-3-6}) and $\textrm{diam}(\hat \Sigma_i(\lambda_0))\leq \zeta'$ for $\lambda\in W_{\delta_\m}(\tilde\Lambda)$ (Lemma \ref{claim-est-A-m-2}), we can get that
$\hat\Sigma=\bigcup_{j=1}^{l_{\m+1}}\acute{\Sigma}_{j}$ is a $(\nu',\bar{\zeta})$-separated decomposition in $W_{\delta_{\m}}(\tilde{\Lambda})$ where
$\bar{\zeta}=m(\zeta'+\nu'+4\nu_{m_{1}})$. 

Moreover, 
by Lemma \ref{claim-est-A-m-2} and Lemma \ref{hebin}, we can get that  $\hat \Sigma=\acute \Sigma_1\cup\cdots\cup\acute\Sigma_{l_{\m+1}}$ is $(M', \delta_\m, c'', m^2r_\m)$-transverse on $\tilde\Lambda$ with $c''=((\frac{\delta_\m}{4m^4r_\m^2(2R')^{m^2}})^{2m^4r_\m}c')^{m^{2m^2r_\m+2}}$. 
\end{pf}

Besides, there exists $L_\m\in \GL(m,\R)$ which is the product of elementary matrices that exchange the rows, with $\|L_\m\|=\|L^{-1 }_\m\|=1$, such that the conjugated block diagonal matrix 
\[
L_\m^{-1} A_{\m,2}'L_\m=:\acute A_{\m,4}=\textrm{diag}\{\acute A_{\m, 4,11},\cdots, \acute A_{\m,4,l_{\m+1}l_{\m+1}}\},\]
 has the property that $\acute A_{\m,4,jj}=\textrm{diag}\{A_{\m,2,ii}: i\in\tilde S_j\}$. Let 
 \[
  A_{\m, 4}=L_\m^{-1} A_{\m,2}L_\m, \ \  F_{\m,4}=L_\m^{-1} F_{\m,2}L_\m.
  \]
   Then we have 
$|\acute A_{\m,4}-A_{\m,4}|_{\delta_\m}\leq 3(m-1)\tilde M_\m\epsilon_\m$, and 
\[ \Sigma(\acute A_{\m, 4}(\lambda))=\hat\Sigma(\lambda),  \ \Sigma(\acute A_{\m, 4,ii}(\lambda))=\acute \Sigma_i(\lambda).
\]

Since 
\[|\acute A_{\m, 4}|_{\delta_\m}\leq | A_{\m, 2}'|_{\delta_\m}\leq \tilde M', \ \ |\acute A_{\m, 4}-A_{\m, 4}|_{\delta_\m}\leq 3(m-1)\tilde M_\m\epsilon_\m,\]
 and 
\[ (6m^4\tilde M' \nu'^{-1})^{m^2+1}(3(m-1)\tilde M_\m\epsilon_\m)^{\frac{1}{2}}<1,\ 2^7m^2 \tilde M'^3R'(3(m-1)\tilde M_\m\epsilon_\m)^{\frac{1}{m}}<1,\]
then by Lemma \ref{norm2}, there exists $Z_{\m, 2}\in C^\omega_{\delta_\m}(\tilde \Lambda, \GL(m,\C))$ such that
\[Z_{\m,2}^{-1}A_{\m,4}Z_{\m,2}=A_{\m+1}=\textrm{diag}\{  A_{\m+1,11}, \cdots,   A_{\m+1, l_{\m+1}l_{\m+1} }\},\]
with $| A_{\m+1,ii}- \acute A_{\m, 4,ii}|_{\delta_\m}< 9m\tilde M'^2\epsilon_\m$, $|Z_{\m,2}-\textrm{Id}|_{\delta_\m}<4m\tilde M_\m^{\frac{1}{2}}\epsilon_\m^{\frac{1}{2}}$. 
Thus, we have 
\[|A_{\m+1}|_{\delta_\m}\leq \tilde M_\m+20m\tilde M_\m^2\epsilon_\m=\tilde M_{\m+1}.\]
Moreover, the decomposition $\Sigma^{(\m+1)}(\lambda)= \Sigma_1^{(\m+1)}(\lambda)\cup\cdots \Sigma_{l_{\m+1}}^{(\m+1)}(\lambda)$ is $(\nu_{\m+1},  \zeta_{\m+1})$-separated in $W_{\delta_\m}(\tilde\Lambda)$, $(2R_{\m+1})^{m^2}, \delta_\m, c_{\m+1},  r_{\m+1})$-transverse on $\tilde\Lambda$ and  for $\lambda\in W_{\delta_\m}(\tilde\Lambda)$, $ \Sigma^{(m+1)}(\lambda)\subseteq D(R_{\m+1})$, where $\nu_{\m+1},$ $\zeta_{\m+1},$ $R_{\m+1},$ $c_{\m+1},$ $r_{\m+1}$ are defined as in $(\mathcal P4.1)_{\m+1}$. 

\begin{Claim}
 There exists a partition $\Pi_{\m+1}(\tilde\Lambda)$ of $\tilde\Lambda$ such that for any $\bar\Lambda\in\Pi_{\m+1}(\tilde\Lambda)$, there exists $\lambda_0\in\bar\Lambda$ that 
 \[
 d_H(\Sigma^{(\m+1)}_i(\lambda), \Sigma_i^{(\m+1)}(\lambda_0))<\nu_{\m+1}, \ \forall \lambda\in W_{\delta_{\m+1}}(\bar\Lambda),
 \]
 where $\delta_{\m+1}={b}^{-1}R_{\m}^{-(2m-1)} \nu_{\m+1}^{m}\delta_{\m}$.
\end{Claim}
 \begin{pf}
Divide $\tilde\Lambda$ into $1+  [\frac{\tilde\Lambda}{\delta_{\m+1}} ]$ intervals with equal length less than $\delta_{\m+1}$, which is the partition of $\tilde\Lambda$, denoted by $\Pi_{\m+1}(\tilde\Lambda)$. For any $\bar\Lambda\in \Pi_{\m+1}(\tilde\Lambda)$, fixing some $\lambda_0\in\bar\Lambda$, we have for any $\lambda\in W_{\delta_{\m+1}}(\bar\Lambda)$, $|\lambda-\lambda_0|<2\delta_{\m+1}$. Since $\delta_{\m+1}< \frac{\nu_{\m+1}^m}{4m^{m+1}R_{\m+1}^{2m-1}}\delta_\m$, then similar as  (\ref{equ-est-different-parameter}), we can obtain the result.
\end{pf}
Now we let  $F_{\m+1}=Z_{\m,2}^{-1}F_{\m,4}Z_{\m,2}$, and $Z_\m=e^{Y_\m}H_\m L_\m Z_{\m,2}$.  Then  the result follows.

 \medskip

\textbf{Case $\mathbf{c: \m=s_{k+1}-1}$.} By Lemma \ref{claim-est-A-m-2} and Lemma \ref{hebin}, $\hat \Sigma(\lambda)$ is $(M',\delta_\m,\tilde c', m^2r_\m)$-transverse on $\tilde \Lambda$ with $$\tilde c'=((\frac{\delta_\m}{4m^4r_\m^2(2R')^{m^2}})^{m^4r_\m}c')^{m^{2m^2r_\m+2}}.$$ 
Denote $\Sigma'(\lambda):=\Sigma(A_{\m,2}(\lambda))$, $R_{\m+1}=R_{\m}+  b R_{\m}^2\tilde M_{\m}^4\epsilon_{\m}^{\frac{1}{m}}$, and 
$$\tilde c''=\tilde c'-2^{8m+2}m^{10m^2+1}R_\m^{3m^2}\tilde M_\m^{m+1}(\frac{2m^2r_\m}{\delta_\m})^{m^2r_\m}\epsilon_\m.$$
 Then due to the fact that 
\[ 64m^2\tilde M'^2R'(3(m-1)\tilde M_\m\epsilon_\m)^{\frac{1}{m}}<1, \ \ |A_{\m,2}- A_{\m,2}'|_{\delta_\m}\leq 3(m-1)\tilde M_\m\epsilon_\m,
\]
 by   Lemma \ref{ndh} with $l=1$, we can obtain that $\Sigma'( \lambda )\subseteq D(R_{\m+1})$ for all $\lambda\in W_{\delta_\m}(\tilde\Lambda)$ and the multiset $\Sigma'(\lambda)$  is $((2R_{\m+1})^2, \delta_\m, \tilde c'', m^2r_\m)$-transverse on $\tilde\Lambda$.
 
  Let  $\delta_{\m+1}, r_{\m+1}$, $\nu_{\m+1}$, $\zeta_{\m+1}$, $\tilde{M}_{\m+1}$ and $c_{\m+1}$ be defined as in $(\mathcal P3)_{\m+1}$.
Then by Proposition \ref{Nth} with $\nu_{\m+1}$, there exists  a partition $\Pi_{\m+1}(\tilde\Lambda)$ for $\tilde{\Lambda}$ with
\[
\#\Pi_{\m+1}(\tilde\Lambda)\leq 1+\frac{|\tilde\Lambda|}{\delta_{\m+1}},
\] 
such that for each $\bar{\Lambda}\in\Pi_{\m+1}(\tilde\Lambda)$ there is a similarity transformation $Z_{\m,3}\in C_{\delta_{\m+1}}^{\omega}(\bar{\Lambda},\mathrm{GL}(m,\C))$ such that
$$
Z_{\m,3}^{-1}A_{\m,2}Z_{\m,3}=A_{\m+1}= \mathrm{diag}\{A_{\m+1,11},\cdots,A_{\m+1,l_{\m+1}l_{\m+1}}\},
$$
with $|Z_{\m,3}|_{\delta_{\m+1}},|Z_{\m,3}^{-1}|_{\delta_{\m+1}},|A_{\m+1,ii}|_{\delta_{\m+1}}\leq \tilde{M}_{\m+1}$.
Then one can check that $A_{\m+1}\in C_{\delta_{\m+1}}^\omega(\bar\Lambda, \GL(m,\C))$ satisfies $(\mathcal P1)_{\m+1}$ and $(\mathcal P2)_{\m+1}$.

Moreover, there exists $\lambda_0\in\bar\Lambda$ such that the decomposition $\Sigma^{(\m+1)}(\lambda_0)=\Sigma^{(\m+1)}_1(\lambda_0)\cup\cdots\cup\Sigma_{l_{\m+1}}^{(\m+1)}(\lambda_0)$ is $(8\nu_{\m+1},\zeta_{\m+1})$-separated and $\forall \lambda\in\ W_{\delta_{\m+1}}(\bar\Lambda)$, we have 
\[d_H(\Sigma^{(\m+1)}_i(\lambda), \Sigma_i^{(\m+1)}(\lambda_0))<\nu_{\m+1}. \]

Let $\Pi_{\m+1}=\cup_{\tilde\Lambda\in\Pi_\m}\Pi_{\m+1}(\tilde\Lambda)$. Then
\begin{eqnarray*}
\#\Pi_{\m+1}&=&\sum_{\tilde\Lambda\in\Pi_\m}\#\Pi_{\m+1}(\tilde\Lambda)\leq \sum_{\tilde\Lambda\in\Pi_\m}(1+\frac{|\tilde\Lambda|}{\delta_{\m+1}})\\
&\leq&  \#\Pi_\m+\frac{bR_\m^{6m}\tilde M_\m^{3m}|\Lambda|}{\nu_{\m+1}^{3m} }\max_{\tilde\Lambda\in\Pi_\m}\frac{1}{\delta_\m(\tilde\Lambda)},
\end{eqnarray*}
since $\delta_\m$ may vary with $\tilde\Lambda\in\Pi_\m$, .

Let 
\[ F_{\m+1}=Z_{\m,3}^{-1}F_{\m,2}Z_{\m,3},\ \  Z_\m=e^{Y_\m}H_\m Z_{\m,3}.
\]
By (\ref{equ-epsilon-add-1}), (\ref{equ-epsilon-add-2}), (\ref{equ-pro-3-con-2}), (\ref{con1}), and the definition of $u_{s_{k+1}}$ in (\ref{equ-def-u-n}), we can check that $2\tilde M_{\m+1}<\epsilon_\m^{-m}$, and the following estimates hold:
\[|F_{\m+1}|_{h_{\m+1},\delta_{\m+1}}\leq \tilde M_{\m+1}^2|F_{\m,2}|_{h_{\m+1},\delta_\m}<\epsilon_\m^{4m8^\m-2m}<\epsilon_{\m+1},\]
\[|Z_\m|_{h_{\m+1},\delta_{\m+1}}^{\frac{\m^2}{2}}|F_{\m+1}|_{h_{\m+1},\delta_{\m+1}}<\epsilon_\m^{1+m(4\cdot 8^\m a_m^p-2^{\m+1}\m^2a_\m^p-\frac{\m^2}{2}-2)}<\epsilon_\m,\]
\[|Z_{\m}|_{h_{\m+1},\delta_{\m+1}}<2\tilde M_{\m+1}\epsilon_\m^{-m2^{\m+2}a_\m^p}<\epsilon_{\m}^{-m(2^{\m+2}a_\m^m+1)}<\epsilon_{\m+2m+2}^{-1},\]
and the same estimates hold for $Z_\m^{-1}$.

\bigskip

To that end,  for $s_k\leq \m<s_{k+1}-1$,  in both cases (Case a and b), we can obtain that
\begin{eqnarray*}
\#\Pi_{\m+1} = \sum_{\tilde\Lambda\in\Pi_\m}\#\Pi_{\m+1}(\tilde\Lambda)\leq \sum_{\tilde\Lambda\in\Pi_\m}(1+\frac{bR_\m^{2m-1}|\tilde\Lambda|}{\nu_{\m+1}^m\delta_{\m }})
 \leq  \#\Pi_{\m}+\frac{bR_\m^{2m-1} |\Lambda|}{\nu_{\m+1 }^{ m}}\max_{\tilde\Lambda\in\Pi_\m}\frac{1}{\delta_\m}.
\end{eqnarray*}
The proof is finished. \qed

\subsection{Proof of Proposition \ref{pro-iterative-lemma-full}:}
 By (\ref{equ-est-epsilon-1}) and (\ref{equ-est-epsilon-2}), inductively we can check that for any $j\geq 1$,
 \begin{equation}\label{equ-thm-2-4}
 \tilde M_j'<2\tilde M_1', \ \ \tilde R_j<2\tilde R_1, \ \ \tilde\delta_j>\tilde\epsilon_j^{\frac{1}{4\tilde r_1'}}, \ \ (2^j |\ln\tilde\epsilon_j|)^{\tau+d}<\tilde\epsilon_j^{-\frac{\varrho}{100\tilde r_1'^2}}.
 \end{equation}
 Then by (\ref{equ-est-epsilon-2}), we can obtain that
 \[
 \frac{(\tilde r_1'+1)!}{(\tilde\delta_{j+1}/2)^{\tilde r_1'+1}}(8m^{4}\tilde R_j \tilde M_j')^{ 3m^2 }\tilde\epsilon_j<\tilde\epsilon_{j}^{-\frac{1}{32}}\tilde\epsilon_{j+1}^{-\frac{\tilde r_1'+1}{4\tilde r_1'}}\tilde\epsilon_j\leq \tilde\epsilon_j^{\frac{1}{5}},
 \]
 which implies that $\tilde C_j<\tilde C_1+2\tilde\epsilon_1^{\frac{1}{5}}<2\tilde C_1$, and $\tilde c_j'>\tilde c_1'-2\tilde\epsilon_1^{\frac{1}{5}}>\frac{\tilde c_1'}{2}$ by (\ref{equ-est-epsilon-1}). Moreover, also by (\ref{equ-est-epsilon-1}), we can check that
 \[
 \tilde c_1'K_1>\tilde \epsilon_1^{\frac{1}{100\tilde r_1'^2}}\tilde\epsilon_1^{-\frac{1}{10\tilde r_1'}}>4,
 \]
 and thus for any $j\geq 1$, $\tilde c_j'>2K_j^{-1}$.
 
 \subsubsection{Selection of the subset $\Lambda^{(j+1)}\subseteq \Lambda^{(j)}$}

 \begin{Lemma}\label{lem-est-resonant-set}
 There exists $  \Lambda^{(j+1)}\subseteq \Lambda^{(j)}$ such that for any $\bar\Lambda^{(j+1)}\in\mathcal C(\Lambda^{(j+1)})$, $\lambda\in W_{\tilde\delta_{j+1}}(\bar \Lambda^{(j+1)})$, $0\leq |k|\leq \tilde N_j$, we have 
 \[
| g_j(\lambda, \la k,\alpha\ra) |>\frac{1}{2K_j},
 \]
 where 
 \[
\Lambda^{(j+1)}=\Lambda^{(j)}\backslash\mathcal R_j(\tilde\Lambda),\ \ 
\mathcal R_j(\tilde\Lambda)=\cup_{i\in J_j(\tilde\Lambda)} I^{(j)}_i,
\]
with $\{I_{i}^{(j)}\}_i\subseteq \Lambda^{(j)}$ being disjoint intervals, $|I_{i}^{(j)}|<\tilde\epsilon_j^{\frac{2}{25\tilde r_1'^2}}$, and   
\begin{eqnarray*}
&&\#J_j(\tilde\Lambda)\leq 2^{\tilde r_1' +d}\tilde N_j^d(\frac{8\tilde C_1}{\tilde c_1'}| \Lambda^{(j)}|+\#\mathcal C (\Lambda^{(j)})),\\
 && \#\mathcal C(\Lambda^{(j+1)})\leq 2\#\mathcal C(\mathcal R^{(j)}(\tilde\Lambda)).
\end{eqnarray*}
 \end{Lemma}
 \begin{pf}
By the assumptions in Proposition \ref{pro-iterative-lemma-full}, for any $\bar\Lambda^{(j)}\in\mathcal C(\Lambda^{(j)})$ and $u\in\T$,  the function $g_j(\lambda,u)$ is $(\tilde C_j, \tilde c_j', \tilde r_1')$-pyartli on  $\bar \Lambda^{(j)}$. Since  $\tilde c_j'>2K_j^{-1}$, then by Lemma \ref{Lem1}, for any $0\leq |k|\leq \tilde N_j$ there is a disjoint union of intervals $\cup_{i\in  J_{j,k}(\bar\Lambda^{(j)})}I_{k,i}( \bar \Lambda^{(j)})\subseteq \bar\Lambda^{(j)}$ such that
\begin{eqnarray}
\label{est-I-length} &&\max_{i\in J_{j,k}(\bar\Lambda^{(j)})}|I_{k,i}(\bar\Lambda^{(j)})|\leq 2(\frac{2K_j^{-1}}{\tilde c_j'})^{\frac{1}{\tilde r_1'}},\\
\label{est-J-number} &&\#J_{j,k}(\bar\Lambda^{(j)})\leq 2^{\tilde r_1' }(\frac{2\tilde C_j|\bar\Lambda^{(j)}|}{\tilde c_j'}+1)<2^{\tilde r_1' }(\frac{8\tilde C_1|\bar\Lambda^{(j)}|}{\tilde c_1'}+1),\\
\nonumber &&|g_j(\lambda,\la k,\alpha\ra)|\geq K_j^{-1}, \ \ \forall \lambda\in \bar \Lambda^{(j)}\backslash\cup_{i\in J_{j,k}(\bar\Lambda^{(j)})}I_{k,i}(\bar\Lambda^{(j)}).
\end{eqnarray}

Let 
\[
\mathcal R_j(\tilde\Lambda )=\cup_{|k|=0}^{\tilde N_j}\cup_{\bar\Lambda^{(j)}\in  \mathcal C(\Lambda^{(j)})}\cup_{i\in J_{j,k}(\bar\Lambda^{(j)})}I_{k,i}( \bar \Lambda^{(j)})=:\cup_{i\in J_j(\tilde\Lambda)}I_i^{(j)},\ \ \Lambda^{(j+1)}=\Lambda^{(j)}\backslash\mathcal R_j(\tilde\Lambda).
\]
Without loss of generality, we assume that $\#J_{j,k}(\bar\Lambda^{(j)})\geq 1$; otherwise, we can take an arbitrary interval $I_{k,1}(\bar\Lambda^{(j)})\subseteq \bar\Lambda^{(j)}$ with $|I_{k,1}(\bar\Lambda^{(j)})|\leq 2(\frac{2K_j^{-1}}{\tilde c_j'})^{\frac{1}{\tilde r_1'}}$. 

Then for any $\bar\Lambda^{(j+1)}\in \mathcal C(\Lambda^{(j+1)})$, $\lambda\in W_{\tilde \delta_{j+1}}(\bar\Lambda^{(j+1)})$, there exists $\lambda'\in \bar\Lambda^{(j+1)}$ such that $|\lambda-\lambda'|<\tilde\delta_{j+1}$.
Since $|g_j(\lambda,u)|_{\tilde\delta_j}\leq (2\tilde R_j)^{m^2}$, then by Cauchy estimate, we obtain that
\[
|g_j(\lambda,u)-g_j(\lambda',u)|\leq |\frac{\partial}{\partial\lambda}g_j(\lambda,u)|_{\tilde\delta_{j+1}}|\lambda-\lambda'|<\frac{(2\tilde R_j)^{m^2}}{\tilde\delta_j/2}\tilde\delta_{j+1}=\frac{1}{2K_j},
\]
which implies that for any $0\leq |k|\leq \tilde N_j$,
\[|g_j(\lambda,\la k,\alpha\ra)|\geq |g_j(\lambda',\la k,\alpha\ra)|-|g_j(\lambda,\la k,\alpha\ra)-g_j(\lambda', \la k,\alpha\ra)|>\frac{1}{2K_j}.\]

Moreover, by (\ref{est-J-number}), 
\begin{eqnarray*}
\#J_i(\tilde\Lambda)&\leq& (2\tilde N_j)^d \sum_{\bar\Lambda^{(j)}\in\mathcal C(\Lambda^{(j)})}2^{\tilde r_1' }(\frac{8\tilde  C_1|\bar \Lambda^{(j)}|}{\tilde c_1'}+1)\\
&\leq &2^{\tilde r_1'+d }\tilde N_j^d\left(\frac{8\tilde C_1}{\tilde c_1'}|\Lambda^{(j)}|+\#\mathcal C(\Lambda^{(j)})\right).
\end{eqnarray*}
Furthermore, by the definition of $\Lambda^{(j+1)}$, and $\#J_{j,k}(\bar\Lambda^{(j)})\geq 1$, we can obtain that 
\[
\#\mathcal C(\Lambda^{(j+1)})\leq \sum_{\bar\Lambda^{(j)}\in \mathcal C(\Lambda^{(j)})}(1+\# \mathcal R_j(\bar\Lambda^{(j)}))
\leq 
2\#\mathcal C(\mathcal R^{(j)}(\tilde\Lambda)),
\]
where $\# \mathcal R_j(\bar\Lambda^{(j)})=\cup_{|k|=0}^{\tilde N_j}\ \cup_{i\in J_{j,k}(\bar\Lambda^{(j)})}I_{k,i}( \bar \Lambda^{(j)})$.

 To that end, by (\ref{est-I-length}), for any interval $I_i^{(j)}\in\mathcal C(\mathcal R_j(\tilde\Lambda))$, we have 
 \[
 |I_i^{(j)}|\leq 2(\frac{2K_j^{-1}}{\tilde c_j'})^{\frac{1}{\tilde r_1'}}\leq 8\tilde\epsilon_j^{\frac{1}{10\tilde r_1'^2}}  \tilde c_1'^{-1}<\tilde\epsilon_j^{\frac{2}{25\tilde r_1'^2}}. \]
\end{pf}

\subsubsection{Eliminate the non-resonant terms}

\begin{Claim}
For any $\bar\Lambda^{(j+1)}\in\mathcal C(\Lambda^{(j+1)})$, $\lambda\in \bar\Lambda^{(j+1)}$, the 
multiset $\tilde\Sigma_j(\lambda)$ is $$(\tilde N_j, (2\gamma^{-1}K_j(2\tilde R_j)^{m^2}\tilde N_j^\tau)^{-1})$$-nonresonant with itself.
\end{Claim}
\begin{pf}
First, for any $0< |k|\leq \tilde N_j$, $\lambda\in W_{\tilde\delta_{j+1}} (\bar\Lambda^{(j+1)})$, $\sigma_{\ell_1}, \sigma_{\ell_2}\in\tilde\Sigma_j(\lambda)$, if $\ell_1\neq\ell_2$, then 
\begin{eqnarray*}
|g_j(\lambda,\la k,\alpha\ra)|&=&|\sigma_{\ell_1}-e^{2\pi \ii \la k,\alpha\ra}\sigma_{\ell_2}|\prod_{\sigma_{i_1},\sigma_{i_2}\in \tilde\Sigma_j(\lambda), \atop{i_1\neq i_2, (i_1, i_2)\neq (\ell_1,\ell_2) }}|\sigma_{i_1}-e^{2\pi \ii \la k,\alpha\ra}\sigma_{i_2}|\\
&\leq& (2\tilde R_j)^{m^2-m-1}|\sigma_{\ell_1}-e^{2\pi\ii\la k,\alpha\ra}\sigma_{\ell_2}|,
\end{eqnarray*}
which, by Lemma \ref{lem-est-resonant-set}, implies that 
\[
|\sigma_{\ell_1}-e^{2\pi\ii\la k,\alpha\ra}\sigma_{\ell_2}|\geq \frac{1}{2K_j(2\tilde R_j)^{m^2-m-1}}> \frac{\gamma}{2K_j(2\tilde R_j)^{m^2}\tilde N_j^\tau}.
\]
Moreover, for any $0< |k|\leq \tilde N_j$, $\lambda\in W_{\tilde\delta_{j+1}} (\bar\Lambda^{(j+1)})$, $\sigma\in\tilde\Sigma_j(\lambda)$,  
\[
|\sigma-e^{2\pi \ii\la k,\alpha\ra}\sigma|
     \geq \frac{1}{\tilde R_j}\cdot \frac{\gamma}{2|k|^\tau}\geq \frac{\gamma}{2\tilde R_j\tilde N_j^\tau}>   \frac{\gamma}{2K_j(2\tilde R_j)^{m^2}\tilde N_j^\tau}.
     \]
The result follows.
\end{pf}

 Then by Definition \ref{def-resonant-structrue} with $l=1$, for any $\bar\Lambda^{(j+1)}\in\mathcal C(\Lambda^{(j+1)})$, the multiset $\tilde\Sigma_j(\lambda)$ satisfies $\mathcal H(\tilde N_j, \tilde N_j, 2\gamma^{-1}K_j(2\tilde R_j)^{m^2}\tilde N_j^\tau)$ on $W_{\tilde\delta_{j+1}} (\bar\Lambda^{(j+1)})$. Moreover, by (\ref{equ-est-epsilon-2}) and (\ref{equ-thm-2-4}), one can check that 
\[
\tilde \epsilon_j<(12m^5\tilde M_j'^2 \gamma^{-1}K_j(2\tilde R_j)^{m^2}\tilde N_j^\tau)^{-2m^2}.
\]
Applying Proposition \ref{ehe} to $\tilde A_j+\tilde F_j$, we get $Y_j, \tilde F_j^{(re)}, \tilde F_{j+1}\in C^\omega_{\tilde h_j,\tilde \delta_{j+1}}(\T^d\times\bar\Lambda^{(j+1)}, gl(m,\C))$, such that 
\[e^{-Y_j^+(\cdot)}(\tilde A_j+\tilde F_j(\cdot))e^{Y_j(\cdot)}=\tilde A_{j+1}+\tilde F_j^{(re)}(\cdot)+\tilde F_{j+1}(\cdot),\]
with $|Y_j|_{\tilde h_j,\tilde\delta_{j+1}}\leq \tilde\epsilon_j^{\frac{1}{2}}$, $|\tilde A_{j+1}-\tilde A_j|_{\tilde\delta_{j+1}}\leq 3\tilde M_j'\tilde\epsilon_j$,  $\tilde F_j^{(re)}\equiv 0$, and
\begin{eqnarray*}
|\tilde F_{j+1}|_{\tilde h_{j+1}, \tilde \delta_{j+1}}\leq  3\tilde M_j'e^{-2\pi \tilde N_j(\tilde h_j-\tilde h_{j+1})}\tilde\epsilon_j <\tilde\epsilon_j^{\frac{3}{2}}=\tilde\epsilon_{j+1},
\end{eqnarray*}
where $Y_j^+(\cdot)=Y_j(\cdot+\alpha)$.

\subsubsection{Verify the properties of $\tilde A_{j+1}$}

By Lemma \ref{ndh} with $l=1$, since $64m^2\tilde M_j'^2\tilde R_j(3\tilde M_j'\tilde\epsilon_j)^{\frac{1}{m}}<1 $, then $\tilde\Sigma_{j+1}(\lambda)\subseteq D(\tilde R_{j+1})$ for $\lambda\in W_{\tilde\delta_{j+1}}(\bar\Lambda^{(j+1)})$, and for any $u\in\T$
 \[|g_{j+1}(\cdot, u)-g_j(\cdot, u)|_{\tilde\delta_{j+1}}< 2^{8m+2}(m^{10}\tilde R_j^3)^{m^2}\tilde M_j'^{m+1}\tilde\epsilon_j,\]
 where $\tilde\Sigma_{j+1}(\lambda)=\Sigma(\tilde A_{j+1}(\lambda))$, and $g_{j+1}(\lambda,u)= \prod_{\sigma_{\ell_1},\sigma_{\ell_2}\in \tilde\Sigma_{j+1}(\lambda) \atop{\ell_1\neq \ell_2}}(\sigma_{\ell_1}-e^{2\pi\ii u}\sigma_{\ell_2})$.
 Then for any $u\in\T$, we can obtain that $|g_{j+1}|_{\tilde\delta_{j+1}}\leq (2\tilde R_{j+1})^{m^2}$, and for any $\lambda\in W_{\frac{\tilde\delta_{j+1}}{2}}(\bar\Lambda^{(j+1)})\cap\R$,
\begin{eqnarray*}
   \lefteqn{ \sup_{0\leq l\leq \tilde r_1'}|\frac{\partial^l}{\partial\lambda^l}g_{j+1}(\lambda,u)| }\\
    &\geq & \sup_{0\leq l\leq \tilde r_1'}|\frac{\partial^l}{\partial\lambda^l}g_j(\lambda,u)|-\frac{\tilde r_1'!}{(\tilde\delta_{j+1}/2)^{\tilde r_1'}}2^{8m+2 } (m^{10}\tilde R_j^3)^{m^2}\tilde M_j'^{m+1}\tilde\epsilon_j>\tilde c_{j+1}'.
\end{eqnarray*}
 Thus, $\tilde\Sigma_{j+1}$ is $((2\tilde R_{j+1})^{m^2}, \tilde\delta_{j+1}, \tilde c_{j+1}', \tilde r_1')$-transvers on $\bar\Lambda^{(j+1)}\in \mathcal C(\Lambda^{(j+1)})$. Furthermore, for any $\lambda\in\bar\Lambda^{(j+1)}$, we also have 
 \begin{eqnarray*}
  \sup_{0\leq l\leq \tilde r_1'+1}|\frac{\partial^l}{\partial\lambda^l}g_{j+1}(\lambda,u)| &\leq& \sup_{0\leq l\leq \tilde r_1'+1}|\frac{\partial^l}{\partial\lambda^l}g_j(\lambda,u)|+\frac{(\tilde  r_1'+1)!}{\tilde\delta_{j+1}^{\tilde r_1'+1}}2^{8m+2} (m^{10}\tilde R_j^3)^{m^2}\tilde M_j'^{m+1}\tilde\epsilon_j\\
 & <& \tilde C_{j+1}.
\end{eqnarray*}
 We finish the proof. \qed



\section{Proof of main results}

In this section, we give the proof of main results.  Before giving the detailed proof,  we  first state an auxiliary lemma, which will be used to give the useful estimates in  KAM step.  Recall that $b\geq (120m)^{8m^3}$, $\kappa\geq m^{2m^2+10}$ are constants only depending on $m$.
\begin{Lemma}\label{pest1}
Suppose $\tilde{M}_{1}>5$, $0<\epsilon_{1}<1$, and $\Theta\subseteq\N$ satisfies $\#(\Theta\cap[s_{k},s_{k+1}))\leq m$.
Then for $n\geq 2$ we  define $R_n$, $\tilde{M}_{n}$, $\delta_{n}$ and $c_{n}$ inductively:
\[  R_{n+1}=R_{n }+bR_{n}^2\tilde M_{n}^4\epsilon_{n}^{\frac{1}{m}};\]
\begin{enumerate}[a)]
\item For $n=s_{k}-1$, we let
\begin{align*}
r_{n+1}&=m^{2}r_{n},\ \ \nu_{n+1}=u_{s_k},\\
\delta_{n+1}&= {b}^{-1}R_n^{-6m}(\nu_{n+1}\tilde M_n^{-1})^{\kappa}\delta_{n}, \\
\tilde{M}_{n+1}&=b(\nu_{n+1}^{-1}\tilde{M}_{n})^{\kappa}, \\
c_{n+1}&=( {b}^{-1}R_{n}^{-1}r_{n}^{-1}\delta_{n}\nu_{n+1}\tilde M_n^{-1}c_{n})^{\kappa^{r_{n}}}-  r_{n}^{\kappa r_{n}}\epsilon_{n}.
\end{align*}
\item For $s_{k}\leq n<s_{k+1}-1$ and $n\notin \Theta$, we have:
\begin{align*}
r_{n+1}&=r_{n},\ \ \nu_{n+1}=\nu_n-b\tilde M_n^4\epsilon_n^{\frac{1}{m}},\\
\delta_{n+1}&=\delta_{n}, \ \ \tilde{M}_{n+1}=\tilde{M}_{n}+20m\tilde M_n^2\epsilon_{n}, \\
c_{n+1}&=c_{n}- {b}R_{n}^{3m^2}\tilde{M}_{n}^{m+2}(2\delta_{n}^{-1}r_{n})^{  r_{n}}\epsilon_{n}.
\end{align*}
\item For $s_{k}\leq n<s_{k+1}-1$ and $n\in \Theta$, we have:
\begin{align*}
r_{n+1}&=m^{2}r_{n},\ \  \nu_{n+1}=\nu_n-b\tilde M_n^4\epsilon_n^{\frac{1}{m}}, \\
\delta_{n+1}&= {b}^{-1}R_n^{-(2m-1)} \nu_{n }^{m}\delta_{n}, \ \ \tilde{M}_{n+1} = \tilde{M}_{n}+20m\tilde M_n^2\epsilon_n, \\
c_{n+1}&=( {b}^{-1}R_{n}^{-1}r_{n}^{-1}\delta_{n}  c_{n})^{\kappa^{r_{n}}}-( {b}R_n  \tilde{M}_{n}\delta_{n}^{-1}r_{n})^{\kappa r_{n}}\epsilon_{n},
\end{align*}
\end{enumerate}
where $b, \kappa$ are constants only depending on $m$. Then there exists $\epsilon_{*}>0$, depending on $m,R_1,\tilde{M}_{1},\delta_{1},r_{1},c_{1}, \gamma,\tau$, such that if $\epsilon_1\leq \epsilon_*$, then  the following estimations hold:
\begin{align}
R_n&<2R_1=:R,\\
\tilde{M}_{n}&\leq (\tilde{M}_{1}\frac{|\ln\epsilon_{1}|}{2\pi h_{1}})^{b_{3}\ln n}e^{b_{3}e^{4\sqrt{n}}},\\
\delta_{n}&>\delta_{1}(\tilde M_1\frac{|\ln\epsilon_{1}|}{2\pi h_{1}})^{-b_{4}\ln n}e^{-b_{4}e^{4\sqrt{n}}},\\
c_{n}&>(c_{1}\frac{2\pi h_{1}}{|\log\epsilon_{1}|})^{n^{b_{5}}}e^{-b_{6}e^{5\sqrt{n}}},
\end{align}
where $b_{3}=b_{3}(R,m,\tau,\gamma)$, $b_{4}=b_{4}(R,m,\tau,\gamma)$, $b_{5}=b_{5}(m,R,\gamma,\tau,\delta_{1},r_1,\tilde M_1)$ and $b_{6}=b_{6}(m,R,\gamma,\tau,\delta_{1},r_1,\tilde M_1)$.
\end{Lemma}
\begin{pf}
First, we observe that
\begin{equation}\label{equ-est-r}
r_{n}\leq m^{2(m+1)\xi_{n}}r_{1}\leq m^{36m}e^{8m(\ln m) \ln^{(8)}(n)}r_{1}\leq d_{3}(m) r_{1}\ln^{(6)}(n).
\end{equation}
Moreover, by the definition of $u_n$ and $\epsilon_n$ in the beginning of Section \ref{sec-4}, there exists $\epsilon_{*,0}=\epsilon_{*,0}(m,R_1,h_1,\gamma, \tau,)$ such that if $\epsilon_1<\epsilon_{*,0}$, we have for any $n\geq 1$
\begin{equation}\label{equ-est-epsilon-u}
\epsilon_n^{\frac{1}{4m}}<u_{n+1}.
\end{equation}
Furthermore, there exists $n_*=n_*(m)\in\N$ such that for any $n\geq n_*$, we have $4^n>\kappa+2$. Therefore, there exists $\epsilon_{*,1}=\epsilon_{*,1}(m, \tilde M_1, R_1,h_1, \gamma, \tau)$ such that if $\epsilon_1<\epsilon_{*,1}$, then for any $n\leq n_*$, we have $b\tilde M_n^4\epsilon_n^{\frac{1}{2m}}<1$ and $\epsilon_n^{\frac{1}{4m}}<\nu_{n+1}$. Then if $\epsilon_1<\min\{\epsilon_{*,0}, \epsilon_{*,1}\}$, we can check that 
 for any $n\geq 1$, 
\begin{equation}\label{equ-est-tilde-m}
b\tilde M_n^4\epsilon_n^{\frac{1}{2m}}<1. 
\end{equation} 
The reason is as the following: By the selection of $\epsilon_{*,1}$, we only need to prove (\ref{equ-est-tilde-m}) holds for $n\geq n_*+1$. Suppose that (\ref{equ-est-tilde-m}) holds for $n=j\geq n_*$ . Then for $n=j+1$,  if $s_k\leq j<s_{k+1}-1$ for some $s_k$, then  we can get that
 \[
 b\tilde M_{j+1}^4\epsilon_{j+1}^{\frac{1}{2m}}=b(\tilde M_j+20m\tilde M_j^2\epsilon_j)^4\epsilon_j^{\frac{4^j}{2m}}<b\tilde M_j^4\epsilon_j^\frac{1}{2m}\cdot 2^4\epsilon_j^{\frac{4^j-1}{2m}}<1.
 \]
 Otherwise, if $j=s_{k+1}-1$, then by (\ref{equ-est-epsilon-u}) and the selection of $n_*$, we can get 
 \[
  b\tilde M_{j+1}^4\epsilon_{j+1}^{\frac{1}{2m}}=b^5\nu_{s_{k+1}}^{-4\kappa}\tilde M_j^{4\kappa}\epsilon_j^{\frac{4^j}{2m}}< \epsilon_j^{-\frac{\kappa}{m}}(b\tilde M_j^4\epsilon_j^{\frac{1}{2m}})^\kappa\epsilon_j^{\frac{4^j-\kappa}{2m}}<1.
 \]
Now, for any $n\geq 1$, there exists $k\in\N$ such that $s_k\leq n<s_{k+1}$. Then by direct calculation we can get: 
\begin{eqnarray*}
\lefteqn{\ln\tilde{M}_{n} \leq \kappa \ln\tilde{M}_{s_{k}-1}+\kappa\ln u_{s_{k}}^{-1} + \ln( 2b)}\\
&\leq&d_4(m,\tau,R,\gamma)\ln (\frac{|\ln \epsilon_1|}{2\pi h_1})+d_5(m,\tau)e^{4\sqrt s_{k}}+\kappa\ln\tilde M_{s_k-1}\\
&\leq& \kappa^k\ln\tilde M_1+
\frac{\kappa^k-1}{\kappa-1}d_4\ln (\frac{|\ln \epsilon_1|}{2\pi h_1})+d_5(\kappa^{k-1}e^{4\sqrt s_1}+\kappa^{k-2}e^{4\sqrt s_2}+\cdots+e^{4\sqrt s_k}).
 \end{eqnarray*}
By the fact that $\sqrt{s_{k+1}}>e^{\sqrt{s_{k}}}$, we obtain that
\begin{eqnarray*}
\kappa^{k-1}e^{4\sqrt s_1}+\kappa^{k-2}e^{4\sqrt s_2}+\cdots+e^{4\sqrt s_k} \leq k_*e^{4\sqrt k_*}\kappa^k+2e^{4\sqrt s_k},
\end{eqnarray*}
where $k_*$ is a constant only depending on $m$. Moreover, by the fact that $\kappa^{\xi_n}\leq \kappa^{9+2\ln d_7}\ln^{(6)}( n)$, where $d_7=d_7(m)$ is a constant only depending on $m$,  we can get
\begin{eqnarray}\label{equ-lem-11-2}
\ln\tilde M_n&\leq& b_3(R,m,\tau,\gamma)(\ln^{(6)}( n) \cdot\ln (\tilde M_1 \frac{|\ln \epsilon_1|}{2\pi h_1} )+ e^{4\sqrt s_k}) \\
\nonumber&\leq & b_3(\ln^{(6)}(n)\cdot \ln(\tilde M_1\frac{|\ln\epsilon_1|}{2\pi h_1})+e^{4\sqrt n})
\end{eqnarray}
implying that
\[\tilde M_n\leq (\tilde M_1 \frac{|\ln \epsilon_1|}{2\pi h_1} )^{b_3\ln^{(6)}( n)}e^{b_3e^{4\sqrt n}}.\]
Then by (\ref{equ-est-epsilon-u}) and (\ref{equ-est-tilde-m}), we can get that for any $n\geq 1$, when $s_k\leq n<s_{k+1}$,
\begin{equation}\label{equ-est-nu-u}
u_{s_k}\geq\nu_n \geq \nu_{s_k}-2\epsilon_{s_k}^{\frac{1}{2m}}=u_{s_k}-2\epsilon_{s_k}^{\frac{1}{2m}}>\frac{u_{s_k}}{2}.
\end{equation}
Moreover, let $\epsilon_1<\epsilon_{*,3}=R^{-4m}$. Then inductively, we can obtain that for any $n\in\N$,
\[R_n\leq R_1\prod_{j=1}^{n-1}(1+\epsilon_j^{\frac{1}{4m}})<R.\]
Now we denote $\tilde b=b (2R)^{6m}$, and  by the definition of $\delta_n$, $\tilde M_n$ and (\ref{equ-lem-11-2}), (\ref{equ-est-nu-u}), we obtain that for $s_k\leq n<s_{k+1}$,
\begin{align}\label{equ-lem-11-3}
  \delta_n &> \tilde b^{-m}\tilde M_{s_k}^{-2}\delta_{s_{k}-1} >\tilde b^{-2m}\tilde M_{s_k}^{-2}\tilde M_{s_{k-1}}^{-2}\delta_{s_{k-1}-1} >\cdots
  \\
\nonumber&> \tilde b^{-\xi_n m}(\tilde M_{s_k} \tilde M_{s_{k-1}}\cdots\tilde M_{s_2})^{-2}u_1^{m^2}\delta_1\\
\nonumber&\geq b_1^{-m^2}\tilde b^{-\xi_n m}(\tilde M_1\frac{|\ln\epsilon_1|}{2\pi h_1})^{-2b_3\xi_n\ln^{(6)}(n)}e^{-2b_3(e^{4\sqrt{ s_{k}}}+e^{4\sqrt {s_{k-1}}}+\cdots +e^{4\sqrt {s_1}})}\delta_1\\
\nonumber&\geq  (\tilde M_1\frac{|\ln\epsilon_1|}{2\pi h_1})^{-b_4\ln^{(5)}(n)}e^{-b_4e^{4\sqrt s_k}}\delta_1\\
\nonumber&\geq  (\tilde M_1\frac{|\ln\epsilon_1|}{2\pi h_1})^{-b_4\ln^{(5)}(n)}e^{-b_4e^{4\sqrt n}}\delta_1.
\end{align}
%
 By the above calculation we can find $\epsilon_{*}'=\epsilon_{*}'(R, m, \tau,\gamma,\tilde{M}_{1},\delta_{1},r_{1},h_1)$ such that if $\epsilon_{1}<\epsilon_{*}'$, then for any $n\geq 1$, we have $4e^{4\sqrt n}<|\ln\epsilon_n|$ and
\begin{align}\label{pest}
 (\tilde{b}\tilde{M}_{n}\nu_{n+1}^{-1}\delta_{n}^{-1}r_{n})^{\kappa r_{n}}\epsilon_{n}<\epsilon_{n}^{\frac{1}{2}}.
\end{align}
Now we will prove the estimation of $c_{n}$. We prove a rough estimation first.
Let $\tilde c_n=(\tilde{b}^{-1}r_{n}^{-1}\delta_{n}\nu_{n+1}\tilde M_{n}^{-1}c_{n})^{\kappa^{r_{n}}}$ if $n=s_k-1$ or $s_k\leq n<s_{k+1}-1$, $n\in\Theta$; otherwise, we let  $\tilde c_n=c_n$.
Suppose $\tilde{c}_{n}>4\epsilon_{n}^{\frac{1}{4}}$. Then we have 
\[c_{n+1}\geq \tilde c_n-(\tilde{b}\tilde{M}_{n }\delta_{n}^{-1}r_{n})^{\kappa r_{n}}\epsilon_{n}>\tilde c_n-\epsilon_n^{\frac{1}{2}}>\frac{\tilde c_n}{2}.\]
Furthermore, we have
\begin{equation}\label{equ-lem-11-1}
\tilde c_{n+1}\geq  (\tilde b^{-1}r_{n+1}^{-1}\delta_{n+1}\tilde M_{n+1}^{-1}\nu_{n+2}c_{n+1})^{\kappa^{r_{n+1}}}> (2^{-1}\tilde b^{-1}r_{n+1}^{-1}\delta_{n+1}u_{n+2}\tilde M_{n+1}^{-1}\tilde c_n)^{\kappa^{r_{n+1}}}.
\end{equation}
Then by the definition of $u_n$, (\ref{equ-est-r}), (\ref{equ-lem-11-2}), and (\ref{equ-lem-11-3}), for $s_k\leq n+1<s_{k+1}-1$, we can get
\begin{align}\label{equ-lem-11-4}
|\ln \tilde c_{n+1}|&\leq \kappa^{r_{n+1}}\left(\ln(4\tilde b)+\ln r_{n+1} +\ln \delta_{n+1}^{-1} +\ln\tilde M_{n+1}+\ln u_{s_k}^{-1}+|\ln \tilde c_n|\right)\\
\nonumber&\leq d_8(R,\gamma,\tau,r_1,\delta_1,\tilde M_1)\ln^{(2)} (n) \left(\ln\frac{|\ln\epsilon_1|}{2\pi h_1}+e^{4\sqrt{s_k}}+|\ln\tilde c_n|\right),
\end{align}
and for $n=s_{k+1}-2$, we have 
\[
|\ln \tilde c_{n+1}|\leq d_8\ln^{(2)} (n) \left(\ln\frac{|\ln\epsilon_1|}{2\pi h_1}+e^{4\sqrt{n}}+|\ln\tilde c_n|\right).
\]
Then there exists $N_*=N_*(R,\gamma,\tau,r_1,\delta_1,\tilde M_1)$, such that for $n\geq N_*$ we have $d_8\ln^{(2)} (n)<\frac{4^n}{8}$, which implies 
\[\tilde c_{n+1}>4\epsilon_{n+1}^{\frac{1}{4}},\]
 for $n\geq N_*$, provided $\tilde c_n>4\epsilon_n^{\frac{1}{4}}$. For $n\leq N_*$, by the definition of $\tilde c_n$, there exists $\epsilon_*''=\epsilon_*''(R,\gamma,\tau,r_1,\delta_1,\tilde M_1, c_1)$ such that if $\epsilon_1<\epsilon_*''$ we have $\tilde c_n>4\epsilon_n^{\frac{1}{4}}$ for $n\leq N_*$. Thus if $\epsilon_1<\min\{\epsilon_{*,1}, \epsilon_{*,2},\epsilon_{*,3}, \epsilon_*', \epsilon_*''\}$, we have $\tilde c_n>4\epsilon_n^{\frac{1}{4}}$ for any $n\in\N$.

Next we will estimate $c_{n}$ more precisely.
Suppose that $s_k-1=: m_0<m_1<\cdots<m_l< s_{k+1}-1=:m_{l+1}$, where $l\leq m$, and $\Theta\cap [s_k, s_{k+1}-1)=\{m_1, m_2, \cdots, m_l\}$. Then for any $m_i<n<m_{i+1}$, by (\ref{pest}) and the fact that $\tilde c_n>4\epsilon_n^{\frac{1}{4}}$, we have
\[\tilde c_n=c_n\geq c_{n-1}-\epsilon_{n-1}^{\frac{1}{2}}\geq \cdots\geq   c_{m_i+1}-2\epsilon_{m_i+1}^{\frac{1}{2}}\geq \tilde c_{m_i}-\epsilon_{m_i}^{\frac{1}{2}}-2\epsilon_{m_i+1}^{\frac{1}{2}}>\frac{\tilde c_{m_i}}{2}.\]
Then by (\ref{equ-lem-11-4}), we obtain that
\begin{eqnarray*}
\lefteqn{|\ln \tilde c_{s_{k+1}-2}|\leq 2|\ln \tilde c_{m_l}|\leq 2d_8 \ln^{(2)} (m_l)(\ln\frac{|\ln\epsilon_1|}{2\pi h_1}+e^{4\sqrt{s_k}}+|\ln\tilde c_{m_l-1}|)}\\
&\leq& 2d_8\ln^{(2)}(s_{k+1})(\ln\frac{|\ln\epsilon_1|}{2\pi h_1}+e^{4\sqrt{s_k}}+2|\ln\tilde c_{m_{l-1}}|)\leq \cdots\\
&\leq& (2d_8\ln s_{k+1})^{m+2}(\ln\frac{|\ln\epsilon_1|}{2\pi h_1}+e^{4\sqrt{s_k}}+ |\ln\tilde c_{s_k-2}|)\leq \cdots\\
&\leq & \Big((2d_8\ln s_{k+1})^{m+2}+(2d_8\ln s_{k+1}\cdot 2d_8\ln s_k )^{m+2}+\cdots \\
&&+(2d_8\ln s_{k+1}\cdots 2d_8\ln s_2 )^{m+2}\Big)(\ln\frac{|\ln\epsilon_1|}{2\pi h_1}+ e^{4\sqrt{s_k}})\\
&&+(2d_8\ln s_{k+1}\cdots 2d_8\ln s_2)^{m+2}| \ln\tilde c_{s_1}|\\
&\leq& k(2d_8)^{m+2}(4^k \ln  s_{k+1})^{2(m+2)}(\ln\frac{|\ln\epsilon_1|}{2\pi h_1c_1}+ e^{4\sqrt{s_k}})\\
&\leq& k(8d_8)^{2(m+2)k}( 2d_2 \sqrt{ s_{k}})^{2(m+2)}(\ln\frac{|\ln\epsilon_1|}{2\pi h_1c_1}+ e^{4\sqrt{s_k}}).
\end{eqnarray*}
Then for any $s_k\leq n<s_{k+1}$, we obtain that
\begin{eqnarray*}
|\ln c_n|&\leq& 2|\ln \tilde c_{n-1}|\leq 2|\ln\tilde c_{s_{k+1}-2}|\\
&\leq& \xi_n(8d_8)^{2(m+2)\xi_n }(2d_2 \sqrt { s_{k}})^{2(m+2)}(\ln\frac{|\ln\epsilon_1|}{2\pi h_1c_1}+ e^{4\sqrt{s_k}})\\
&\leq&d_9\ln^{(6)}(n)s_k^{m+2}(\ln\frac{|\ln\epsilon_1|}{2\pi h_1c_1}+e^{4\sqrt {s_k}})\\
&\leq&   n^{b_5}\ln\frac{|\ln\epsilon_1|}{2\pi h_1c_1}+ b_6e^{5\sqrt n}.
\end{eqnarray*}
We finish the proof.
\end{pf}

\subsection{Proof of Theorem \ref{thm-1}.}

Since $A$ is analytic in $\lambda\in\Lambda$, there exists $\tilde \delta>0$ such that $A\in C_{\tilde\delta}^\omega(\Lambda, GL(m,\C))$. By (\ref{equ-g-analytic-1}) and (\ref{trans-func-1}), for any $u\in\T$, the function $g(\lambda,u)$ is analytic in $W_{\tilde\delta}(\Lambda)$. Then by the non-degeneracy of $A$ on $\Lambda$, there exists $0<\delta\leq \tilde\delta$ such that 
 $A(\lambda)$ is non-degenerate on $W_{\frac{\delta}{2}}(\Lambda)\cap \R$ with $ \frac{c}{2}$ and $r\in\N^+$. Let 
 \begin{eqnarray*}
 &&r_0=r,\ \ c_0=c/2, \ \  \delta_0=\delta, \ \ h_0=h,\\
  &&  R_0=\max\{|A|_{\delta}, |A^{-1}|_{\delta}\},\ \  M_0=(2R_0)^{m^2}, \ \ \tilde M_0=R_0, 
 \end{eqnarray*}
  and $\varepsilon:=|F|_{h,\delta}$.
By the definition of $R_0$, we obtain that for any $\lambda\in W_{\delta}(\Lambda)$, $u\in\T$, $\Sigma(A(\lambda))\subseteq D(R_0)$, and $|g(\lambda, u)|\leq (2R_0)^{m^2}$. Together with the fact  $A(\lambda)$ is non-degenerate on $W_{\frac{\delta}{2}}(\Lambda)\cap\R$, we  obtain that $\Sigma(A(\lambda))$ is $(M_0, \delta_0, c_0, r_0)$-transverse on $\Lambda$. 

Let
\begin{eqnarray*}
&& h_1=h_0,\ \   r_1=r_0, \ \ R_1=R_0, \ \ M_1=(2R_1)^{ m^2}, \\
&& \nu_1=u_1,\ \ \zeta_1=10m\nu_1, \ \
\delta_1= b^{-1}((R_0^2\tilde M_0)^{-1}\nu_1)^{3m}\delta_0, \\
&& \tilde M_1=   b (\nu_1^{-1}\tilde M_0)^{m^2(m+2)},\  c_1=(  b R_0\nu_1^{-1}\delta_0^{-1}\tilde M_0)^{-r_0m^3(m+6)}c_0.
\end{eqnarray*}
By the definition of $u_1$ in (\ref{equ-def-u-n}), there exists $\tilde\epsilon_{*,1}=\tilde\epsilon_{*,1}(m,R_0, \tau,\gamma, h_1)>0$ such that if $\varepsilon<\tilde\epsilon_{*,1}$, we have $\nu_1=u_1<\textrm{const}\cdot \tilde M_0$. 

Recall that $b\geq \acute b$, where $\acute b$ is the constant in Lemma \ref{lem-similar-transform}.
Then by Proposition \ref{Nth} with $\nu'=\nu_1$, there exists a partition $\Pi_1$ of $\Lambda$ such that for any $\tilde\Lambda\in \Pi_1$, there exists $S_0\in C^\omega_{\delta_1}(\tilde\Lambda, \mathrm{GL}(m,\C))$ such that $\forall\lambda\in W_{\delta_1}(\tilde\Lambda)$
\[S_0^{-1}(\lambda)A(\lambda)S_0(\lambda)=\textrm{diag}\{A_{1,11}(\lambda), \cdots, A_{1,l_1l_1}(\lambda)\}=:A_1(\lambda),\]
with $|A_{1,ii}|_{\delta_1}, |S_0|_{\delta_1}, |S_0^{-1}|_{\delta_1}\leq \tilde M_1$, the decomposition $\Sigma^{(1)}(\lambda)=\Sigma^{(1)}_1(\lambda)\cup\cdots\cup\Sigma^{(1)}_{l_1}(\lambda)$ is $(M_1, \delta_1, c_1, r_1)$-transvers on $\tilde\Lambda$ and $(\nu_1,\zeta_1)$-separated for $\lambda\in W_{\delta_1}(\tilde\Lambda)$, where 
\[\Sigma^{(1)}(\lambda):=\Sigma(A_1(\lambda)) \subseteq D(R_1), \ \ \Sigma^{(1)}_i(\lambda):=\Sigma(A_{1,ii}(\lambda)).
\]
 Moreover, we have
\[\#\Pi_1\leq   b(R_0^2\tilde M_0\nu_1^{-1})^{3m}\frac{|\Lambda|}{\delta_0}+1,\]
and  there exists $\lambda_0\in \tilde\Lambda$ such that 
$\Sigma^{(1)}_1(\lambda_0)\cup\cdots\cup\Sigma^{(1)}_{l_1}(\lambda_0)$ is $(8\nu_1, \zeta_1)$-separated, and for any $\lambda\in W_{\delta_1}(\tilde\Lambda)$, $\#\Sigma_i(\lambda)=\#\Sigma_i(\lambda)$, and 
\[
d_H(\Sigma^{(1)}_i(\lambda),\Sigma_i^{(1)}(\lambda_0))<\nu_1, \ \forall1\leq i\leq l_1.
\] 
Thus, we have $\mathfrak P(1)$ holds for $A_1$.

Let $Z_0=S_0$ and 
$F_1=S_0^{-1}F S_0$.
Then there exists $\tilde\epsilon_{*,2}=\tilde\epsilon_{*,2}(m,R_0, \tau,\gamma, h_1)>0$ such that if $\varepsilon<\tilde\epsilon_{*,2}$,  
\[
|F_1|_{h_1,\delta_1}\leq  \tilde M_1^2\varepsilon=    b^2((b_1\frac{|\ln \varepsilon|}{4\pi h_1})^\tau e^{b_2e^4}R_0)^{2m^2(m+2)}\varepsilon<\varepsilon^{\frac{1}{2}}=:\epsilon_1.
\]
Moreover, we have
\[
|Z_0|_{\delta_1}, |Z_0^{-1}|_{\delta_1}<\varepsilon^{-\frac{1}{2}}=\epsilon_1^{-1}, \ \ \ |Z_0|_{\delta_1} ^{\frac{1}{2}}|F_1|_{h_1,\delta_1},  |Z_0^{-1}|_{\delta_1}^{\frac{1}{2}}|F_1|_{h_1,\delta_1}<1.
\]
In addition, there exists $\tilde\epsilon_{*,3}=\tilde\epsilon_{*,3}(d,R_0, m,\tau, \gamma,h_0, \delta_0,r_0, c_0)>0$, such that if $\varepsilon<\tilde\epsilon_{*,3}$, then (\ref{equ-pro-3-con-2})-(\ref{con1}) hold for $\tilde m=n=1$. 

Furthermore, by  Lemma \ref{pest1} and the definition of $K_n, N_{n,p}$, we can get that there exists  $\tilde\epsilon_*=\tilde\epsilon_*(d,m,\gamma,\tau,R, \delta_1,h_1, r_1, c_1)$,  such that if $\varepsilon^{\frac{1}{2}}=\epsilon_1<\tilde\epsilon_*$, then for any $n\geq 1$, the inequalities (\ref{equ-pro-3-con-2})-(\ref{con1}) hold.

 Let $\epsilon_*=\min\{\epsilon_{*,1}, \epsilon_{*,2}, \epsilon_{*,3}, \tilde\epsilon_*^2\}$, and $\varepsilon<\epsilon_*$. Then we can apply Proposition \ref{pro-iterative} inductively. 
Fix $\epsilon>0, \varsigma>0$. There exists $N^*_1=N^*_1(\epsilon, d,R_0, m,\tau, \gamma,h_0, \delta_0,r_0)$ and $N_2^*=N_2^*(\varsigma )$ such that if $n\geq N_1^*$, then $\epsilon_n<\epsilon$, and if $n\geq N_2^*$, then 
\[\sum_{j=N_2^*}^\infty \frac{2}{j^2}<\frac{\varsigma}{2}.\]
 Furthermore, for the above $N_2^*$, there exists $N_3^*=N_3^*(\varsigma)>N_2^*$ such that if $n\geq N_3^*$, then $$\epsilon_n^{\frac{\varsigma}{2}}<\prod_{j=2m+2}^{2m+N_2^*+1}\epsilon_j.$$ 
 
 Now let $N=\max\{N_1^*, N_3^*\}$. Then there exists a partition $\Pi_{N-1}$ of $\Lambda$, such that for any $\tilde\Lambda\in\Pi_{N-1}$, there is a sequence of transformations $Z_{j}\in C^\omega_{h_{j+1}, \delta_{j+1}}(\T^d\times\tilde\Lambda, \mathrm{GL}(m,\C))$ with $0\leq j\leq N-1$, such that
\[Z_{j}^{-1}(\cdot+\alpha,\lambda)\cdots Z_{ 0}^{-1}(\cdot+\alpha,\lambda)(A(\lambda)+F(\cdot, \lambda))Z_{ 0}(\cdot,\lambda)\cdots Z_{ j}(\cdot,\lambda)=A_{j+1}(\lambda)+F_{j+1}( \cdot, \lambda), \]
with estimates 
\begin{eqnarray*}
&& |F_{j+1}|_{h_{j+1},\delta_{j+1}}<\epsilon_{j+1},\ \  |F_{j+1} |_{h_{j+1},\delta_{j+1}}\leq |F_{j } |_{h_{j },\delta_j},\\
&&|Z_{j}|_{h_{j+1},\delta_{j+1}}, |Z_{ j}^{-1}|_{h_{j+1},\delta_{j+1}}<\epsilon_{j+2m+2}^{-1},\\
&&|Z_{  j}|_{h_{j+1},\delta_{j+1}}|F_{j+1}|_{h_{j+1},\delta_{j+1}}^{\frac{2}{j^2}},|Z_{ j}^{-1}|_{h_{j+1},\delta_{j+1}}|F_{j+1} |_{h_{j+1},\delta_{j+1}}^{\frac{2}{j^2}}<1.
\end{eqnarray*}
Let 
\[ 
B=Z_{0}\cdots Z_{N-1},\ \ \Pi=\Pi_{N-1},\ \ \eta=\delta_{1}(\tilde M_1\frac{|\ln\epsilon_{1}|}{2\pi h_{1}})^{-b_{4}\ln N}e^{-b_{4}e^{4\sqrt{N}}}>0, 
\]
where $b_4=b_4(R,m,\gamma,\tau)$ is defined as in Lemma \ref{pest1}. Then by Lemma \ref{pest1}, we have $\delta_N\geq \eta$. Moreover, 
\[
B^{-1}(\cdot+\alpha,\lambda)(A(\lambda)+F(\cdot, \lambda))B(\cdot,\lambda)=A_{N}(\lambda)+F_N(\cdot, \lambda )=:\tilde A(\lambda)+\tilde F(\cdot, \lambda),
\]
with $|\tilde F |_{h/2, \eta}\leq \epsilon_N<\epsilon$, and
\begin{eqnarray*}
\lefteqn{|B|_{h/2,\eta}|\tilde F|^\varsigma_{h/2,\eta}
\leq |Z_{ 0}|_{h_1,\delta_1}\cdots|Z_{ N-1}|_{h_N,\delta_N}|F_N|_{h_N,\delta_N}^{\varsigma} }\\
&\leq &( \prod_{j=2m+2}^{2m+N_2^*+1}\epsilon_j^{-1})\epsilon_N^{\frac{\varsigma}{2}} (\prod_{j=N_2^*}^{N-1}|Z_{ j}|_{h_{j+1},\delta_{j+1}})|F_N|_{h_N,\delta_N}^{\frac{\varsigma}{2}}\\
&\leq & \prod_{j=N_2^*}^{N-1}\big(|Z_{ j}|_{h_{j+1},\delta_{j+1}}|F_N|_{h_N,\delta_N}^{\frac{2}{j^2}}\big)\leq \prod_{j=N_2^*}^{N-1}\big(|Z_{ j}|_{h_{j+1},\delta_{j=1}}|F_{j+1}|_{h_{j+1},\delta_{j+1}}^{\frac{2}{j^2}}\big)<1.
\end{eqnarray*}
We can also get   $|B^{-1}|_{h/2,\eta}|\tilde F |^\varsigma_{h/2,\eta}<1$ similarly, and finish the proof.
\qed

\subsection{Proof of Theorem \ref{Thm2}:}

The main idea of the proof is the following: we will first apply the stratified quantitative almost reducibility  (Proposition \ref{pro-iterative}) to make the perturbation small enough with the prescribed estimates,  then we  apply Proposition \ref{pro-iterative-lemma-full} to obtain full measure reducibility set, and finally we estimate the Hausdorff dimension of the parameter set for irreducible cocycles.

\subsubsection{Reduce the perturbation with all parameters}

Let $|F|_{h,\delta}<\epsilon_*$, where $\epsilon_*$ is defined as in Theorem \ref{thm-1}. By Proposition \ref{pro-iterative}, there exists a sequence of partitions $\{\Pi_n\}_{n\in\N}$ of $\Lambda$ such that on each $\tilde\Lambda\in\Pi_n$, there exist  $B_n\in C^\omega_{  h_{n }, \delta_{n }}(\T^d\times \tilde\Lambda, \mathrm{GL}(m,\C))$, $A_{n}\in C^\omega_{\delta_{n}}(\tilde\Lambda, \mathrm{GL}(m,\C))$ and $F_{n}\in C_{h_{n},\delta_{n}}^\omega(\T^d\times\tilde\Lambda, {gl}(m,\C))$ such that
\[B_n^{-1}(\cdot+\alpha)(A+F(\cdot))B_n(\cdot)= A_n+  F_n(\cdot),\]
with $| F_n|_{  h_n,  \delta_n}\leq  \epsilon_n$, and $|  A_n|_{ \delta_n}\leq \tilde{  M}_n$, $A_n=\textrm{diag}\{A_{n,11},\cdots, A_{n,l_nl_n}\}$, $\Sigma^{(n)}(\lambda)=\Sigma^{(n)}_1(\lambda)\cup\cdots\cup\Sigma^{(n)}_{l_n}(\lambda)$ is $(2R_n)^{m^2}, \delta_n, c_n, r_n)$-transverse on $\tilde\Lambda$, and $\Sigma^{(n)}(\lambda)\subseteq D(R_n)$ for $\lambda\in W_{\delta_n}(\tilde\Lambda)$.
Moreover, for $n=s_k$, we have $$\#\Pi_n\leq \#\Pi_{n-1}+\frac{bR_{n-1}^{6m}\tilde M_{n-1}^{3m}|\Lambda|}{\nu_n^{3m}}\max_{\tilde\Lambda\in\Pi_{n-1}}\frac{1}{\delta_{n-1}(\tilde\Lambda)},$$ and for $s_k<n<s_{k+1}$, we have
$$\#\Pi_n\leq \#\Pi_{n-1}+\frac{bR_{n-1}^{2m-1}|\Lambda| }{\nu_n^{ m}}\max_{\tilde\Lambda\in\Pi_{n-1}}\frac{1}{\delta_{n-1}(\tilde\Lambda)}.$$ Then by Lemma \ref{pest1} and the definition of $\nu_n$, we can obtain that for $s_k\leq n<s_{k+1}$,
\begin{equation}\label{equ-thm-2-1}
\#\Pi_n\leq 1+n\frac{bR^{6m}\tilde M_{s_k-1}^{3m}|\Lambda|}{(u_{s_k}/2)^{3m}}\cdot\max_{\tilde\Lambda\in\Pi_{n-1}} \frac{1}{\delta_{n-1}(\tilde\Lambda)}<\frac{\tilde b\tilde M_{n-1}^{3m}|\Lambda| n}{u_{s_{k}}^{3m}}\max_{\tilde\Lambda\in\Pi_{n-1}} \frac{1}{\delta_{n-1}(\tilde\Lambda)},
\end{equation}
where $\tilde b=b(2R)^{6m}$.

\subsubsection{Remove the set of parameters that resonant may occur}

Now for any $0<\varrho<1, \tilde\varepsilon>0$, let
$\eta_{\tilde\varepsilon}=\tilde\varepsilon^{\frac{8}{3}}.$
 Due to the fact that 
 $$ \epsilon_\m^{\frac{\varrho}{100m^4r_\m^2}}\leq \epsilon_1^{\frac{\varrho}{100m^4r_1^2d_3^2}2^{\m^2-2\m}},$$
  by (\ref{equ-est-r}),  there exists $m_{*,1}\in\N$ such that for any $\m\geq m_{*,1}$,
\begin{equation}\label{equ-thm-2-2}
 \epsilon_\m^{\frac{\varrho}{100m^4r_\m^2}}
 <\min\{ \frac{\tilde\varepsilon^{\frac{1}{3}}}{2} , \ (9R_\m)^{-(\tau+d)},\  (2|\ln\epsilon_\m|)^{-(\tau+d)}\}.
 \end{equation}
 Moreover, by (\ref{equ-thm-2-1}) and Lemma \ref{pest1}, there exists $m_{*,2}\in\N$ such that for any $\m\geq m_{*,2}$,
\begin{equation}\label{equ-thm-2-3}
\frac{\#\Pi_\m 2^{m^2r_\m+d+3}(|\Lambda|+1)}{\gamma(\pi h_\m)^{d+\tau}}\cdot \frac{(m^2r_\m+1)!(2R_\m)^{m^2}}{\delta_\m^{m^2r_\m+1} c_\m'}<\epsilon_\m^{-\frac{\varrho}{100m^4r_\m^2}},
\end{equation}
where $  c_\m'=((\frac{\delta_\m}{4m^4r_\m^2(2R_\m)^{m^2}})^{m^4r_\m}c_\m)^{ m^{2m^2r_\m+2}}$.

Let $m_*=\max\{m_{*,1}, \ m_{*,2}\}$. For 
 any $\tilde\Lambda\in \Pi_{m_*}$, we consider the properties of $A_{m_*}\in C^\omega_{\delta_{m_*}}(\tilde\Lambda, \mathrm{GL}(m,\C))$, $F_{m_*}\in C^\omega_{h_{m_*}, \delta_{m_*}}(\T^d\times\tilde\Lambda, gl(m,\C))$. By Lemma \ref{hebin}, the multiset $\Sigma^{(m_*)}(\lambda)$ is $((2  R_{m_*})^{m^2},  \delta_{m_*},   c_{m_*}',   r_{m_*}')$-transverse on $\tilde\Lambda $, where $   r_{m_*}'=m^2  r_{m_*}$, and $\Sigma^{(m_*)}(\lambda):=\Sigma(A_{m_*}(\lambda))$. Moreover, 
let 
\[
g_{m_*}(\lambda,u)=\prod_{\sigma_i,\sigma_j\in \Sigma^{(m_*)}(\lambda), \atop{i\neq j}}(\sigma_i-e^{2\pi\ii u}\sigma_j).
\]
 Because $|g_{m_*}(\lambda, u)|_{ \delta_{m_*}}\leq (2  R_{m_*})^{m^2}$, then for any $\lambda\in\tilde\Lambda$, $u\in\T^1$, we have
\[\sup_{0\leq l\leq  r_{m_*}'+1}|\frac{\partial^l g_1(\lambda,u)}{\partial\lambda^l}|\leq \frac{(  r_{m_*}'+1)!}{ \delta_{m_*}^{  r_{m_*}'+1}}(2  R_{m_*})^{m^2}=:C_{m_*}.\]
Then by (\ref{equ-thm-2-2}) and (\ref{equ-thm-2-3}), we can apply Proposition \ref{pro-iterative-lemma-full} to $A_{m_*}, F_{m_*}$ inductively and get the sets $\Lambda^{(j+1)}=\Lambda^{(j)}\backslash\mathcal R_j(\tilde\Lambda)$ for $j\in\N$, where $\Lambda^{(1)}=\tilde\Lambda$, $\mathcal R_j(\tilde\Lambda)=\cup_{i\in J_j(\tilde\Lambda)}I^{(j)}_i$, with  
\begin{eqnarray}
\label{est-length-I-j}&&|I_{i}^{(j)}|< \epsilon_{m_*}^{\frac{2}{25 r_{m_*}'^2}\cdot (\frac{3}{2})^{j-1}},\\
\label{est-number-J-j}&&\#J_j(\tilde\Lambda)\leq 2^{ r_{m_*}' +d}\tilde N_j^d(\frac{8 C_{m_*}}{ c_{m_*}'}| \Lambda^{(j)}|+\#\mathcal C (\Lambda^{(j)})),\\
 \label{est-number-connected}&& \#\mathcal C(\Lambda^{(j+1)})\leq 2\#\mathcal C(\mathcal R^{(j)}(\tilde\Lambda)).
\end{eqnarray}

Let 
\[
 \mathcal R_j=\cup_{\tilde\Lambda\in\Pi_{m_*}} R_j(\tilde\Lambda), \ \ \mathcal R=\cup_{j=1}^\infty \mathcal R_j. 
 \]
 Then
 for any $\lambda\in \tilde\Lambda\backslash\mathcal R(\tilde\Lambda)$, there exists $\tilde B_\lambda\in C^\omega_{\frac{h_{m_*}}{2}}(\T^d, \mathrm{GL}(m,\C))$ with  $|\tilde B_\lambda|_{h_{m_*}/2}<2$, such that
\[
\tilde B_{\lambda}^{-1}(\cdot+\alpha) ( A_{m_*}(\lambda) +  F_{m_*}(\cdot,\lambda))\tilde B_\lambda(\cdot)=\tilde A(\lambda),
\]
where $\tilde A(\lambda)\in \mathrm{GL}(m,\C)$ and $\|\tilde A(\lambda)-  A_{ m_*}(\lambda)\|<6\tilde M_{m_*} \epsilon_{m_*}$.  

Let
$B_\lambda=B_{m_*}(\lambda)\tilde B_\lambda$.
Then we have $B_\lambda\in C^\omega_{h/4}(\T^d,\mathrm{GL}(m,\C))$, and for any $\lambda\in \Lambda\backslash\mathcal R$,
\[
B_\lambda^{-1}(\cdot+\alpha)(A(\lambda)+F(\cdot,\lambda))B_\lambda(\cdot)=\tilde{A}(\lambda).
\]

\subsubsection{Estimate the Hausdorff dimension of the removed set}

Denote $\mathcal S\subseteq \Lambda$ as the set of $\lambda$ such that for any $\lambda\in \mathcal S$, the cocycle $(\alpha, A(\lambda)+F(\cdot, \lambda)))$ is not reducible to some $\tilde A(\lambda)\in \GL(m,\C)$ with $\tilde A(\lambda)$ having simple eigenvalues. Then we have the following observation:

\begin{Claim} \label{claim-set-S}
$\mathcal{S}\subseteq \mathcal R=\cup_{j=1}^\infty \mathcal R_j$,
\end{Claim}
\begin{pf}
Indeed,  by Lemma \ref{ndh} with $l=1$ for fixed $\lambda$, we  obtain 
\[
d_H(  \Sigma^{(m_*)}(\lambda), \tilde\Sigma(\lambda))\leq 12m^2 \tilde M_{m_*}^2 \epsilon_{m_*}^{\frac{1}{m}}<\epsilon_{m_*}^{\frac{1}{2m}},
\]
where $\tilde\Sigma(\lambda):=\Sigma(\tilde A(\lambda))$.
Moreover, by Lemma \ref{lem-est-resonant-set}, 
for any  
$\lambda\in\Lambda \backslash \mathcal R $, we have $|g_{m_*}(\lambda,0)|\geq \epsilon_{m_*}^{ \frac{1}{10r_{m_*}'}}/2$.
Then for any $\sigma_{\ell_1},\sigma_{\ell_2}\in \Sigma^{(m_*)}(\lambda)$ with $\ell_1\neq\ell_2$, we have
\[
\epsilon_{m_*}^{ \frac{1}{10r_{m_*}'}}/2\leq|g_{m_*}(\lambda,0)|=\prod_{\sigma_i, \sigma_j\in \Sigma^{(m_*)}(\lambda), \\ \atop{i\neq j} }|\sigma_i-\sigma_j|\leq (2 R_{m_*})^{m^2-m-1}|\sigma_{\ell_1}-\sigma_{\ell_2}|,
\]
which implies that
\[|\sigma_{\ell_1}-\sigma_{\ell_2}|\geq \frac{\epsilon_{m_*}^{ \frac{1}{10r_{m_*}'}}}{2(2  R_{m_*})^{m^2-m-1} } > \epsilon_{m_*}^{\frac{1}{5m^2  r_{m_*}}}.\]
Therefore, for any $\tilde\sigma_{\ell_1}, \tilde\sigma_{\ell_2}\in \Sigma(\tilde A(\lambda))$ with $\ell_1\neq\ell_2$, we have
\[|\tilde\sigma_{\ell_1}-\tilde\sigma_{\ell_2}|\geq \epsilon_{m_*}^{\frac{1}{5m^2 r_{m_*}}}-2 \epsilon_{m_*}^{\frac{1}{2m}}>0,\]
meaning that $\tilde A(\lambda)$ has simple eigenvalues, and thus $\Lambda\backslash \mathcal R\subseteq \Lambda\backslash\mathcal S$.
\end{pf}

%
%
%
%

Recall the Hausdorff dimension of a subset is defined as in Definition \ref{def-hausdorff}. Now we give Hausdorff $\varrho$-dimensional measure estimate of  $\mathcal{S}$ for any $0<\varrho<1$, and finish the whole proof.  We first estimate the number of the removed sets. 
 \begin{Claim} \label{claim-est-R}
$ \#\mathcal C(\mathcal R_j)\leq \#\Pi_{m_*}\cdot( 2^{ r_{m_*}'+d+3} C_{m_*}(|\Lambda|+1) c_{m_*}'^{-1})^j\prod_{l=1}^j\tilde N_l^d.
$
 \end{Claim}

  \begin{pf}
Recall that $\mathcal C(\mathcal R_j)$ denotes the connected component of the set $\mathcal R_j$. By (\ref{est-number-J-j}) and (\ref{est-number-connected}), we can check that 
  \begin{eqnarray*}
\#\mathcal C(\mathcal R_j)&\leq& \sum_{\tilde\Lambda  \in\Pi_{m_*}}\#\mathcal C(R_j(\tilde\Lambda))\\
&\leq& \sum_{\tilde\Lambda  (\supseteq\Lambda^{(j)})\in\Pi_{m_*}}2^{ r_{m_*}' +d}\tilde N_j^d(\frac{8 C_{m_*}}{ c_{m_*}'}| \Lambda^{(j)}|+\#\mathcal C (\Lambda^{(j)}))\\
&\leq &2^{ r_{m_*}'+d }\tilde N_j^d\left(\frac{8 C_{m_*}}{ c_{m_*}'}|\Lambda|+2  \#\mathcal C(\mathcal R_{j-1})\right)\\
&\leq &  2^{ r_{m_*}'+d+3}(|\Lambda|+1)  C_{m_*}  c_{m_*}'^{-1}\tilde N_j^d \#\mathcal C(\mathcal R_{j-1}).
\end{eqnarray*}
Consequently, we can estimate: 
\begin{eqnarray*}
\#\mathcal C(\mathcal R_j)&\leq& ( 2^{  r_{m_*}'+d+3}  C_{m_*}(|\Lambda|+1)  c_{m_*}'^{-1})^{j-1}\tilde N_j^d\cdots \tilde N_2^d\#\mathcal C(\mathcal R_1)\\
&\leq& \#\Pi_{m_*} \cdot(2^{  r_{m_*}'+d+3}  C_{m_*}(|\Lambda|+1)  c_{m_*}'^{-1})^j\prod_{l=1}^j\tilde N_l^d.
\end{eqnarray*}
  \end{pf}

Note by (\ref{equ-thm-2-2}) and (\ref{est-length-I-j}), for any 
  interval $I\in\mathcal C(\mathcal R_j)$ with $j\geq 1$,  we have
\begin{eqnarray*}
|I|\leq \epsilon_{m_*}^{\frac{2}{25  r_{m_*}'^2}\cdot(\frac{3}{2})^{j-1}}<\tilde\varepsilon^{\frac{8}{3}}=\eta_{\tilde\varepsilon}.
\end{eqnarray*}
Then by Claim \ref{claim-set-S} and Claim \ref{claim-est-R}, combining (\ref{equ-thm-2-2}), (\ref{equ-thm-2-3}) and (\ref{equ-thm-2-4}), we can obtain that
\begin{eqnarray*}
H_{\eta_{\tilde \varepsilon}}^\varrho(\mathcal S)&\leq& \sum_{j=1}^\infty\sum_{I\in\mathcal C(\mathcal R_j)}|I|^\varrho\leq \sum_{j=1}^\infty\#\mathcal C(\mathcal R_j) \epsilon_{m_*}^{\frac{2\varrho}{25  r_{m_*}'^2}\cdot(\frac{3}{2})^{j-1}} \\
 &\leq&  \sum_{j=1}^\infty(\#\Pi_{m_*}2^{\tilde r_1'+d+2}C_{m_*}(|\Lambda|+1)  c_{m_*}'^{-1})^j(\prod_{l=1}^j\tilde N_l^d)\epsilon_{m_*}^{\frac{2\varrho}{25  r_{m_*}'^2}\cdot(\frac{3}{2})^{j-1}}\\
 &<&\sum_{j=1}^\infty \epsilon_{m_*}^{-\frac{\varrho}{100  r_{m_*}'^2}j} \epsilon_{m_*}^{(2-3(\frac{3}{2})^{j-1})\frac{\varrho}{100  r_{m_*}'^2}}\epsilon_{m_*}^{\frac{2\varrho}{25  r_{m_*}'^2}\cdot(\frac{3}{2})^{j-1}}\\
 &\leq&\sum_{j=1}^\infty \epsilon_{m_*}^{\frac{3\varrho}{100 r_{m_*}'^2}\cdot(\frac{3}{2})^{j-1}}<2 \epsilon_{m_*}^{\frac{3\varrho}{100  r_{m_*}'^2}}<\tilde\varepsilon.
  \end{eqnarray*}
Since $\tilde\varepsilon>0$ is arbitrary, then $H^{\varrho}(\mathcal S)=0$. Also by the arbitrariness of $0<\varrho<1$, we obtain that $\textrm{dim}_H(\mathcal S)=0$.
The proof is finished.\qed

\subsection{Proof of  Theorem \ref{main}  and \ref{Thm1-1}:}

Our proof is based on Aubry duality. Note IDS is invariant  under Aubry-duality (c.f. \eqref{Thm1-2-1}).
Thus in order to prove Theorem \ref{main} and Theorem \ref{Thm1-1}, we only need to prove the absolutely continuity of IDS for the corresponding dual operator:

 \begin{Theorem}\label{Thm1}
Let $ \alpha\in DC_d(\gamma,\tau)$ and $W\in C^\omega(\T^d,\R)$.There exists $\varepsilon_0=\varepsilon_0(d, \gamma,\tau,\hat V_\ell,W)$, such that if $|\varepsilon|<\varepsilon_0$, then the IDS for (\ref{op-1}) is absolutely continuous for any $E\in\R$.
\end{Theorem}

The basic observation is that if $\varepsilon$ is small, then the corresponding cocycle $(\alpha, A_\varepsilon)$ of (\ref{op-1}) can be viewed as perturbation of constant cocycle 
$(\alpha, A(E))$ where 
\begin{equation}\label{equ-matrix-a}
A(E)=\begin{pmatrix} -\frac{\hat V_{\ell-1}}{\hat V_{\ell}} & \cdots & -\frac{\hat V_{1}}{\hat V_{\ell}} & \frac{E-\hat V_0 }{\hat V_{\ell}} & -\frac{\hat V_{-1}}{\hat V_{\ell}} & \cdots & -\frac{\hat V_{-\ell+1}}{\hat V_{\ell}} & -\frac{\hat V_{-\ell}}{\hat V_{\ell}} \\ 1 &&&&&&&
\\& \ddots &&&&&&
\\&& 1 &&&&&
\\&&& 1&&&&
\\&&&& 1&&&
\\&&&&& \ddots&&
\\&&&&&& 1&0
\end{pmatrix},
\end{equation}
and it is non-degenerate:
\begin{Lemma}\label{lem-4-1}
Let $\tilde\Lambda=[a,b]$. Then $A(E)$  satisfies the non-degeneracy conditions (\ref{def-non-degeneracy}) on $\tilde\Lambda$  with some $c>0, r\in\N^+$.
\end{Lemma}
\begin{pf} We first prove that for any $u\in\T$, $g(E,u)=0$ holds for at most finitely many $E\in\tilde\Lambda$. Indeed, we have 
\begin{Claim} 
For any $u\in\T$, $\sigma_{\ell_1}, \sigma_{\ell_2}\in \Sigma(A(E))=:\Sigma(E)$ with $\ell_1\neq \ell_2$, there are only finitely many $E\in\tilde\Lambda$, such that 
\[
\sigma_{\ell_1}(E)=e^{2\pi \ii u}\sigma_{\ell_2}(E).
\]
\end{Claim}
\begin{pf}
Denote
\[
f_E(z)=\textrm{det}(zI-A(E))=\prod_{\sigma_i\in\Sigma(E)}(z-\sigma_i).\]
By direct calculation, we can get that
\begin{equation}\label{equ-f-E}
\hat V_\ell f_E(z)=\sum_{j=-\ell}^\ell\hat V_j z^{j+\ell}-E z^\ell.
\end{equation}

\textbf{Case a: $u\in\T\backslash\Z$.}
On one hand, if $\sigma_{\ell_1}(E)=e^{2\pi \ii u}\sigma_{\ell_2}(E)$, then we have 
 \begin{equation}\label{equ-zero-g-u}
 \prod_{\sigma_i \in\Sigma(E) \atop{i\neq \ell_2}}(\sigma_i-e^{2\pi \ii u}\sigma_{\ell_2})=0,
 \end{equation}
  implying 
$f_E(e^{2\pi \ii u}\sigma_{\ell_2})=0$. Moreover, we have  $f_E(\sigma_{\ell_2})=0$.
Then $\sigma_{\ell_2}(E)$ is a solution of the equation 
\begin{equation}\label{equ-lem-4-1-1}
f_E(z)=f_E(e^{2\pi \ii u}z)=0.
\end{equation}

On the other hand, because det$A(E)\neq 0$,  if $z_0$ is a solution of (\ref{equ-lem-4-1-1}), then $z_0\neq 0$ , and by (\ref{equ-f-E}), it satisfies the following equation
\[\sum_{j=-\ell}^\ell (1-e^{2\pi \ii ju})\hat V_j z^{j+\ell}=0.\]
Since the above equation is polynomial in $z$ with degree $2\ell$, it only has finitely many solutions,  implying (\ref{equ-lem-4-1-1}) only has finitely many solutions. 
Since by (\ref{equ-f-E}), for $z\neq 0$ and $E\neq E'$, we have $f_E(z)\neq f_{E'}(z)$, then (\ref{equ-lem-4-1-1}) has different solutions for different $E$. Therefore, there are only finitely many $E$ such that (\ref{equ-lem-4-1-1}) has a solution. 

Therefore, there are only finitely many $E\in\tilde\Lambda$, such that 
$
\sigma_{\ell_1}(E)=e^{2\pi \ii u}\sigma_{\ell_2}(E).
$

\medskip

\textbf{Case b: $u\in\Z$.} If $\sigma_{\ell_1}(E)=e^{2\pi \ii u}\sigma_{\ell_2}(E)$, then  $\sigma_{\ell_2}(E)$ is a zero of $f_E(z)$ with multiplicity at least 2, meaning that it is a solution of the following equation
\begin{equation}\label{equ-lem-4-1-2}
f_E(z)=f_E'(z)=0.
\end{equation}
On the other hand, if $z_0$ is a solution of (\ref{equ-lem-4-1-2}), then $z_0\neq 0$ and it satisfies
\[\sum_{j=-\ell}^\ell \frac{j}{\ell}\hat V_jz^{j+\ell}=0,\]
which also has finitely many solutions. This also implies that there are only finitely many $E\in\tilde\Lambda$, such that 
$
\sigma_{\ell_1}(E)=e^{2\pi \ii u}\sigma_{\ell_2}(E)$. We finish the proof.
\end{pf}

Therefore, by the definition of $g(E,u)$, for any $u\in\T$, $g(E,u)=0$ holds for at most finitely many $E\in\tilde\Lambda$. This implies the result. The reason is the following: Otherwise, 
for $\forall n\in\N^{+}$, there exist $u_{n}\in\T$, $E_{n}\in\tilde\Lambda$ such that
$$
\max_{0\leq l\leq n}|\frac{\partial^{l}g(E_{n},u_{n})}{\partial E^{l}}|<\frac{1}{n}.
$$
Because $\tilde\Lambda\times\T$ is compact, then there exists a subsequence of $\{(E_n, u_n )\}_n$ that converges, say to $(  E_0, u_0)\in\T\times\tilde\Lambda$. 
By (\ref{equ-g-0}), (\ref{trans-func-1}), Lemma \ref{lem-g-neq-0}, and the fact that det$A(E)=\frac{\hat V_{-\ell}}{\hat V_\ell}$, we can obtain that  
$g(E, u)$ is a polynomial in $E$ and $e^{2\pi \ii u}$. Thus, $ \frac{\partial^{l} g(E, u )}{\partial E^{l}}$ is continuous in $E$ and $u$, implying that 
for any $l\geq 0$, 
\[
\frac{\partial^{l}g(E_{0},u_{0})}{\partial E^{l}}=0.
\]
Since $g(E, u_{0})$ is a polynomial in $E$, then
$g(E, u_{0})\equiv 0$ for $E\in\tilde\Lambda$, which is a contradiction.  
\end{pf}


\textbf{Proof of  Theorem \ref{Thm1}:}
Let $\Sigma_{\varepsilon,\alpha}\subseteq \Lambda$ be the spectrum of the operator (\ref{op-1}),  where $\Lambda\subseteq\R^1$ is a bounded interval. Let $F(\theta)=(F_{ij}(\theta))_{1\leq i,j\leq 2\ell}$, where $F_{1\ell}(\theta)=-\frac{\varepsilon W(\theta)}{ \hat  V_\ell}$, and other elements $F_{ij}=0$. Then one can rewrite  $(\alpha, A_\varepsilon):=(\alpha, A(E)+F(\cdot) )$, and by Lemma \ref{lem-4-1},  $A(E)$ satisfies the non-degenerate condition on $\Lambda$ with some $r\in\N^+$, $c>0$.

By Theorem \ref{Thm2}, there exists $\tilde\varepsilon_1=\tilde\varepsilon_1(\alpha,d, V, W)>0$ such that for $|\varepsilon|\leq \tilde\varepsilon_1$, there exists $\mathcal S\subseteq\Lambda$ such that for any $E\in\Lambda\backslash\mathcal S$, there exists $B_E\in C^\omega(\T^d, \mathrm{GL}(2\ell, \C))$ such that
$$
B_E^{-1}(\cdot+\alpha)(A(E)+F(\cdot))B_E(\cdot)=D_E,
$$
where $D_E=\mathrm{diag}\{\lambda_{1},\cdots,\lambda_{2\ell} \} $ with $\lambda_{i}\neq\lambda_{j}$ for $i\neq j$. For $0<\epsilon<\max\{\|B_E\|_{0}^{-4}, \|B^{-1}_E\|_{0}^{-4}\}$, we have
$$
B_E^{-1}(\cdot+\alpha)(A(E+i\epsilon)+F(\cdot))B_E(\cdot)=D_E+\acute F_E(\cdot),
$$
where $\|\acute F_E\|_{0}\leq\frac{\|B_E\|_{0} \|B_E^{-1}\|_0}{\hat V_\ell}\epsilon=:C(E)\epsilon\ (<\hat V_\ell^{-1}\epsilon^{1/2})$. 

For any $E\in \C$, the $\ell$-th iteration of the cocycle $(\alpha, A_\varepsilon(E,\cdot))$  is a symplectic cocycle \cite{HP}. Then the Lyapunov exponents of the cocycle $(\alpha, A_\varepsilon(E,\cdot))$ appear in pairs $\pm \gamma_j (j=1,\cdots, \ell)$, implying $\hat \gamma=\sum_{j=1}^\ell \gamma_j$, where $\hat \gamma$ is the fibred entropy of the corresponding cocycle. Therefore, we have  
\begin{align*}
\hat \gamma (E+i\epsilon)&=\lim_{n\rightarrow\infty}\frac{1}{n}\int_{\T^d}\ln\|\Lambda^{\ell}(D_E(I+D_E^{-1} \acute F_E))(\theta;n)\|d\theta \\
&\leq\lim_{n\rightarrow\infty}\frac{1}{n}\int_{\T^d}\ln\|\Lambda^{\ell}D_E\|^{n}d\theta+
\lim_{n\rightarrow\infty}\frac{1}{n}\int_{\T^d}\ln\|\Lambda^{\ell}(I+D_E^{-1}\acute F_E)\|^{n}d\theta\\
&\leq\hat \gamma (E)+\ell\ln\|I+D_E^{-1}\acute F_E\|_{0}\\
&\leq \hat\gamma (E)+\ell\tilde C(E)\epsilon.
\end{align*}
Thus by  Thouless formula (c.f. \eqref{thou}), we have
\begin{align*}
\notag\hat\gamma (E+i\epsilon)-\hat\gamma (E)&=\frac{1}{2}\int\ln(1+\frac{\epsilon^{2}}{(E-E')^{2}})d\hat{\mathcal N}(E')\\
\label{Tlf}&>\frac{\ln 2}{2}(\hat{\mathcal N}(E+\epsilon)-\hat{\mathcal N}(E-\epsilon)).
\end{align*}
Hence,  for any $E\in\Sigma\backslash\mathcal{S}$, $\hat {\mathcal N}(E)$ is Lipschitz continuous. So if we decompose $\hat{\mathcal N}=\hat{\mathcal N}_{ac}+\hat{\mathcal N}_{s}$, we know $\mathcal{S}$ is a support of $\hat{\mathcal N}_{s}$.  To complete the proof, we recall the H\"older continuity of  $\hat{\mathcal{N}}(E)$:

\begin{Lemma}[\cite{GYZ}]\label{lem-4-2}
For any $0<\eta<\frac{1}{2\ell}$, $\hat{\mathcal{N}}(E)$ is $(\frac{1}{2\ell}-\eta)$-H\"older continuous for $|\varepsilon|\leq\varepsilon_*$ in the sense that
\[|\hat{\mathcal{N}}(E)-\hat {\mathcal N}(E')|\leq C_* |E-E'|^{\frac{1}{2\ell}-\eta},\]
for any $E, E'\in \R$, where $\varepsilon_*=\varepsilon_*(\alpha, V, W,\eta)$ and $C_*=C_*(\alpha, V, W)$.
\end{Lemma}

Now due to $\mathcal{S}$ is a set with Hausdorff  dimension zero, then  for $\forall\zeta>0$ we can find  a cover of $\mathcal{S}$, denoted by $\{U_{i}\}_{i=1}^{\infty}$, such that
$$
\Sigma_{i=1}^\infty\textrm{diam} (U_i)^{\frac{1}{4\ell}}<C_*^{-1}\zeta.
$$
We let $|\varepsilon|\leq\varepsilon_0=\min\{\tilde\varepsilon_1, \varepsilon_*\}$, where $\varepsilon_*$ is defined as in Lemma \ref{lem-4-2}. Then by Lemma \ref{lem-4-2} with $\eta=\frac{1}{4\ell}$, we have
$$
\hat{\mathcal N}(\mathcal{S})\leq\Sigma_{i=1}^\infty\hat{\mathcal N}(U_{i})\leq C_*\Sigma_{i=1}^\infty\textrm{diam}(U_i)^{\frac{1}{4\ell}}<\zeta.
$$
Then by the arbitrariness of $\zeta>0$, we get $\hat{\mathcal N}_{s}(\mathcal{S})=0$ and the result follows.

\section*{Acknowledgements}

J. Wang was partially supported by National Key R\&D
Program of China (2021YFA100 1600), NSFC grant (11971233), the Outstanding Youth Foundation of Jiangsu Province (No. BK20200074), and Qing Lan Project of Jiangsu province.
J. You and  Q. Zhou were partially  supported by
 National Key R\&D Program of China (2020YFA0713300)  and  Nankai Zhide Foundation.   J. You was also partially supported by NSFC grant (11871286).  Q.Zhou was  supported by NSFC grant (12071232), the Science Fund for Distinguished Young Scholars of Tianjin (No. 19JCJQJC61300).


\begin{thebibliography}{99}

\bibitem{Ai}
M. Aizenman; G. Graf, \textit{Localization bounds for an electron gas. }{\it J. Phys. A: Math. Gen.}  \textbf{31} (1998), 6783.


\bibitem{Amor}
S. Hadj Amor. \textit{H\"older continuity of the rotation number for quasi-periodic co-cycles in $SL(2,\R)$}. {\it Commun. Math. Phys.} 287(2) (2009), 565-588.

\bibitem{A0}
A. Avila, \textit{Global theory of one-frequency Schr\"{o}dinger operators.} {\it Acta Math.} \textbf{215} (2015), 1-54.

\bibitem{A2}
A. Avila, \textit{The absolutely continuous spectrum of the almost Mathieu operator.} arXiv:0810.2965.

\bibitem{A1}
A. Avila, \textit{Almost reducibility and absolute continuity I.} arXiv:1006.0704.

\bibitem{A3}
A. Avila, \textit{KAM, Lyapunov exponents and the spectral dichotomy for one-frequency Schr\"odinger operators.} Preprint.

%
%

\bibitem{AD2}
A. Avila; D. Damanik, \textit{Absolute continuity of the integrated density of states for the almost Mathieu operator with non-critical coupling.} {\it Invent. Math.} \textbf{172}(2) (2008), 439-453.



\bibitem{AJ}
A. Avila; S. Jitomirskaya, \textit{The Ten Martini Problem.} {\it Ann. of Math.}
 \textbf{170} (2009), 303-342.

\bibitem{aj1} 
A. Avila; S. Jitomirskaya, \textit{Almost localization and almost reducibility.}  {\it J. Eur. Math. Soc.} {\bf12} (2010), 93-131.



\bibitem{ALSZ}
A. Avila; Y. Last; M. Shamis; Q. Zhou,  \textit{On the abominable properties of the Almost Mathieu operator with well approximated frequencies}, to appear in {\it Duke Math. J.}.
	

\bibitem{AYZ}
A. Avila; J. You; Q. Zhou, \textit{Sharp phase transitions for the almost Mathieu operator.} {\it Duke Math. J.} \textbf{166}(14) (2017), 2697-2718.

\bibitem{AYZ2}
A. Avila; J. You; Q. Zhou, \textit{Dry ten Martini problem in the non-critical case.} Preprint.

\bibitem{AOS} J.~E.~Avron; D.~Osadchy; R.~Seiler,  \textit{A topological look at the quantum Hall effect.} {\it Physics today}  (2003), 38-42. 



\bibitem{BDGL} 
I. Binder; D. Damanik; M.Goldstein; M. Lukic, \textit{Almost periodicity in time of solutions of the
KdV equation},  {\it Duke Math. J.} \textbf{167}(14) (2018), 2633-2678.



\bibitem{Bh}
R. Bhatia,   \textit{Perturbation Bounds for Matrix Eigenvalues.} {\it Society for Industrial and Applied Mathematics} (2007).


\bibitem{B2002}
J. Bourgain, 
 \textit{On the spectrum of lattice Schr\"odinger operators with
deterministic potential.II.} {\it
J. Anal. Math., Dedicated to the memory of
Tom Wolff}, \textbf{88} (2002), 221-254.



\bibitem{B0}
J. Bourgain,  \textit{Green's function estimates for lattice Schr\"odinger operators and applications.} {\it Annals of Mathematics Studies. Princeton University Press, Princeton, NJ.} \textbf{158} (2005).


\bibitem{BG}
J. Bourgain; M. Goldstein, \textit{On non-perturbative localization with quasi-periodic potential.} {\it Ann. of Math.} \textbf{152} (2000), 835-879.


\bibitem{BJ}
J. Bourgain; S. Jitomirskaya, \textit{Absolutely continuous spectrum for 1D quasi-periodic operators.} {\it Invent. Math.} \textbf{148} (2002), 453-463.

%



\bibitem{CCYZ}
A. Cai; C. Chavaudret; J. You; Q. Zhou, \textit{Sharp H\"older continuity of the Lyapunov exponent of finitely differentiable quasi-periodic cocycles.} {\it Math. Z.} \textbf{291}(3) (2019), 931-958.



\bibitem{CSZ}
H. Cao; Y. Shi; Z. Zhang,  \textit{Quantitative Green's function estimates for lattice quasi-periodic Schr\"odinger operators.} arXiv:2209.03808.



\bibitem{CSZ1}
H. Cao; Y. Shi; Z. Zhang,  \textit{Localization and regularity of the integrated density of states for Schr\"odinger operators on $\Z^d$ with $C^2$-cosine like quasi-periodic potential}. arXiv:2303.01071


\bibitem{Cha2}
C. Chavaudret, \textit{Strong almost reducibility for analytic and gevrey quasi-periodic cocycles.} {\it Bull. Soc. Math. France} \textbf{141} (2011), 47-106.

\bibitem{ChD}
V. Chulaevsky; E. Dinaburg, \textit{Methods of KAM-theory for long-range quasi-periodic operators on $\Z^\mu$. Pure Point Spectrum.}  {\it Commun. Math. Phys.} \textbf{153} (1993), 559-577.

\bibitem{DGL1}
D. Damanik; M. Goldstein; M. Lukic, \textit{The spectrum of a Schr\"odinger operator with small quasi-periodic potential is homogeneous}, {\it J. Spec. Theory} {\bf 6} (2016), 415-427.

\bibitem{dgsv}
D. Damanik; M. Goldstein; W. Schlag; M. Voda,   \textit{Homogeneity of the spectrum for quasi-periodic Schr\"odinger operators}, {\it J. Eur. Math. Soc.} \textbf{20} (2018), 3073-3111.


%

\bibitem{DS}
E. Dinaburg; Ya. G. Sinai, \textit{The one-dimentional Schr\"{o}dinger equation with a quasi-periodic potential.} {\it Funct. Anal. Appl.} \textbf{9} (1975), 279-289.

\bibitem{Eli92}
L. H. Eliasson, \textit{Floquet solutions for the 1-dimensional quasi-periodic Schr\"{o}dinger equation.} {\it Commun. Math. Phys.} \textbf{146} (1992), 447-482.

\bibitem{Eli97}
L. H. Eliasson, \textit{Discrete one-dimensional quasi-periodic Schr\"{o}dinger operators with pure point spectrum.} {\it Acta Math.} \textbf{179} (1997), 153-196.

\bibitem{Eli98}
L. H. Eliasson, \textit{Reducibility and point spectrum for linear quasi-periodic skew-products. } {\it Proceedings of the International Congress of Mathematicians.} Vol. 2, (1998), 779-787.



\bibitem{Eli02}
L. H. Eliasson, \textit{Perturbations of linear quasi-periodic system.}
In {\it Dynamical Systems and Small Divisors} (pp. 1-60). Springer, Berlin, Heidelberg (2002).

\bibitem{Eli11}
L. H. Eliasson,  \textit{Reducibility for linear quasi-periodic differential equations.}  {\it Winter School, St Etienne de Tine} (2011).

\bibitem{Eli17}
L. H. Eliasson, \textit{Almost reducibility for the quasi-periodic linear wave equation.} The conference ``In
Memory of Jean-Christophe Yoccoz" 2017.  \href{https://www.college-de-france.fr/agenda/colloque/la-memoire-de-jean-christophe-yoccoz/almost-reducibility-for-the-quasi-periodic-linear-wave-equation.}{Viedo}


\bibitem{FJN}
R. Fabbri; R. Johnson; C. N\'{u}\={n}ez,   \textit{Rotation number for non-autonomous linear Hamiltonian systems I: Basic properties} 
{\it Z. angew. Math.  Phys.} \textbf{54} (2003), 484-502.


\bibitem{FSW}
J. Fr\"{o}hlich; T. Spencer; P. Wittwer, \textit{Localization for a class of one dimensional quasi-periodic Schr\"{o}dinger operators.} {\it Commun. Math. Phys.} \textbf{132} (1990), 5-25.

\bibitem{GJZ}
L. Ge; S. Jitomirskaya; X. Zhao, \textit{Stability of the non-critical spectral properties I: arithmetic absolute continuity of the integrated density of states.}   {\it Commun. Math. Phys.} (2023).
 https://doi.org/10.1007/s00220-023-04724-7

\bibitem{GY} L. Ge; J. You.  \textit{Arithmetic version of Anderson localization via reducibility.}  {\it Geom. Funct. Anal.} {\bf30}(5), 1370-1401 (2020).


\bibitem{GYZ}
L. Ge; J. You; X. Zhao, \textit{H\"older regularity of the integrated density of states for quasi-periodic long-range operators on $\ell^2(\Z^d)$.}  {\it Commun. Math. Phys.} \textbf{392}(2) (2022), 347-376.

\bibitem{GYZo}
L. Ge; J. You; Q. Zhou, \textit{Exponential dynamical localization: Criterion and applications.}  {\it  Ann. Scient. Ec. Norm. Sup.} \textbf{56} (2023), 91-126.

\bibitem{GS0}
M. Goldstein; W. Schlag, \textit{H\"older continuity of the integrated density of states for quasi-periodic Schr\"odinger equations and averages of shifts of subharmonic functions. } {\it Ann. of Math.} 154(1) (2001), 155-203.


\bibitem{GS}
M. Goldstein; W. Schlag, \textit{Fine properties of the integrated density of states and a quantitative separation property of the Dirichlet eigenvalues.} {\it Geom. Funtc. Anal.} \textbf{18}(3) (2008), 755-869.

%

\bibitem{HP}
A. Haro; J. Puig, \textit{A Thouless formula and Aubry duality for long-range Schr\"odinger skew-products.} {\it Nonlinearity.} \textbf{26} (2013)1163-1187.

\bibitem{Ha} 
P.G. Harper,  \textit{Single band motion of conduction electrons in a uniform magnetic
field.} \textit{Proc. Phys. Soc. London A.,} \textbf{68} (1955), 874-892.

\bibitem{HeY}
H. Her; J. You, \textit{Full measure reducibility for generic one-parameter family of quasi-periodic linear systems.} {\it J. Dynam. Differential Equations} \textbf{20}(4) (2008), 831-866.


%




\bibitem{J}
S. Jitomirskaya, \textit{Metal-insulator transition for the almost Mathieu operator.} {\it Ann. of Math.} \textbf{150} (1999), 1159-1175.

\bibitem{JK}
S. Jitomirskaya; I. Kachkovskiy, \textit{$L^2$-reducibility and localization for quasiperiodic operators.} {\it Math. Res. Lett.} \textbf{23} (2016), 431-444.
%

\bibitem{JLiu}
S. Jitormiskya; W. Liu, \textit{Universal hierarchical structure of quasiperiodic eigenfunctions.}  {\it Ann. of Math.} \textbf{187} (3) (2018), 721-776.

%





\bibitem{K}
S. Kotani,  \textit{Generalized floquet theory for stationary Sch\"odinger operators in one dimension. }
{\it Chaos Solitons Fractals} \textbf{8} (1997), 1817–1854. 

  
\bibitem{K3}
R. Krikorian, \textit{R\'eductibilit\'e des syst\`emes produits-crois\'es \`a valeurs dans des groupes compacts (French) [Reducibility of compact-group-valued skew-product systems]}. {\it Ast\'erisque} \textbf{259} (1999), 1-216.


\bibitem{K2}
R. Krikorian,
\textit{R{\'e}ductibilit{\'e} presque partout des flots fibr{\'e}s quasi-p{\'e}riodiques {\`a} valeurs dans des groupes compacts.}
 {\it Ann. Sci. \'Ecole Norm. Sup. (4)} \textbf{32}(2) (1999), 187-240.



\bibitem{LYZZ}
M. Leguil; J. You; Z. Zhao; Q. Zhou, \textit{Asymptotics of spectral gaps of quasi-periodic Schr\"odinger operators.} arXiv:1712.04700.

\bibitem{MY}
S. Marmi; J. C. Yoccoz, \textit{Some open problems related to small divisors.} {\it Dynamical Systems and Small Divisors.} Springer, Berlin, Heidelberg (2002), 175-191.

\bibitem{MP}
J. Moser; J. P\"{o}schel, \textit{An extension of a result by Dinaburg and Sinai on quasi-periodic potentials.} {\it Comment. Math. Helvetici.} \textbf{59} (1984), 39-85.

\bibitem{OA}
D.~Osadchy;  J.~E.~Avron,  \textit{Hofstadter butterfly as quantum phase
diagram.} {\it J.\ Math.\ Phys.}\ \textbf{42} (2001), 5665-5671.

\bibitem{Pa}
Y. Pan, \textit{Renormalization of symplectic quasi-periodic cocycles.} In preparation. 

\bibitem{Pe} 
R. Peierls, \textit{Zur Theorie des Diamagnetismus von Leitungselektronen.} \textit{Z. Phys.} \textbf{80} (1933), 763-791.

\bibitem{Py}
A.S. Pyartli, \textit{Diophantine approximations on submanifoids of euclidian spaces.} {\it Funkt. Anal. i. Priloz.}  \textbf{3} (1969), 303-306.

%
%

\bibitem{Pu}
J. Puig, \textit{A nonperturbative Eliasson's reducibility theorem.} {\it Nonlinearity} \textbf{19} (2006), 355-376.

\bibitem{R}
A. Rauh, \textit{Degeneracy of Landau levels in chrystals.} \textit{Phys. Status Solidi B} \textbf{65} (1974), 131-135.

%

\bibitem{Sin}
Ya. G. Sinai, \textit{Anderson localization for one-dimensional difference Schr\"{o}dinger operator with quasi-periodic potential.} {\it J. Stat. Phys.} \textbf{46} (1987), 861-909.

\bibitem{sy1} 
M. Sodin; P. Yuditskii,  \textit{Almost periodic Sturm-Liouville operators with Cantor homogeneous
spectrum.}  {\it Comment. Math. Helv.} {\bf70}(4) (1995), 639-658.

\bibitem{sy2} 
M. Sodin; P. Yuditskii, \textit{Almost periodic Jacobi matrices with homogeneous spectrum, infinite dimensional Jacobi inversion, and hardy spaces of character-automorphic functions.}  {\it J. Geom. Anal.} {\bf7} (1997), 387-435.

%
%

\bibitem{you}
J. You, \textit{Quantitative almost reducibility and its applications.} {\it In
Proceedings of the International Congress of Mathematicians.} Vol. 3, (2018),  2113–2135. 


\end{thebibliography}
\end{document}